\documentclass[preprint,10pt]{elsarticle}
\usepackage[margin=1in]{geometry}
\pdfoutput=1 

\usepackage{amssymb,amsmath,amsthm,amsbsy}
\usepackage{physics}
\usepackage[pdftex,bookmarksnumbered]{hyperref}

\usepackage{amssymb,amsmath,amsthm,amsbsy}
\usepackage{systeme}
\usepackage{physics}
\usepackage{hyperref} 
\usepackage{graphicx}
\usepackage{epstopdf}
\usepackage{epsfig}
\usepackage{multicol}
\usepackage{subfigure}
\usepackage{enumitem}
\usepackage{multirow}
\usepackage{algorithm,algorithmic}
\usepackage{bm}
\usepackage{setspace}
\usepackage{color}
\usepackage[mathlines,running]{lineno}

\journal{Computer Methods in Applied Mechanics and Engineering}

\begin{document}
\makeatletter
\def\ps@pprintTitle{%
  \let\@oddhead\@empty
  \let\@evenhead\@empty
  \let\@oddfoot\@empty
  \let\@evenfoot\@oddfoot
}
\makeatother


\begin{frontmatter}

\title{Multi-fidelity surrogate modeling using long short-term memory networks}

\author[add1]{Paolo Conti\corref{cor1}}
\ead{paolo.conti@polimi.it}
\cortext[cor1]{Corresponding author.}

\author[add2]{Mengwu Guo}
\ead{m.guo@utwente.nl}

\author[add3]{Andrea Manzoni}
\ead{andrea1.manzoni@polimi.it}

\author[add4]{Jan S. Hesthaven}
\ead{Jan.Hesthaven@epfl.ch}

\address[add1]{Department of Civil Engineering, Politecnico di Milano}
\address[add2]{Department of Applied Mathematics, University of Twente}
\address[add3]{MOX -- Department of Mathematics, Politecnico di Milano}
\address[add4]{Institute of Mathematics, \'{E}cole Polytechnique F\'{e}d\'{e}rale de Lausanne}

\begin{abstract}
 When evaluating quantities of interest that depend on the solutions to differential equations, we inevitably face the trade-off between accuracy and efficiency. Especially for parametrized, time-dependent problems in engineering computations, it is often the case that acceptable computational budgets limit the availability of high-fidelity, accurate simulation data. Multi-fidelity surrogate modeling has emerged as an effective strategy to overcome this difficulty. Its key idea is to leverage many low-fidelity simulation data, less accurate but much faster to compute, to improve the approximations with limited high-fidelity data. In this work, we introduce a novel data-driven framework of multi-fidelity surrogate modeling for parametrized, time-dependent problems using long short-term memory (LSTM) networks, to enhance output predictions both for unseen parameter values and forward in time simultaneously -- a task known to be particularly challenging for data-driven models. We demonstrate the wide applicability of the proposed approaches in a variety of engineering problems with high- and low-fidelity data generated through fine versus coarse meshes, small versus large time steps, or finite element full order versus deep learning reduced-order models. Numerical results show that the proposed multi-fidelity LSTM networks not only improve single-fidelity regression significantly, but also outperform the multi-fidelity models based on feed-forward neural networks.
\end{abstract}


\end{frontmatter}

\section{Introduction}
The construction of surrogate models is of paramount importance for multi-query, real-time simulations governed by parameterized, time-dependent partial differential equations (PDEs), as the surrogate models provide an efficient approximation of both parametric variation and time evolution in output quantities of interest. On the other hand, we often encounter the situation where many sources of data -- large in number, easily accessible or fast to compute, but not perfectly accurate -- are available. These low-fidelity (LF) data are ideally well correlated to the high-fidelity (HF) quantities that we aim to evaluate, and thus are expected to provide useful information in surrogate modeling. However, the LF data cannot guarantee a satisfactory credibility, because they are often generated by less reliable computations, such as coarse numerical discretizations, simplified physical assumptions, data fitting, and
reduced-order modeling. To overcome this technical hurdle, multi-fidelity (MF) surrogate modeling methods have been developed for an effective data fusion from multiple sources. In particular, such methods are designed to detect trends from the ample LF data and enhance the predictive accuracy where the HF data are scarce. By exploiting the possibility to fast compute large amounts of LF data over the time-parameter domain of interest, the MF methods approximate the correlation between fidelity levels and hence infer the corresponding HF output quantities with limited HF data. A successful MF strategy should not only provide high-quality model predictions, but also achieve an overall efficiency improvement that is guaranteed by a substantially reduced cost of HF evaluations. As recognized, MF surrogate modeling is especially useful for real-time assessment and monitoring, in which the running time of multi-query simulations over time and parameters should be much shorter than the operation time on the real-world asset to allow decision support. 

Several MF strategies have been proposed to model the correlation among the data from different fidelity levels and found many applications in various areas of scientific computing \cite{peherstorfer2018survey}. A widely used MF surrogate modeling technique is co-kriging \cite{Hagan2000,alvarez2012kernels}, in which vector-valued Gaussian processes are used for the regression with multiple data sources. Such a non-parametric Bayesian method has become a popular tool of MF data fusion because of its good flexibility and intrinsic uncertainty quantification. However, Gaussian process regression suffers from the curse of dimensionality and thus significantly impacts the generalization performance of co-kriging in high-dimensional problems. As well known, artificial neural networks (NNs) are capable of handling high-dimensionality \cite{han2018solving,beck2019machine}. Combined with their remarkable flexibility and multi-purpose nature, all these features made the NN-based models overwhelmingly successful in computational science and engineering, especially in the newly emerging area of \textit{scientific machine learning} \cite{baker2019workshop}. For instance, deep NNs are used to solve forward and inverse problems governed by PDEs \cite{raissi2019physics,karniadakis2021physics}; convolutional NNs have promoted massive advancements in image recognition \cite{traore2018deep}; deep auto-encoders provide a successful strategy for manifold learning and model reduction \cite{champion2019sindy,Lee2020autoencoder,fresca2020comprehensive,fresca2021poddlrom}; and transformers \cite{vaswani2017attention} and recurrent NNs, especially long short-term memory (LSTM) networks \cite{hochreiter1997long, gers2000learning}, have shown their effectiveness in time series analysis, e.g., for speech recognition \cite{ graves2005bidirectional, sundermeyer2012lstm}. Importantly, NNs appear to be a promising candiate for MF surrogate modeling, because their strong expressive power for nonlinearity should be able to detect and represent the correlation among MF data sets, even with high-dimensional inputs. In particular, \cite{Meng} tested a deep NN model with weighted linear and nonlinear components in the bi-fidelity correlation; a multi-step NN model was proposed and incorporated into Monte Carlo sampling in \cite{Motamed}; 
deep NNs were embedded into the kernel functions of co-kriging in \cite{raissi2016deep}; Gaussian process latent variables were structured into a multi-layer network for MF modeling in \cite{cutajar2019deep};
a multi-step Bayesian NN model was developed for MF, physics-informed deep learning \cite{meng2020multi}; and LF information, such as initial guess of PDE solutions, was equipped to improve the training of physics-informed NNs \cite{demo2021extended}.
More recently, MF data fusion has been explored with transfer learning \cite{song2021transfer}, deep operator networks \cite{howard2022multifidelity,lu2022multifidelity}, convolutional auto-encoders \cite{partin2022multifidelity}, and reinforcement learning \cite{khairy2022multifidelity}.
Several multi-layer, multi-fidelity NN models were introduced in \cite{guo2022multi}, inspired by the modeling assumptions of co-kriging combined with the interlink between multi-layer NNs and Gaussian processes \cite{lee2017deep,guo2021brief}. In addition, multi-fidelity NN surrogate models have been implemented in several engineering problems, such as structural health monitoring \cite{torzoni2021health, torzoni2023deep} and aerodynamics \cite{wingopt}.

Considering the demonstrated power of LSTM networks in time-series analysis, it is a natural and promising choice to embed LSTM units into the MF surrogate models for time-parameter-dependent problems. Though there have been a few LSTM-based MF techniques proposed towards specific applications, such as autonomous driving \cite{jain2018multi} and turbulent flow simulation \cite{geneva2020multi}, a general MF methodology that non-intrusively fuses a hierarchy of temporal, parametric data with LSTM networks is still lacking in the literature. To fill this gap, three data-driven MF surrogate models are proposed in this work using LSTM networks, generically for the approximation of both time evolution and parameter dependency in output quantities of interest. The proposed multi-fidelity NN models are intended as an extension of those in \cite{guo2022multi} to time-dependent problems. The choice of employing non-intrusive NN models ensures a remarkable flexibility when estimating quantities of interest for which we have no other information than the input-output data pairs. In general, however, due to a lack of incorporated physical model information, the predictive accuracy of data-driven, non-intrusive techniques may deteriorate dramatically outside the training data coverage. Poor generalization may especially be expected when extrapolating over time. In spite of this, our proposed models mitigate these challenges by combining MF strategies with LSTM networks, and present a good robustness in generalization performance. On one hand, the LSTM layers in our proposed models not only fit the data, but also recognize the underlying pattern of temporal evolution, thus enabling an effective prediction forward in time. On the other hand, the MF techniques leverage the dense sampling of cheap LF data to capture the time-parameter dependency towards accurate HF estimations. In this sense, the integration of LF data can be regards as a `regularization' that incorporates additional system information to the HF approximation. 

A major goal of this work is to highlight the importance of embedding LSTM networks into multi-fidelity schemes. Thus, we compare the proposed multi-fidelity LSTM models with both the single-fidelity LSTM networks that are trained with either LF or HF data, and the MF models based on feed-forward NNs (rather than LSTM networks). To verify the effectiveness, applicability, and generality of the proposed models, we exemplify them in a diverse collection of benchmark problems, exploiting LF and HF levels of different nature. We consider \emph{(i)} a nonlinear diffusion-reaction system, namely, a parameterized  one-dimensional FitzHugh-Nagumo membrane model that describes a simplified problem of action-potential propagation in an excitable cell; \emph{(ii)} the temporal evolution of drag and lift coefficients in a fluid flow past a cylinder, governed by the unsteady Navier-Stokes equations; and \emph{(iii)} a nonlinear dynamical system of the Lotka-Volterra type that describes the prey-predator interaction among three species. We consider output quantities with increased difficulty to approximate their time-parameter dependency. Throughout these numerical examples, we show that the proposed MF models are capable of predicting at unseen parameter locations and forward in time simultaneously. We address the following two cases for temporal extrapolation: \emph{(a)} the LF data cover the whole time interval while the HF data are only available in a shorter time window, thus requiring the inference of HF outputs from the LF data outside the HF data coverage (see Sect.~\ref{section: LV}), and \emph{(b)} the LF and HF data cover the same limited portion of the time interval, and the output predictions have to be carried out for new parameter values \texttt{AND} future states, at which we have no information of any fidelity at all (see Sect.~\ref{section: NS}). To the best of our knowledge, the proposed models are the first techniques that exploit both the multi-fidelity NNs' capability in parametric generalization and the LSTM networks' power in temporal pattern recognition, all in a general data-driven framework. Moreover, we utilize Bayesian optimization for hyperparameter tuning \cite{bergstra2011algorithms}, especially for the identification of optimal NN architectures.

This paper is structured as follows. In Sect.~\ref{section: NNs}, we introduce the proposed multi-fidelity LSTM network models and their main features. In Sect.~\ref{section: results_intro}, we present three numerical examples, including their governing physical models, quantities of interest, and bi-fidelity data sources. The results of these numerical tests are presented and discussed in Sects.~\ref{section: Pulse}, \ref{section: NS}, and \ref{section: LV}. {\color{black}{Extension to more than two fidelity levels is discussed in Sect.~\ref{section: Extension} and, finally, conclusions are drawn in Sect.~\ref{section: conclusions}}}. We also provide a hyperparameter summary in the Appendix.

\section{Multi-fidelity LSTM surrogate models}
\label{section: NNs}
The central task of this work is the surrogate modeling of certain time-parameter-dependent quantities of interest, i.e., the approximation of the map from the time and parameter inputs $\vb{x}=(t,\vb*{\mu})\in [t_0,T]\times \mathcal{P}=\mathcal{D}$ to the out quantities $\vb{f}$. Here $t_0$ and $T$ are the start and end instances of the time interval of interest, $\mathcal{P}$ is the parameter domain, and $\mathcal{D}$ denotes the full input domain. Towards this end, we present several types of NN architectures for the construction of MF surrogate models. These models are built upon the data of different fidelity levels:
\begin{itemize}
    \item \textbf{High-fidelity (HF) data}: These data represent the best achievable accuracy and are obtained from detailed numerical simulations, e.g., the finite element approximations with suitably refined spatio-temporal discretizations. For parametrized, nonlinear, time-dependent differential problems, collecting these data can be so demanding that it may clash with computational budget restrictions. In particular, the multi-query evaluations at a large number of parameter locations may result in extremely expensive, or even impracticable, cost.
    \item \textbf{Low-fidelity (LF) data}: These data are less accurate but more accessible, inexpensive or faster to compute. For instance, LF data can be obtained from coarser discretizations, less strict convergence criteria, reduced-order/surrogate models, and/or simplified modeling assumptions.
\end{itemize}

Multiple levels of fidelity can be identified by adjusting the trade-off between accuracy and computational cost, or by considering different sources of data. In this work, we focus on bi-fidelity problems, i.e., we aim to estimate time-parameter-dependent HF quantities of interest with very limited HF observations, while leveraging abundant LF observations. Regarding this modeling task, it is useful to introduce the following notation:
\begin{itemize}
    \item The LF parameter set is given by $\mathcal{P}_\texttt{LF}=\{\bm{\mu}^{(j)}: 1\leq j \leq N^\mu_{\texttt{LF}}\}$, where $\bm{\mu}^{(j)}\in \mathbb{R}^{p_\mu}$, $j = 1, ..., N^\mu_{\texttt{LF}}$, are the parameter instances of the LF observations, and ${p_\mu}$ is the dimension of the parameter space.
    \item We consider $N^t_\texttt{LF}$ equally spaced instants $\{t_n\}_{n=1}^{N^t_\texttt{LF}}$ over the time interval $[t_0,T_\texttt{LF}]$, where $t_0$ is the initial time and $T_\texttt{LF}$ is the time instant associated with the last LF measurement.
	\item The NN approximation of LF functions is denoted by $\vb{f}_\texttt{LF}$.
    \item The LF training set is given by $\mathcal{T}_\texttt{LF} = \{(\vb{x}_\texttt{LF}^{(i)},\vb{y}_\texttt{LF}^{(i)}):1\leq i \leq N_\texttt{LF}\}$, where $N_\texttt{LF} = N^t_\texttt{LF} N^\mu_\texttt{LF}$ is the total number of time-parameter combinations. Here $\vb{x}_\texttt{LF}^{(i)}\in \mathbb{R}^{p_{\text{in}}}$ collects each time-parameter location of the LF data inputs, $p_{in} = p_\mu + 1$, while $\vb{y}^{(i)}_\texttt{LF}\in \mathbb{R}^{p_{\text{out}}}$ consists of the corresponding LF observations of the $p_\text{out}$ quantities of interest that we aim to estimate.
	\item The notation for the HF counterparts is given by replacing the subscript $\texttt{LF}$ with $\texttt{HF}$.
\end{itemize}
Typically, one has $N^\mu_{\texttt{HF}} \ll N^\mu_{\texttt{LF}}$ based on the assumption that the availability of HF data is much limited due to their significantly higher computational cost compared to the LF data.

When dealing with sequential data, it is typically a good practice to organize the data into subsequences.
For this reason, we group training data  in batch tensors of shapes $n_\text{batch}\times K \times p_\text{in}$ and $n_\text{batch}\times K \times p_\text{out}$ for input and output data, respectively, in which $n_\text{batch}$ is the batch size and $K$ is the length of batch subsequences. Thus, for each parameter instance $\bm{\mu}$, we have a time-parameter input subsequence $\{\vb{x}_n\}_{n=1}^{K} = \{(t_n,\bm{\mu})\}_{n=1}^{K}$ and a corresponding output subsequence $\{\vb{y}_n\}_{n=1}^{K}$ containing the quantities of interest. Among the commonly used techniques for  sequential data processing, recurrent NNs  --  in particular LSTM networks  --  represent the state-of-the-art in a wide range of applications. In this work, we employ LSTM cells as the main block to construct multi-fidelity NN surrogate models. We briefly present the structure and main features of an LSTM unit in the following subsection.

\subsection{LSTM unit}
\begin{figure}[t]
    \centering
    \includegraphics[width=0.75\linewidth]{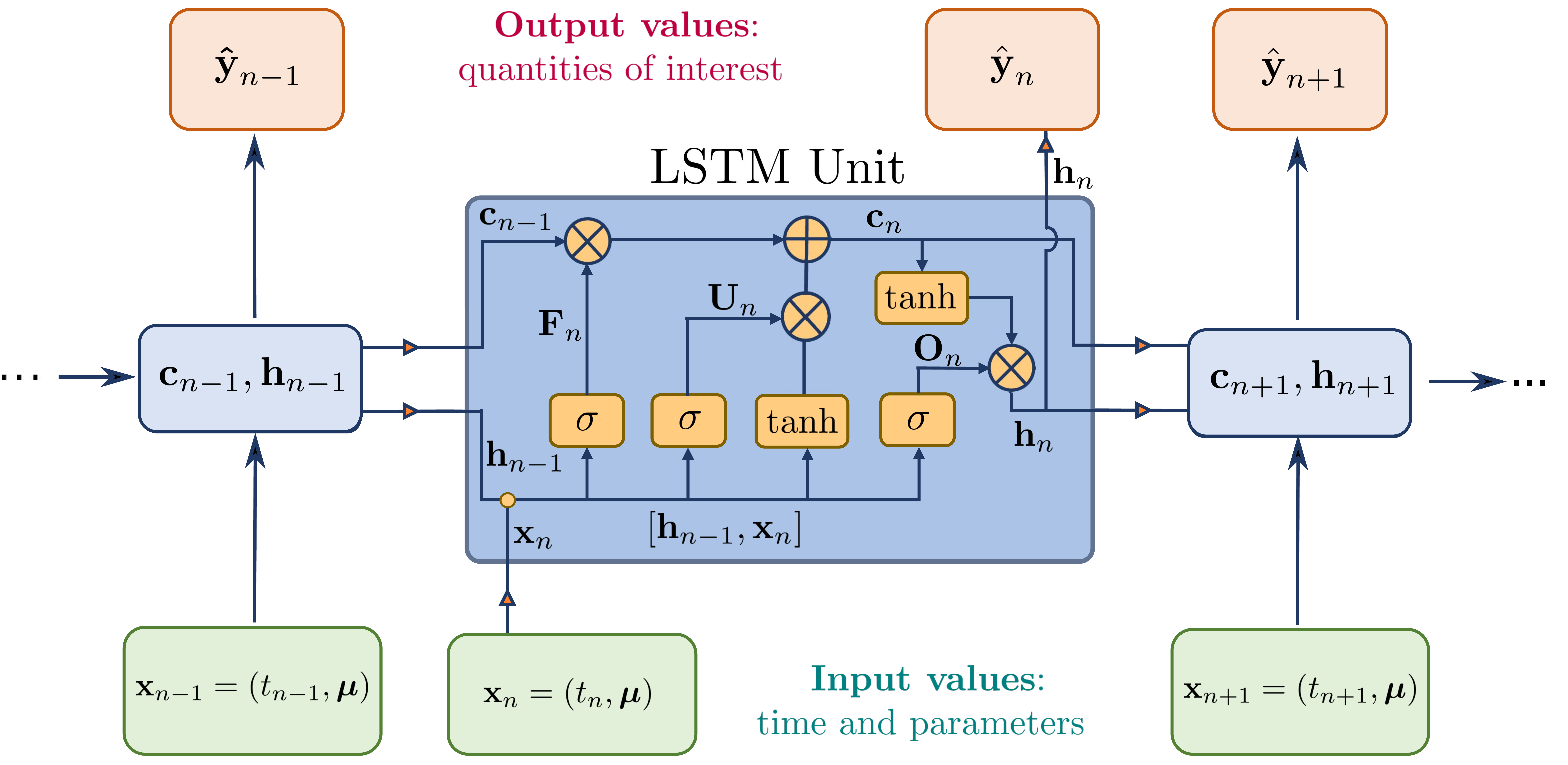}
    \caption{The {\color{black}visualization} of an LSTM unit with its input-output setting \cite{olah2015understanding,  wiki:xxx}.}
    \label{fig: LSTM}
\end{figure}

Commonly in all recurrent networks, there is a \textit{recurrent state} in an LSTM cell (see Fig. \ref{fig: LSTM}), here denoted by $\vb{h}$. The main feature of an LSTM unit is the presence of a \textit{cell state} $\vb{c}$, which accounts for long-term dependencies. At each time-step $n$, the quantities $(\vb{c}_n, \vb{h}_n)$ are computed from the previous step $(\vb{c}_{n-1}, \vb{h}_{n-1})$ together with the input values $\vb{x}_n$, namely the time instant $t_n$ and parameter values $\bm{\mu}$ in our case. These variables $(\vb{x}_n, \vb{h}_{n-1}, \vb{c}_{n-1})$ are passed to a three-fold \textit{gate mechanism} with a \textit{ forget gate} $\vb{F}_n$, an \textit{update gate} $\vb{U}_n$ and an \textit{output gate} $\vb{O}_n$, to  produce new values for the cell and recurrent states $(\vb{c}_n, \vb{h}_n)$ and to compute the output quantity $\vb{\hat{y}}_n$. Finally, $(\vb{c}_n, \vb{h}_n)$ are propagated forward through such a recurrent mechanism. The gate mechanism is the core of an LSTM unit, designed to allow a refined memory management and overcome major drawbacks of recurrent networks, such as exploding and vanishing gradients. The  \textit{gates} in such an mechanism regulate the contribution of current and previous time-step information to determine how the states $\vb{c}$ and $\vb{h}$ change over time. In particular, the cell and recurrent states are updated as follows:
\begin{equation} \label{eq: cnhn}
    \vb{c}_n = \vb{F}_n \circ \vb{c}_{n-1} + \vb{U}_n \circ \Tilde{\vb{c}}_n \,,\quad 
    \vb{h}_n = \vb{O}_n \circ \tanh{\vb{c}_n}\,, 
\end{equation}
in which the operator $\circ$ denotes the element-wise product, and $\Tilde{\vb{c}}_n$ represents the new candidate cell state to replace ${\vb{c}}_{n-1}$, computed as
\begin{equation}
    \Tilde{\vb{c}}_n = \tanh{(\vb{W}_c [\vb{h}_{n-1}, \vb{x}_{n}] + \vb{b}_c)}.
\end{equation}
The gates $\vb{F}_n$, $\vb{U}_n$, $\vb{O}_n$ are defined as
\begin{equation}
        \vb{F}_n = \sigma(\vb{W}_f [\vb{h}_{n-1}, \vb{x}_{n}] + \vb{b}_f)\,, \quad 
        \vb{U}_n = \sigma(\vb{W}_u [\vb{h}_{n-1}, \vb{x}_{n}] + \vb{b}_u)\,, \quad
        \vb{O}_n = \sigma(\vb{W}_o [\vb{h}_{n-1}, \vb{x}_{n}] + \vb{b}_o)\,.
\end{equation}

Here, $\{\vb{W}_c,\vb{W}_f,\vb{W}_u,\vb{W}_o\}$ and $\{\vb{b}_c, \vb{b}_f,\vb{b}_u,\vb{b}_o\}$ are the trainable parameters -- weights and biases -- of the LSTM unit, and $\sigma$ denotes the sigmoid activation function. Each gate takes both the concatenation of the current input parameters $\vb{x}_n$ and the previous recurrent state $\vb{h}_{n-1}$ as input variables, and provides an output vector of values between 0 and 1 that represent how much information is preserved from the pre-gate variables. The extremes 1 and 0 represent a complete preservation of information and a total discard, respectively. 
We note from \eqref{eq: cnhn} that the current cell state $\vb{c}_{n}$ is obtained by the sum of the previous cell state and the new candidate, weighted respectively by the forget and update gates, through which the contribution of past and present information is effectively managed. 
Once we have updated the cell state $\vb{c}_n$, the recurrent state $\vb{h}_n$ is updated accordingly. Different from the cell state $\vb{c}$, which is an internal variable of the LSTM recursive mechanism, the recurrent state $\vb{h}$ also serves as the output of an LSTM unit and is passed to the next unit as input. The updated $\vb{h}$ of the output layer should represent the estimation of the quantities of our interest $\vb{\hat{y}}$.
The discrepancy between the network output $\vb{\hat{y}}$ and the available data $\vb{y}$ is measured by a mean squared error (MSE) loss function,  and it is minimized to determine all the trainable parameters of the network model. For further details, we refer to \cite{olah2015understanding,hochreiter1997long,gers2000learning}.

\subsection{Multi-fidelity models}
\label{sect: mf_models}

\begin{figure}[t!]
\centering
\subfigure[2-step LSTM network model]{
	\label{subfig:2step}
	\includegraphics[width=0.7\linewidth]{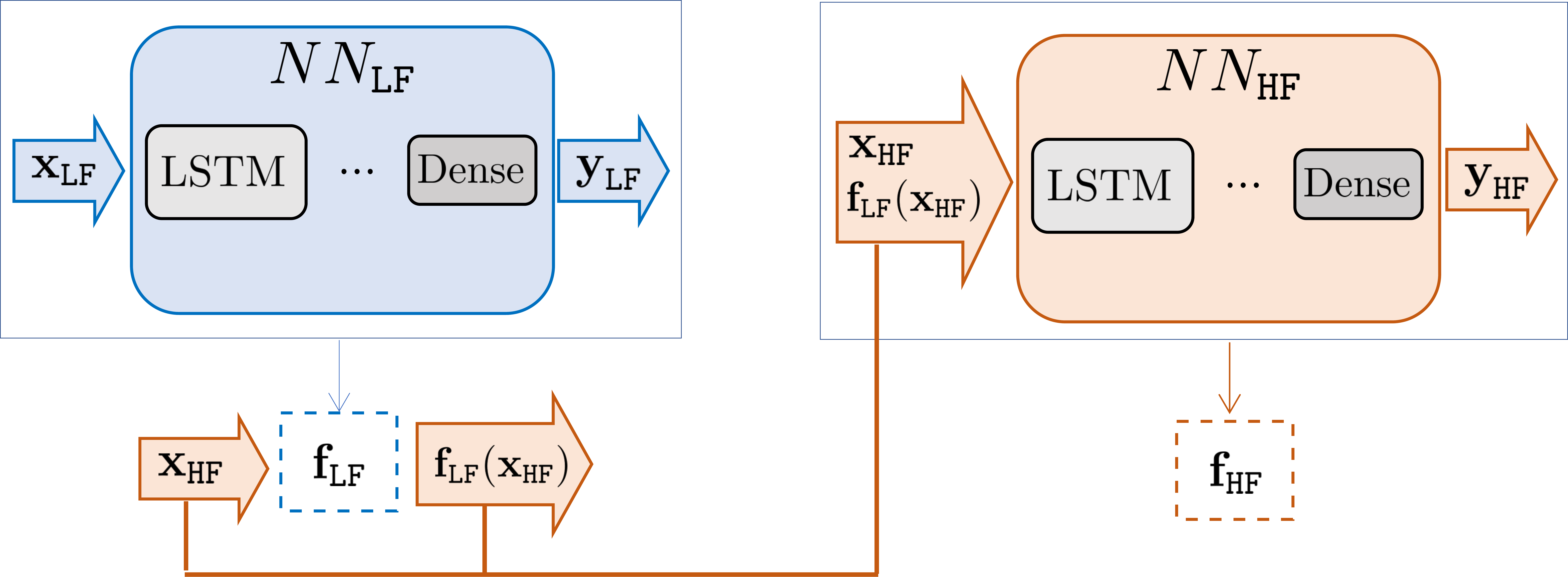} 
    } 
\centering
\subfigure[3-step LSTM network model]{
	\label{subfig:3step}
	\includegraphics[width=.95\linewidth]{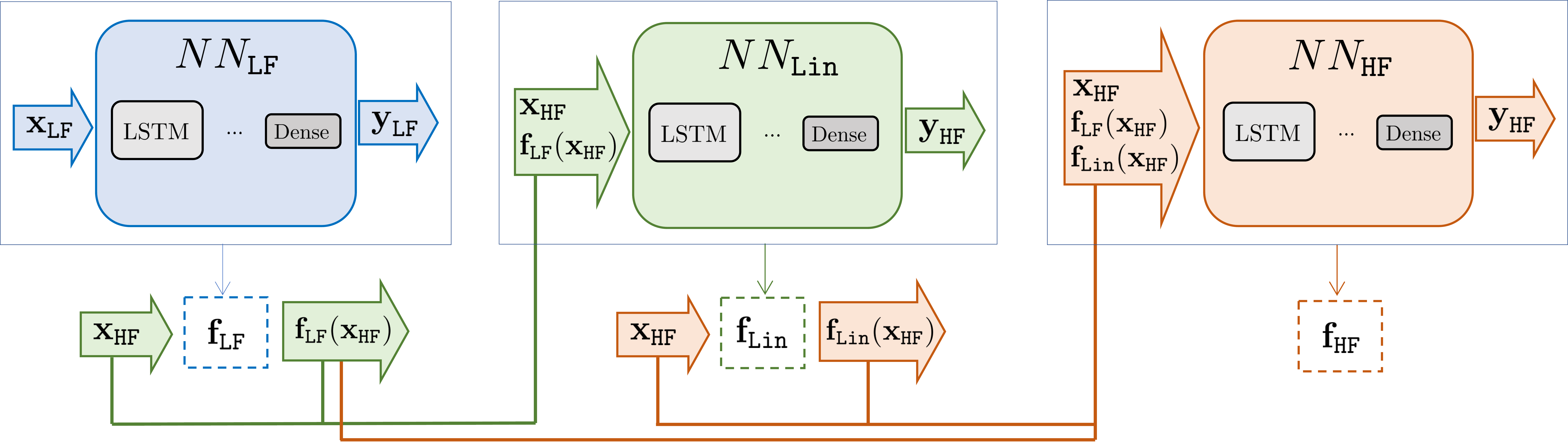} }  
\centering
\subfigure[Intermediate LSTM network model]{
	\label{subfig:Intermediate}
	\includegraphics[width=0.6\linewidth]{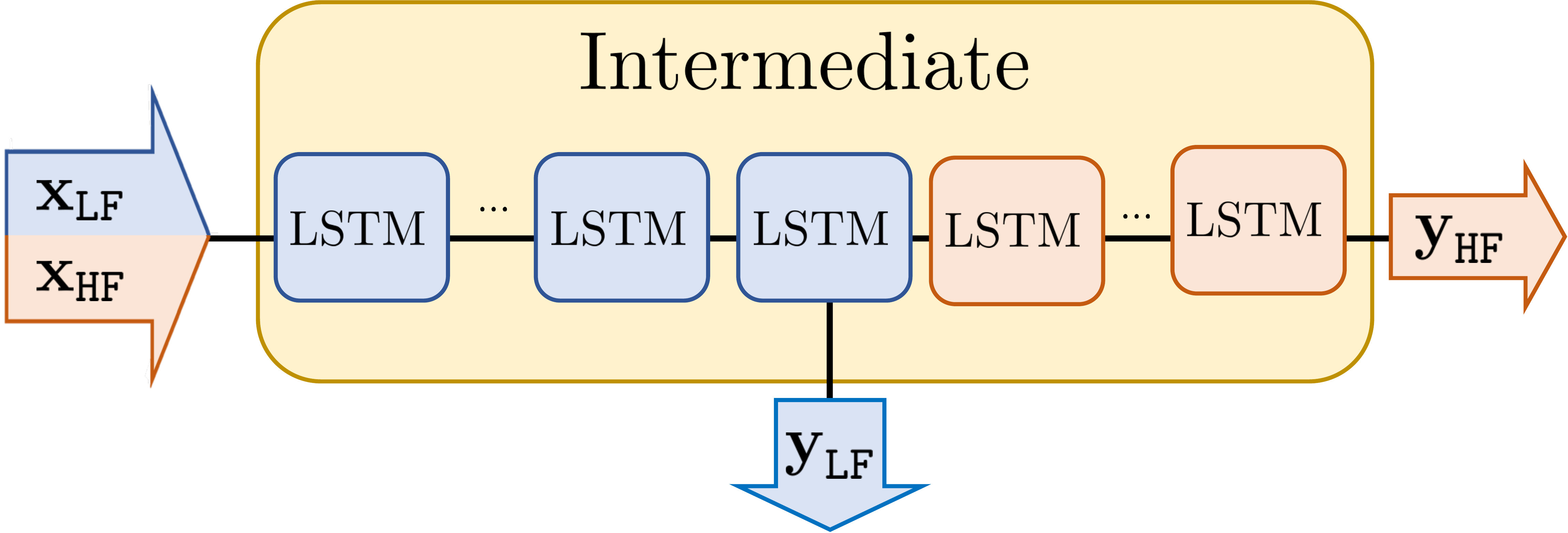} } 
\caption{Multi-fidelity LSTM network models.}
\label{fig:mfLSTM}
\end{figure}

Now we introduce a variety of multi-fidelity NN architectures for the approximation of time-parameter-dependent output quantities $\vb{f}=\vb{f}(\vb{x})=\vb{f}(t,\vb*{\mu})$. Targeting an effective treatment of time-dependency, the models described below extend the strategies proposed in \cite{guo2022multi} through an incorporation of LSTM networks. The proposed MF architectures are shown in Fig. \ref{fig:mfLSTM} and induce three different paradigms:
        
\begin{enumerate}
    \item \textbf{``2-step'' LSTM model}. The 2-step LSTM model is composed of two distinct NNs, $NN_\texttt{LF}$ and $NN_\texttt{HF}$, each consisting of a sequence of LSTM layers followed by dense layers.
    The first network $NN_\texttt{LF}$ (shown in blue in Fig.\ref{subfig:2step}) is trained on the LF data $\mathcal{T}_\texttt{LF}$ to learn the LF function $\vb{f}_\texttt{LF}$. 
    Once $NN_\texttt{LF}$ is trained, we compute the LF outputs $\vb{f}_\texttt{LF}(\vb{x}_\texttt{HF})$ for each HF training input time sequence $\vb{x}_\texttt{HF}$ through $NN_\texttt{LF}$. 
    Next, we train the second network $NN_\texttt{HF}$ (shown in red in \ref{subfig:2step}) to approximate $\vb{f}_\texttt{HF}$, using $[\vb{x}_\texttt{HF},\vb{f}_\texttt{LF}(\vb{x}_\texttt{HF})]^\text{T}$ as input data and the available HF evaluations $\vb{y}_\texttt{HF}$ as output data.
    
    \item \textbf{``3-step'' LSTM model}. This architecture is an upgrade of the 2-step LSTM model by considering an extra concatenation level, generated by a third network $NN_\texttt{Lin}$ (shown in green in Fig.\ref{subfig:3step}).  
    $NN_\texttt{Lin}$ is trained with the same input-output data as $NN_\texttt{HF}$ in the 2-step LSTM model. It does not use any nonlinear activation and thus reads the outputs as linear combinations of the inputs. Therefore, $NN_\texttt{Lin}$ captures linear correlations between the HF and LF data sets.
    Then, the third and last network $NN_\texttt{HF}$ exploits both the HF training set $\mathcal{T}_\texttt{HF}$ and the outputs of the previous NNs to approximate the HF function $\vb{f}_\texttt{HF}$, i.e., $NN_\texttt{HF}$ is trained with inputs $[\vb{x}_\texttt{HF}, \vb{f}_\texttt{LF}(\vb{x}_\texttt{HF}), \vb{f}_\texttt{Lin}(\vb{x}_\texttt{HF})]^\text{T}$ and output $\vb{y}_\texttt{HF}$.
    
    \item \textbf{``Intermediate'' LSTM model}. Different from the previous two models, the intermediate LSTM model only has a single network that simultaneously performs vector-valued learning of $[\vb{f}_\texttt{HF}(\vb{x}),\vb{f}_\texttt{LF}(\vb{x})]^\text{T}$. As shown in Fig. \ref{subfig:Intermediate}, a single input layer processes the time-parameter instances at both LF and HF levels, i.e., $\vb{x}_\texttt{LF}$ and $\vb{x}_\texttt{LF}$. Moreover, the corresponding outputs are placed at different locations: while $\vb{y}_\texttt{HF}$ is the final output layer of the network,  $\vb{y}_\texttt{LF}$ is located at an intermediate LSTM layer, for which this model is named.
    To learn $\vb{f}_\texttt{LF}$ and $\vb{f}_\texttt{HF}$ simultaneously, the corresponding MSEs, $\mathcal{L}_\texttt{LF}$ and $\mathcal{L}_\texttt{HF}$, respectively, are weighted by a hyperparameter $\alpha \in [0,1]$. Hence one can perform a single optimization minimizing the following loss function:
    \begin{equation}
        \mathcal{L} = \alpha \mathcal{L}_\texttt{HF} + (1-\alpha) \mathcal{L}_\texttt{LF}.
        \label{eq: intermediate_loss}
    \end{equation}
    in which the hyperparameter $\alpha$ regulates the contributions from the two fidelity levels.
\end{enumerate}

{\color{black}{
The mapping $(t,\vb*{\mu})\mapsto \vb{f}(t,\vb*{\mu})$ from time-parameter inputs to output quantities of interest can feature complex behaviors and thus be highly nonlinear. In this case, a large HF dataset is typically required to ensure predictive accuracy. By incorporating LF information into HF networks, the proposed MF models aim at simplifying the estimation of the HF input-output mapping so that it can be learnt 
from a small HF dataset. The HF part of the network models mainly learn the correlations between the two fidelity levels, which simplifies and accelerates the HF approximation, and prevents convergence to sub-optimal solutions. Instead of directly computing the expensive HF evaluations at new time-parameter instances, the proposed MF models seek to efficiently infer the HF outputs through the approximated fidelity correlations with cheaply obtained LF information. This not only avoids the high computational cost of HF evaluations, but also extends accurate HF predictions to the entire time-parameter domain of interest (LF data coverage), which can be considerably larger than the HF data coverage.
}}

For each model, the numbers of LSTM and dense layers and the number of neurons per layer can be determined through hyperparameters optimization (HPO), which automatically identifies the optimal network setting along with other hyperparameters, such as the learning rate, optimization algorithm and batch size \cite{hyperopt, guo2022multi}. 

We note that the intermediate LSTM model carries out simultaneous regressions on the LF and HF data, while the 2-step and 3-step LSTM models, referred to as \emph{multi-level} models, approximate the LF and HF functions sequentially through a hierarchy of separately trained networks. This affects how the contribution of each fidelity level is weighted by the model setting, and how much trust is put in the LF data for the estimation of HF quantities. While the role of the LF data is adjusted explicitly by the hyperparameter $\alpha$ in the intermediate LSTM model, this is left implicit in the multi-level models and determined in the training process. 

{\color{black}{The \emph{multi-level} models fuse multiple datasets by including additional networks. However, the computational burden of training the second- and/or third-step NNs is not heavy, because the correlation between the two fidelity levels should be uncomplicated enough to be learned with small amounts of HF data and rather simple network architectures. 
The 3-step model is expected to present advantages over the 2-step model in the case where considerable linear components exist in the correlation between HF and LF data. In addition to the $NN_\texttt{LF}$ outputs, the $NN_\texttt{Lin}$ outputs provide  additional, important features for simplifying the HF approximation by $NN_\texttt{HF}$. Conversely, if the correlation between the two fidelity levels is predominantly nonlinear, the 2-step model is sufficient.}}

All the presented models can be extended to more than two fidelity levels, either by increasing the number of ``steps'' in the multi-level approaches or by locating all the fidelity levels except for the HF in different hidden layers of the intermediate LSTM model. {\color{black}{More details are provided in Sect.~\ref{section: Extension}}}. To assess the approximation capabilities of the proposed models, however, we limit ourselves to the bi-fidelity case  and compare the proposed models with state-of-the-art regression techniques on a diversified collection of numerical benchmarks, which are introduced in the following section. 

\section{Numerical tests}
\label{section: results_intro}
In Sections \ref{section: Pulse}, \ref{section: NS}, and \ref{section: LV}, we apply the proposed multi-fidelitly LSTM models to the following numerical examples:
\begin{enumerate}[label=(\Roman*)] 
    \item  In Sect.~\ref{section: Pulse} we analyze the propagation of an electrical signal in excitable cells described by a one-dimensional nonlinear PDE-ODE (coupled) system. The problem is parameterized by a coefficient that determines the amplitude and the steepness of the action potential front. The action potential at a given spatial location is then considered as the quantity of interest. HF data are generated from highly accurate finite element approximations, while LF data come from a reduced-order model based on deep learning.
    
    \item The estimation of drag and lift coefficients for a fluid flow around a cylindrical obstacle as functions of the Reynolds number is considered in Sect.~\ref{section: NS}. Data are generated by numerical approximation of the unsteady Navier-Stokes equations, and the fidelity levels are defined by the quality of spatial and temporal discretizations.
    
    \item The last example, as discussed in Sect. \ref{section: LV}, considers  a Lotka-Volterra system that represents a three-population  prey-predator nonlinear interaction. The system is characterized by a parameter that regulates the amplitude and the frequency of the oscillatory pattern  for each population size. Data are generated by numerical time integration, and the distinction between the fidelity levels is given by the size of time steps.
\end{enumerate}
To assess the effectiveness and generality of the proposed models, the numerical examples are chosen to involve different types of governing equations and MF data generations.


In example (I), the quantity of interest -- a point-wise action potential -- presents a single wave front. In this case, we carry out MF regression with respect to both the input parameter and time using the proposed LSTM models. Their performance is compared with those of the MF approximation based only on feed-forward NNs and the LSTM regression trained solely with single-fidelity data. 

Example (II) considers the drag and lift coefficients as output quantities. In the full developed state of the fluid flow, these quantities exhibit a periodic oscillatory pattern that varies with respect to the Reynolds number. From the regression point of view, this task is much more challenging than capturing the single peak of the solution in example (I). We repeat the comparison with other regression techniques and, in addition, test the robustness of the proposed LSTM models by varying the size of HF data set. Moreover, we evaluate the simultaneous parametric interpolation (over the parameter domain) and predictive extrapolation (beyond the training time interval), assessing the generalization accuracy while narrowing the training time window.

Finally, in example (III), we extend the proposed methods to the evaluation of vector-valued output quantities.
We again test our models' extrapolation performance over time, yet with a more challenging task than example (II), as the outputs exhibit aperiodic oscillations. To track the aperiodic pattern, we let the LF solution evolve up to the final time of interest, while the HF data only cover a shorter time window. 

Such diversity of the numerical examples aims to emphasize the wide applicability of the proposed MF strategies guaranteed by their non-intrusive nature, i.e., the surrogate models are directly learned from data without requiring access to the governing physical systems or numerical solvers. Instead of being sorted by increasing complexity of the governing physical systems, the numerical examples are listed in an order with increased difficulty of  time-parameter-dependent surrogate modeling, e.g., from interpolation to extrapolation, from periodic to aperiodic behavior capture, and from scalar-valued to vector-valued function approximation.

For all the numerical examples, the comparison with the single-fidelity NNs trained on either LF or HF data set aims at highlighting the benefits of the MF modeling. Moreover, we compare the proposed LSTM architectures with the MF models solely based on feed-forward networks, such as those in \cite{guo2022multi}, to stress the critical role played by LSTM layers for enabling good generalization properties with respect to both time and parameters. All the different network models are evaluated on a test set $\mathcal{T}_\text{test}=\{(\vb{x}_i, \vb{y}_i)\}_{i=1}^{N_\text{test}}$ that covers the whole time-parameter domain of interest. The goodness of fit is measured by the following MSE:
\begin{equation}
\text{MSE}_\text{test} = \frac{1}{N_\text{test}}\sum_{i=1}^{N_\text{test}}||\vb{y}_i-\vb{\hat{y}}_i||_2^2\,,
\end{equation}
where $\{\vb{\hat{y}}_i\}_{i=1}^{N_\text{test}}$ are the predicted output values by the network models.

\section{Numerical example (I): Propagation of electrical signal}
\label{section: Pulse}
The first benchmark problem arises from computational biology and deals with the propagation of electrical potential in excitable cells, described by the following FitzHugh-Nagumo membrane model \cite{fitzhugh1961impulses, nagumo1962active}:
\begin{equation}
    \begin{aligned}
        \mu\frac{\partial{\nu}}{\partial{t}} - \mu^2 \frac{\partial^2\nu}{\partial x^2} + I_{ion}(\nu) + \omega &=0, \qquad  && x \in (0,L), \, t \in (0,T)\,,  \\
        \frac{\partial{\omega}}{\partial{t}} + (\gamma\omega-b\nu) &= 0, \qquad  &&  x \in (0,L), \, t \in (0,T)\,,  \\
        \frac{\partial{\nu}}{\partial{x}}(0,t) = -i_0(t), \quad  \frac{\partial{\nu}}{\partial{x}}(L,t) &=0, \qquad  && t \in (0,T)\,,  \\
        \nu(x,0)=0, \quad \omega(x,0)&=0, \qquad  && x \in (0,L) \,.
        \label{eq: time_sytem}
    \end{aligned}
\end{equation}
\noindent where $\nu$ is the cardiac transmembrane electrical potential (excitation variable), $\omega$ is the recovery variable which accounts for the refractoriness of heart cells, and $t$ denotes a rescaled time.
As given in \cite{pagani2018numerical}, we choose $I_{ion}(\nu) = \nu(\nu-0.1)(\nu-1)$, $T=2$, $L=1$, $\gamma = 2$, and $b=0.5$.
This problem is parametrized by $\mu \in \mathcal{P} = [0.005,0.05]$, and its solution is characterized by the traveling wave fronts whose shapes are determined by the parameter $\mu$. The model can be extended to two- or three-dimensional spatial domains to represent the propagation of electrical stimuli on thin portions of tissue, or even on realistic geometries at the organ scale. 
For such a PDE-ODE coupled system -- even more substantially, for two- or three-dimensional cases -- full order solvers may be computationally demanding, as they usually involve very small time steps to properly detect signal propagation, as well as fine spatial meshes because of the steep fronts of the action potential. 
Therefore, we take advantage of MF strategies to achieve a reasonable approximation with improved efficiency and controlled accuracy.
For the case at hand, the quantity of interest is the solution value of the electrical potential $\nu$ at $\bar{x} = 0.5$, regarded as a function of the parameter $\mu$ and time $t$. 

\subsection{Multi-fidelity setting}
\begin{figure}[t!]
	\centering
	\subfigure[LF solution]
	{\includegraphics[width=0.41\linewidth]{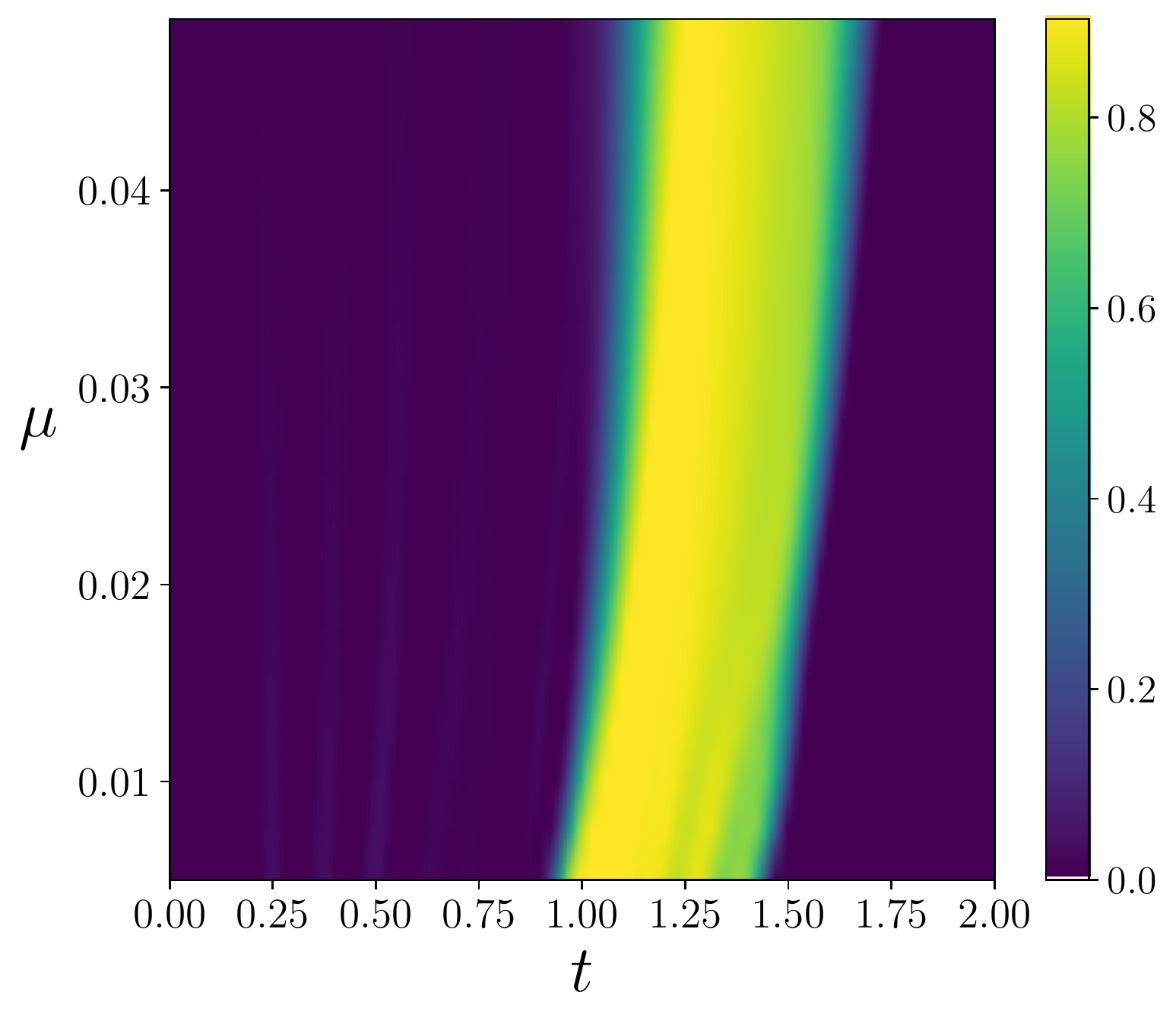}}
	\qquad
	\subfigure[HF solution]
	{\includegraphics[width=0.41\linewidth]{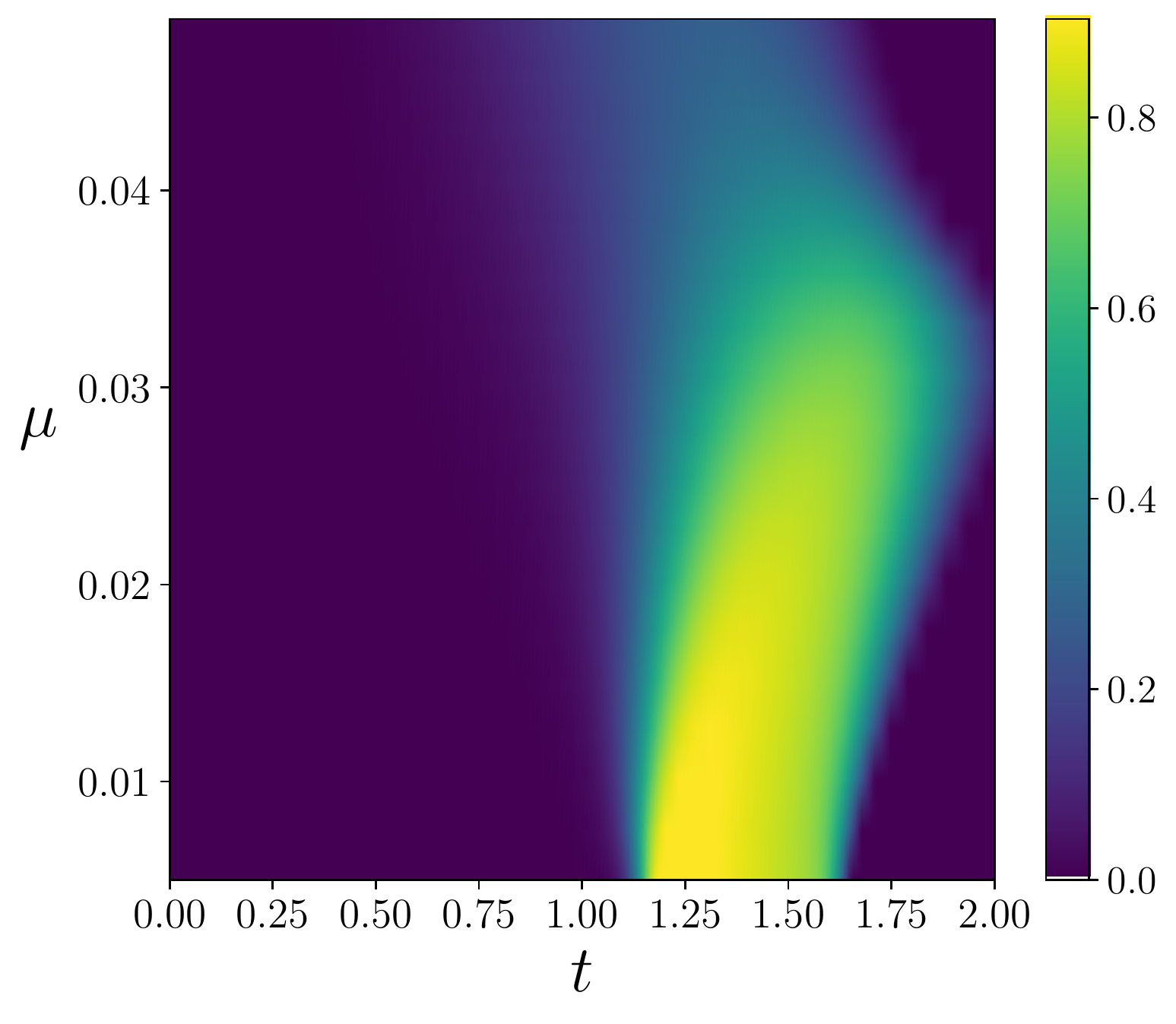}}
    \caption[PDE Low- and high-fidelity models.]{LF and HF solutions in example (I).}
    \label{fig: LF_HF_time}
\end{figure}

The HF model is obtained by a linear finite element discretization over the spatial domain $\Omega = (0,L)$ with the first-order semi-implicit scheme for time integration, and the total number of degrees of freedom is 1024. 
The LF model is constructed through a deep-learning-based reduced order modeling technique, POD-DL-ROM \cite{fresca2021poddlrom}, which efficiently constructs non-intrusive reduced-order models for nonlinear parametrized time-dependent problems starting with a prior dimensionality reduction through the proper orthogonal decomposition.
2 degrees of freedom (latent variables) are adopted in the LF reduced-order model throughout this example. We refer to \cite{pagani2018numerical} and \cite{fresca2020deep, fresca2021poddlrom} for further details about the construction of the HF and the LF models, respectively. Solutions for the quantity of interest at the two fidelity levels are shown in Fig.~\ref{fig: LF_HF_time}.

Training data are obtained by solving either the LF or the HF model for parameter-time locations sampled over $\mathcal{P}\times[0,T]$. We consider $N^\mu_\texttt{HF} = 4$ (resp. $N^\mu_\texttt{LF} = 25$) parameter values uniformly spaced over $\mathcal{P}$ for the sampling of HF (resp. LF) training data. For both fidelity levels, we choose $\Delta t=0.1$ as time step size, such that $N^t_\texttt{HF} = N^t_\texttt{LF} = 20$. We consider full-length time sequences by taking $K = 20$.
In the testing stage, we employ a test set consisting of the HF evaluations on a uniform tensor grid over $\mathcal{P}\times [0,T]$ with $N_\text{test} = 18$ parameter locations and 20 time instants.

\subsection{Results and discussions}
\begin{figure}[t!]
	\centering
	\hspace*{-3mm}
	\subfigure[LF feed-forward]
	{\includegraphics[width=0.253\linewidth]{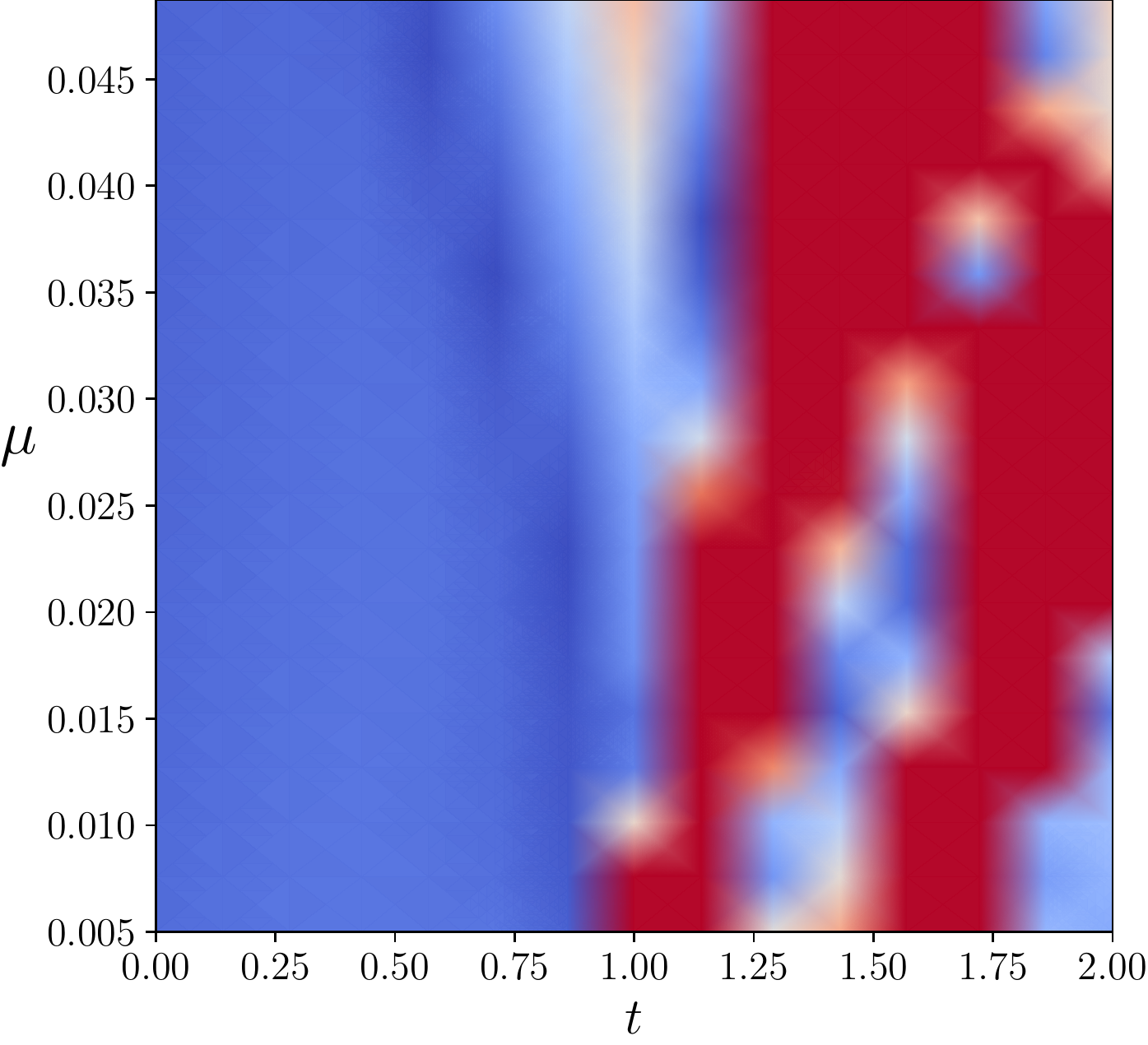}
	\label{fig: err_pulse_LF_FF}}
	\subfigure[HF feed-forward]
	{\includegraphics[width=0.229\linewidth]{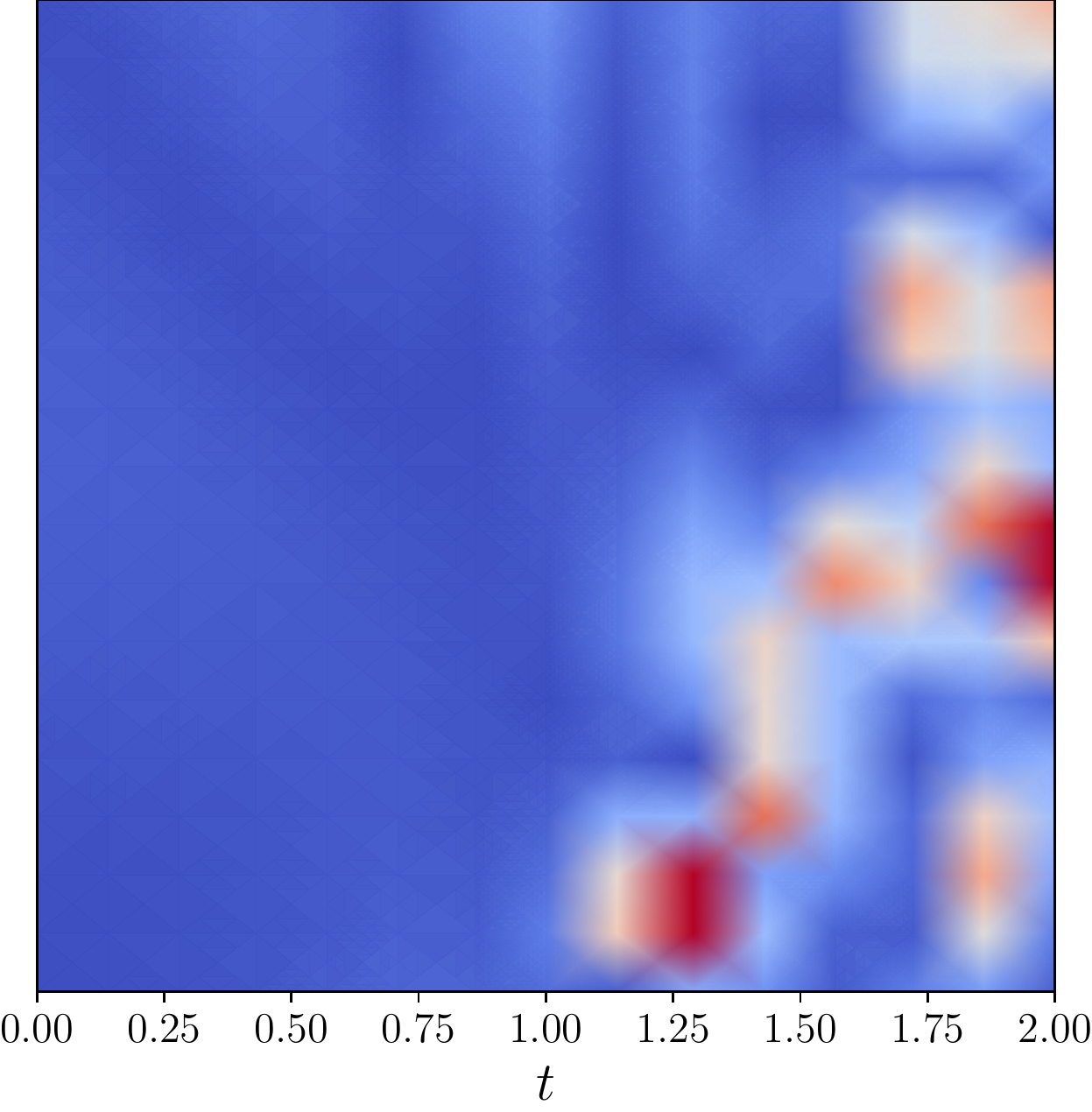}
	\label{fig: err_pulse_HF_FF}}
	\subfigure[LF LSTM]
	{\includegraphics[width=0.229\linewidth]{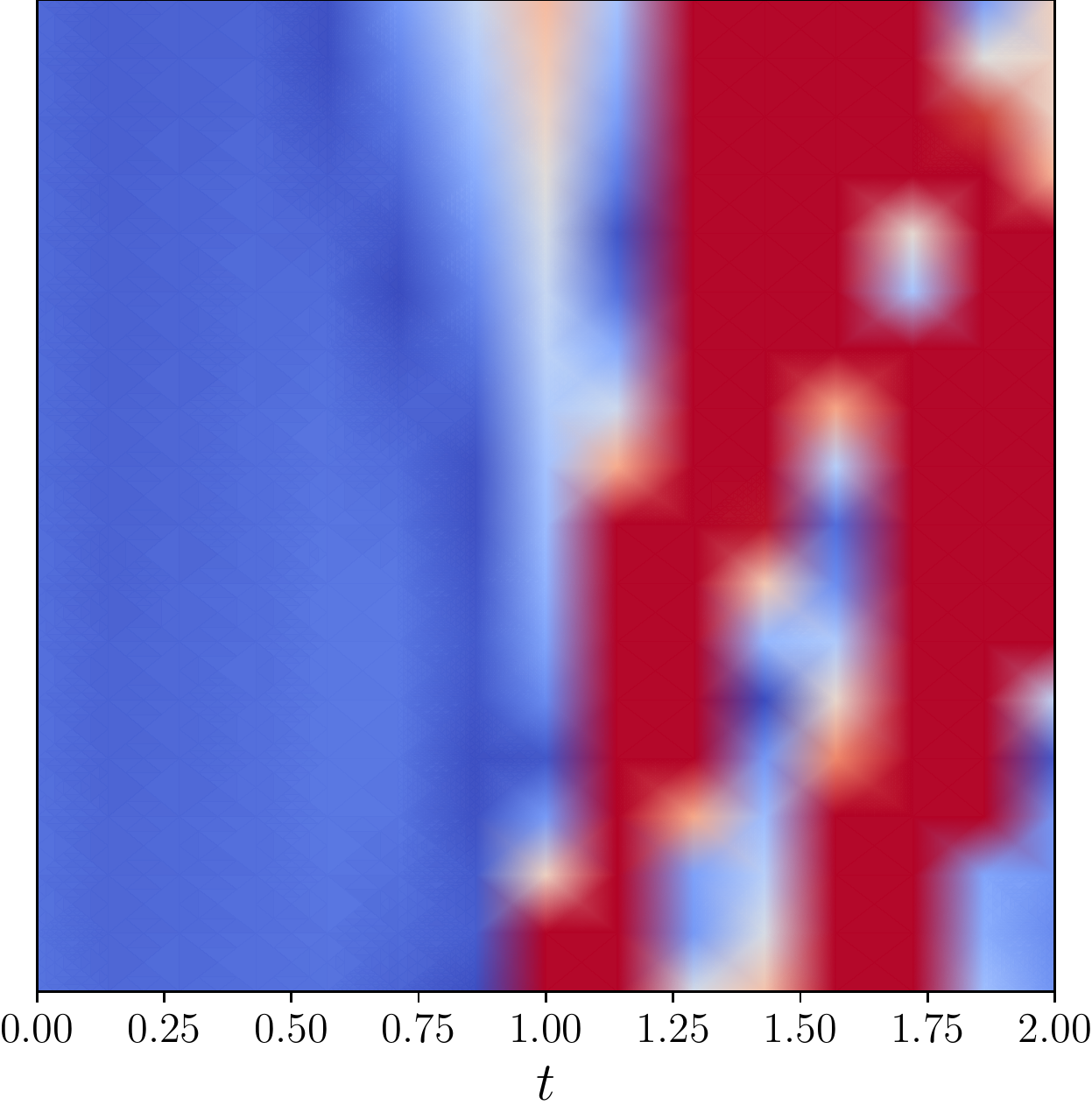}
	\label{fig: err_pulse_LF_lstm}}
	\subfigure[HF LSTM]
	{\includegraphics[width=0.249\linewidth]{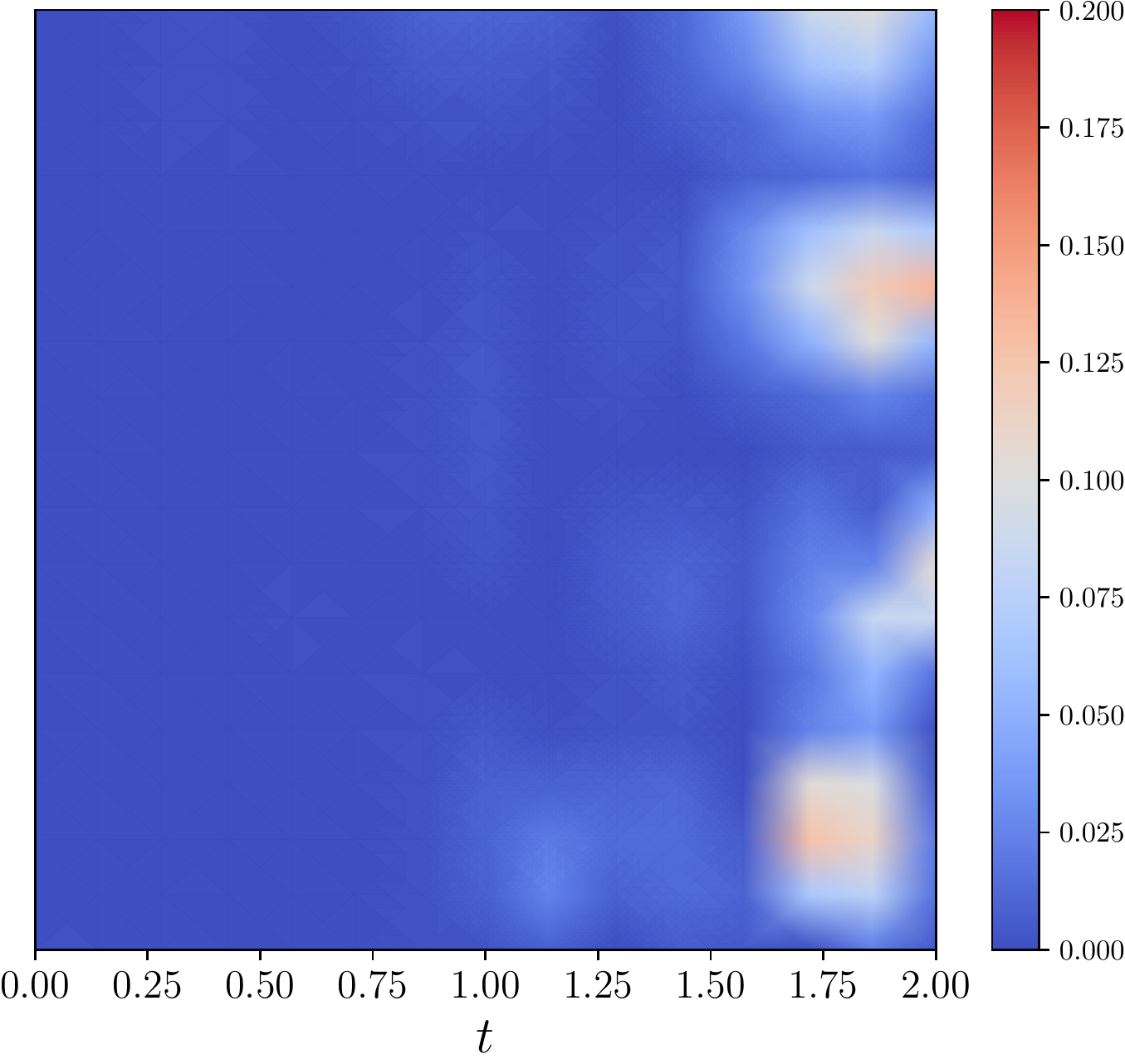}
	\label{fig: err_pulse_HF_lstm}}
    \hspace*{-3mm}
	\subfigure[3-step feed-forward]
	{\includegraphics[width=0.253\linewidth]{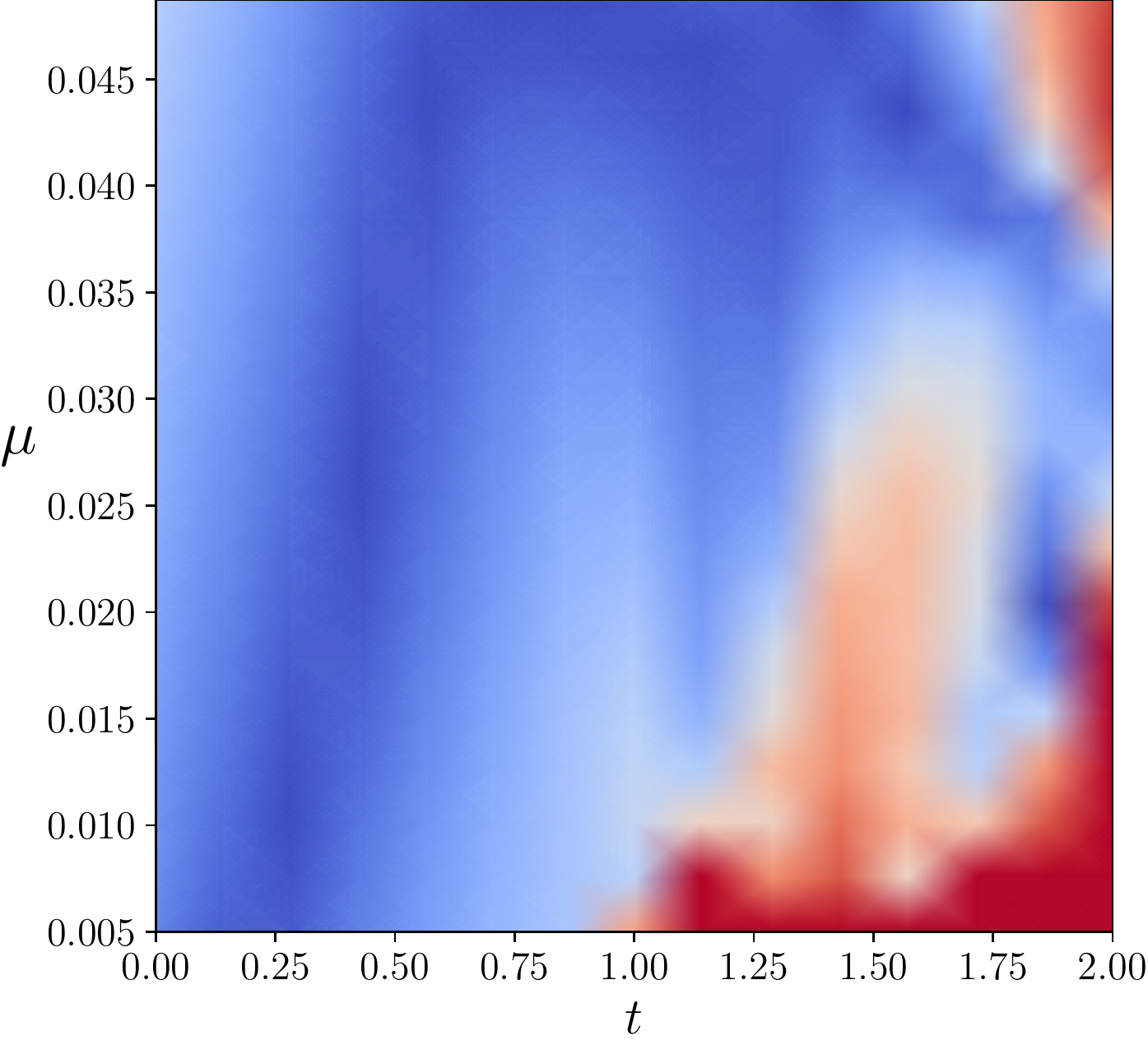}
	\label{fig: err_pulse_MF_FF}}
	\subfigure[MF Intermediate]
	{\includegraphics[width=0.229\linewidth]{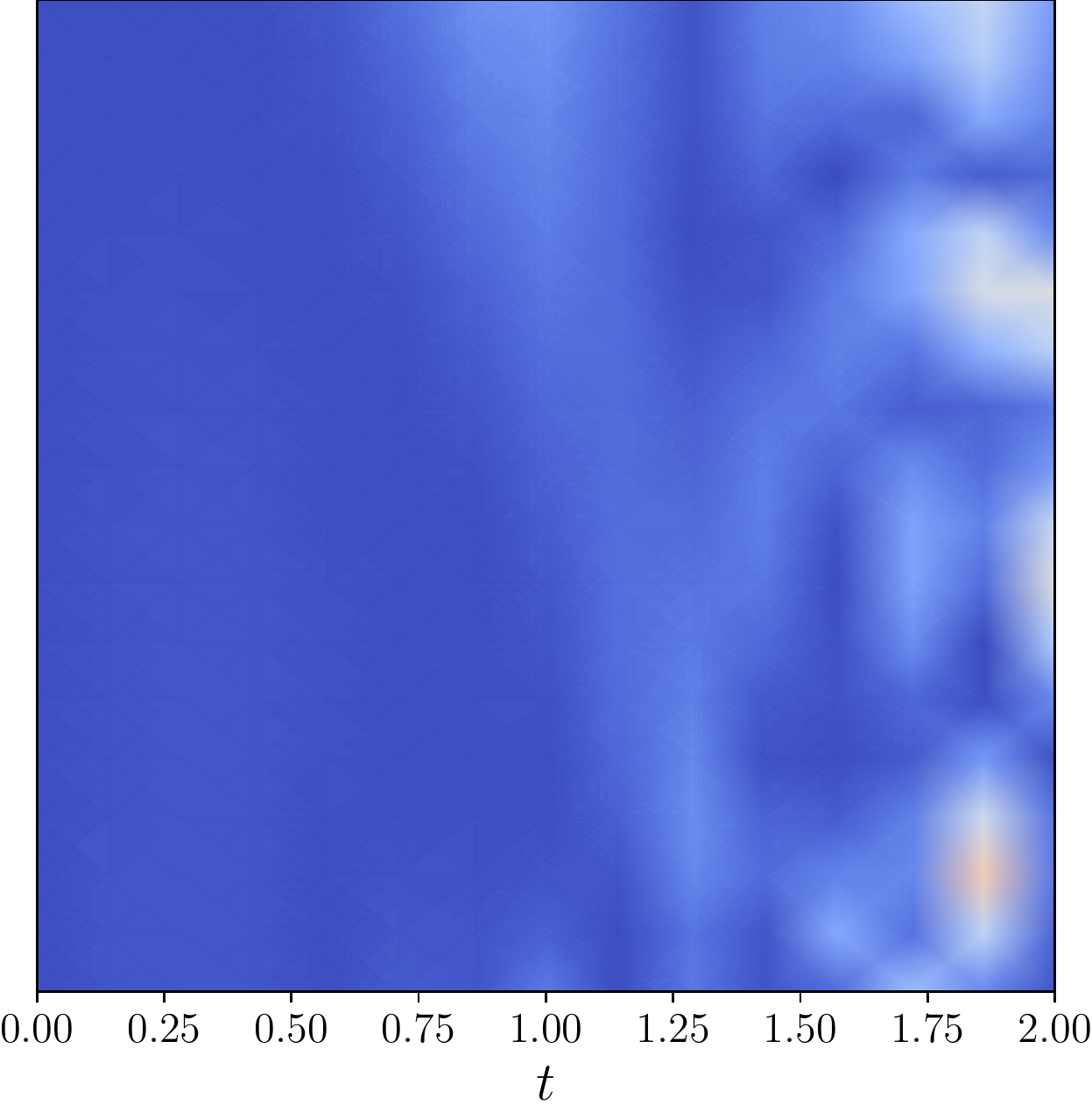}
	\label{fig: err_pulse_inter}}
	\subfigure[MF 2-step LSTM]
	{\includegraphics[width=0.229\linewidth]{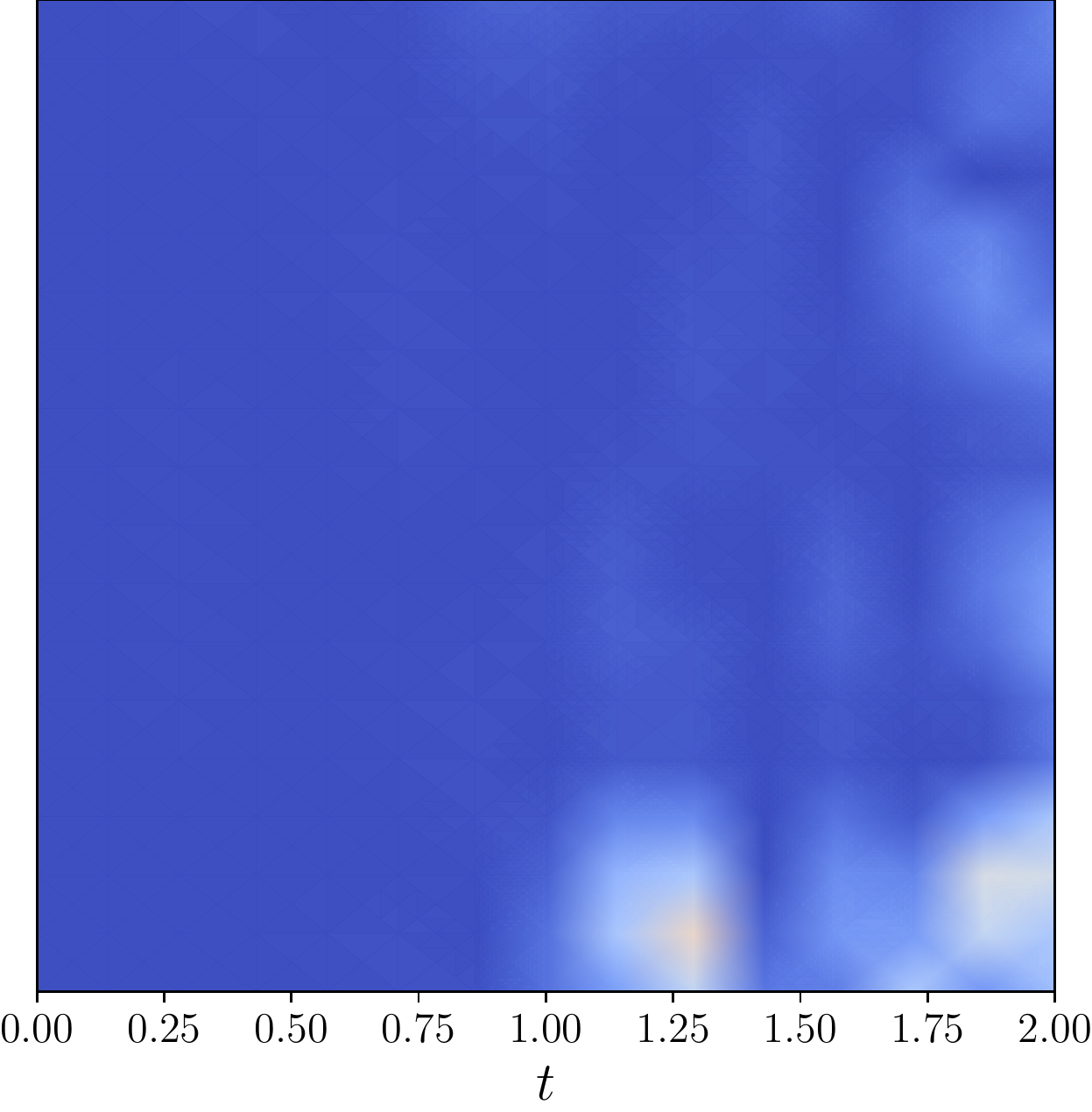}
	\label{fig: err_pulse_2step}}
	\subfigure[MF 3-step LSTM]
	{\includegraphics[width=0.249\linewidth]{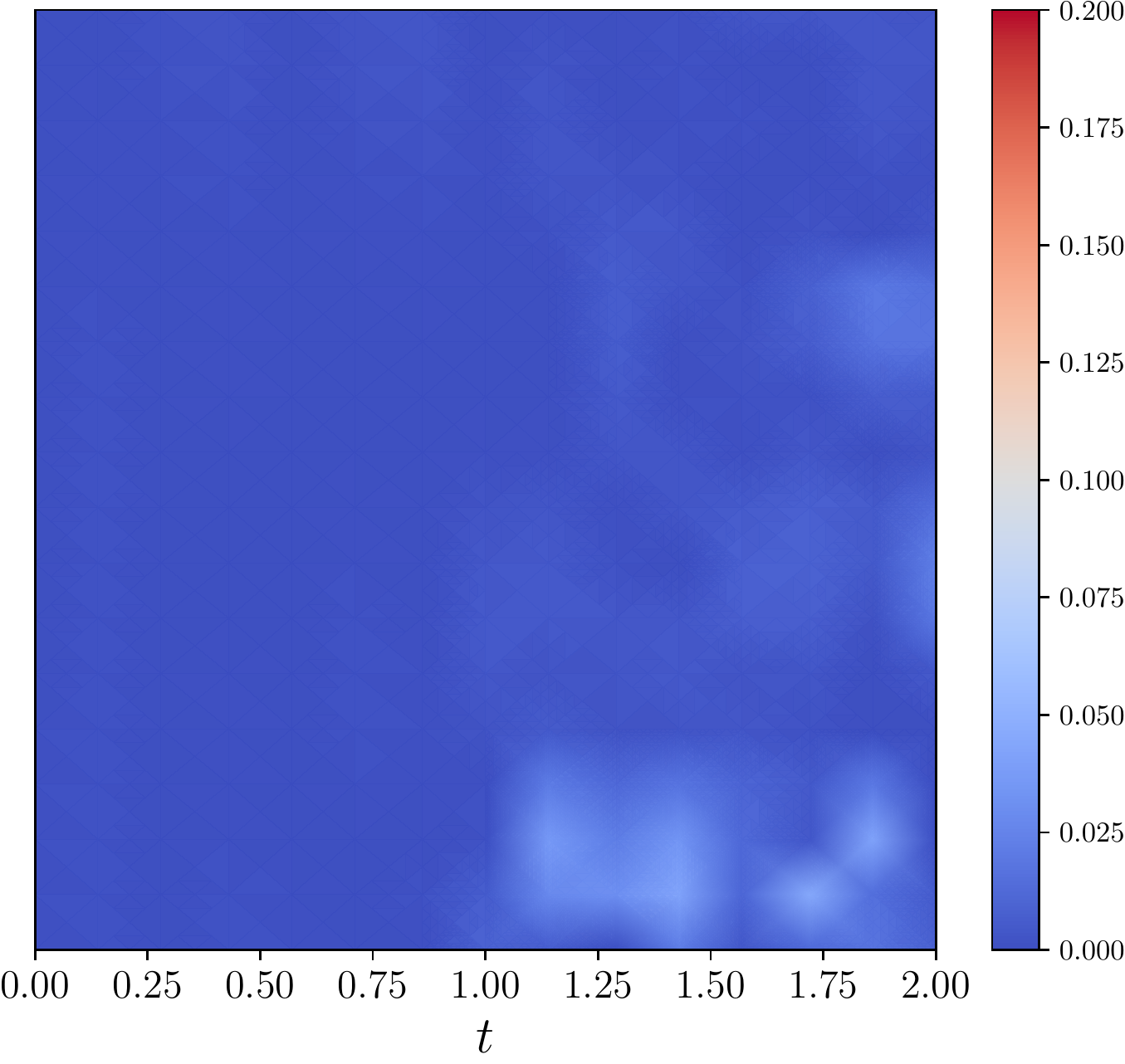}
	\label{fig: err_pulse_3step}}

	\caption{Absolute values of the errors in different NN model predictions of the electrical potential $\nu$ at  $\bar{x} = 0.5$.}
    \label{fig: errors_pulse}
\end{figure}

Our goal in this subsection is to demonstrate the advantages of the proposed MF LSTM surrogates over single-fidelity regressions. As shown in Table \ref{tab: MF_Pulse} and Fig. \ref{fig: errors_pulse}, we train eight different NN models and compare their predictions on the test set. LF and HF feed-forward networks are the single-fidelity models with only feed-forward dense layers, trained on $\mathcal{T}_\texttt{LF}$ and $\mathcal{T}_\texttt{HF}$, respectively. LF LSTM and HF LSTM are their respective extensions incorporating LSTM layers. Intermediate, 2-step and 3-step LSTM networks are the proposed MF models with LSTM layers, while the 3-step feed-forward network is the counterpart of the latter without LSTM layers. The corresponding test MSE values are collected in Table \ref{tab: MF_Pulse}, and the absolute discrepancy values between the model predictions and the HF test values are illustrated in Fig. \ref{fig: errors_pulse}.
This allows us to better understand how the generalization quality varies over the parameter-time domain $\mathcal{P}\times[0,T]$. 

We first note that, in the single-fidelity modeling with LF data $\mathcal{T}_\texttt{LF}$, i.e., the LF feed-forward and the LF LSTM cases,  prediction errors are very large. Although a quite large training set is available,  predictions are poor due to the low accuracy of the LF data, hence the possible advantage of changing architecture is limited by the data quality. In fact, using LSTMs does not yield any remarkable improvement, and both the LF feed-forward and the LF LSTM models show comparable predictive capability.
In the single-fidelity modeling with HF data $\mathcal{T}_\texttt{HF}$, the HF LSTM model results in a clear improvement compared to the HF feed-forward model. As seen in Fig \ref{fig: err_pulse_HF_lstm}, however, there are regions where the discrepancy between HF and predicted solutions is large. Roughly equispaced over the parameter range, these regions correspond to the parameter configurations  for which we have no HF information. 
This implies that, although incorporating LSTM layers has resulted in an improvement when approximating the time dependency, the single-fidelity HF LSTM model has poor generalization with respect to $\mu$ due to the limited availability of training data.

On the other hand, we note that the MF LSTM models (intermediate, 2-step, and 3-step) achieve much better accuracy and allow to successfully generalize with respect to both parameter and time. These models, especially the 3-step LSTM (Fig. \ref{fig: err_pulse_3step}), show uniformly low prediction errors over the time-parameter domain. We stress that these advantages result from the coupling of MF strategies with LSTM-based architectures. In fact, when replacing LSTM layers with dense ones in the best performing MF LSTM model -- the 3-step LSTM, we notice a dramatic deterioration of the predictive capability, as can be clearly seen in both Table~\ref{tab: MF_Pulse} and Fig. \ref{fig: err_pulse_MF_FF}. 
In conclusion, by leveraging the high quality of limited HF data and the intensive domain exploration of the LF data, the MF LSTM models are shown to be efficient tools for accurately predicting time-parameter-dependent quantities of interest that exhibit complex behaviors, e.g., the steep fronts of action-potential in this example.

\begin{table}[t!]
	\renewcommand{\arraystretch}{1.25}
	\centering%
		\caption{Mean squared errors (MSE) calculated on the test set in example (I). }
	\begin{tabular}[b]{c|c|c|c|c}
    \hline
	{{\textbf{Model}}} &{{LF feed-forward}} &  {{HF feed-forward}} &    {{LF LSTM}} & {{HF LSTM}} \\
		\hline\hline
		{{\textbf{Test MSE}}} & ${1.06 \times 10^{-1} }$ & ${2.99 \times 10^{-3}}$ & ${1.11 \times 10^{-1} }$ & ${7.87 \times 10^{-4} }$ \\
		\hline
    \multicolumn{4}{c}{\vspace{-0.3cm}}\\
    \hline
    	{{\textbf{Model}}} &
	    {{MF 3-step feed-forward}} & {{MF Intermediate}} & {{MF 2-step LSTM}} & {{MF 3-step LSTM}} \\
    \hline\hline
        {{\textbf{Test MSE}}} &
    	     ${1.08\times 10^{-2}}$ &	${5.87\times 10^{-4}}$ & ${3.61 \times 10^{-4} }$ & ${6.73 \times 10^{-5} }$\\
		\hline
    \end{tabular}
	\label{tab: MF_Pulse}
\end{table}

\section{Numerical example (II): Fluid flow around a cylinder}
\label{section: NS}
We now consider the estimation of both drag and  lift coefficients associated with a viscous, incompressible fluid flow around a cylinder. The problem is described by the following unsteady Navier-Stokes equations:
\begin{equation}
    \begin{aligned}
           \rho \frac{\partial \vb{v}}{\partial t} - \rho \vb{v}\cdot \nabla\vb{v} - \nabla \cdot \bm{\sigma}(\vb{v},p) &=  \vb{0} \qquad &&(\vb{x},t)\in \Omega \times (0,T) \, ,  \\
    \nabla \cdot \vb{v} &= 0 \qquad &&(\vb{x},t)\in \Omega \times (0,T) \, ,
    \end{aligned}
    \label{eq: NS}
\end{equation}
in which $\vb{v}(\vb{x},t)$ and $p(\vb{x},t)$ are respectively the velocity and pressure fields, $\rho = 1.0 \text{ kg/m}^3$ is the fluid density,  $\bm{\sigma}(\vb{v},p)=-p\vb{I}+2\nu\bm{\epsilon}(\vb{v})$ is the stress tensor, $ \bm{\epsilon}(\vb{v}) = \frac{1}{2}(\nabla\vb{v}+\nabla^\text{T}\vb{v}) $ is the strain tensor, and $\nu$ is the kinematic viscosity. 
The domain, $\Omega = (0, 2.2) \times (0, 0.41) \text{\textbackslash} B_r(0.2,0.2)$ with $r = 0.05$, represents a 2D channel, while the omitted disc, $B_r$, is the obstacle.
Moreover, we prescribe the following boundary and initial conditions:
\begin{equation}
  \begin{aligned}
     \vb{v} &= \vb{0} \qquad &&(\vb{x},t)\in \Gamma_{D_1} \times (0,T) \, , \\
    \vb{v} &= \vb{h} \qquad &&(\vb{x},t)\in \Gamma_{D_2} \times (0,T)  \, , \\ 
    \bm{\sigma}(\vb{v},p)\vb{n} &= \vb{0} \qquad &&(\vb{x},t)\in \Gamma_{N} \times (0,T)  \, , \\ 
    \vb{v}(\vb{x}, 0) &= \vb{0} \qquad && \,\vb{x} \in \Omega \, ,
  \end{aligned}
\end{equation}
i.e., a no-slip condition on $\Gamma_{D_1}$, a parabolic inflow  
\begin{equation}
\vb{h}(\vb{x},t) = \left(\frac{4U(t)x_2(0.41-x_2)}{0.41^2},0\right)\,,  \qquad \text{with }
U(t)=\left\{
  \begin{array}{@{}ll@{}}
    0.75(1-\cos{(\pi t)}), & t <1 \\
    1.5, & t \geq 1
  \end{array}\right.
\end{equation}
on the inlet $\Gamma_{D_2}$, an open boundary condition on the outlet $\Gamma_N$, and a homogeneous initial condition.

Varying the parameter $\nu$ results in a changing Reynolds number. In this specific problem, the Reynolds number can be defined as $Re = LU_\text{mean} / \nu$, in which the average free stream velocity for the parabolic inflow is $U_\text{mean} = 2 \cdot  U_\text{max}/3 = 2 \cdot 1.5 /3= 1 $, and $L =2r = 2 \cdot 0.05 = 0.1$ represents the characteristic length, and thus $Re = 0.1 / \nu $. See also \cite{fresca2021poddlrom,fmfluids1} for further details on this setting.

We are interested in estimating the drag and lift coefficients, denoted by $C_D$ and $C_L$, respectively. They are two important dimensionless quantities defined as follows:
\begin{equation}
    C_D = \frac{2|F_x|}{\rho U^2_\text{mean}}\,,\qquad C_L = \frac{2|F_y|}{\rho U^2_\text{mean}}\,.
\end{equation}
Here $F_x$ and $F_y$ denote the two components of the total force $\vb{F}$ acting on the cylinder, written as
\begin{equation}
    \vb{F} = \oint_{\partial B_r} \left[-p\vb{I}+\nu(\nabla \vb{v} + \nabla^\text{T} \vb{v})\right]\cdot \vb{n}\,\mathrm{d}s\,,
\end{equation}
in which $\vb{n}$ denotes the outer normal vector along the boundary $\partial B_r$. 

Our goal is to predict the evolution of drag and lift coefficients over time as the Reynolds number $Re$ changes. We consider $\mu=Re\in [70,160]$, a range in which the flow is unsteady, and we confine our surrogate modeling within the time interval $t\in[14.5,15.0]$, during which the flow becomes fully developed and presents a periodic behavior.
\begin{figure}[b]
    \centering
    \includegraphics[width=0.8\textwidth]{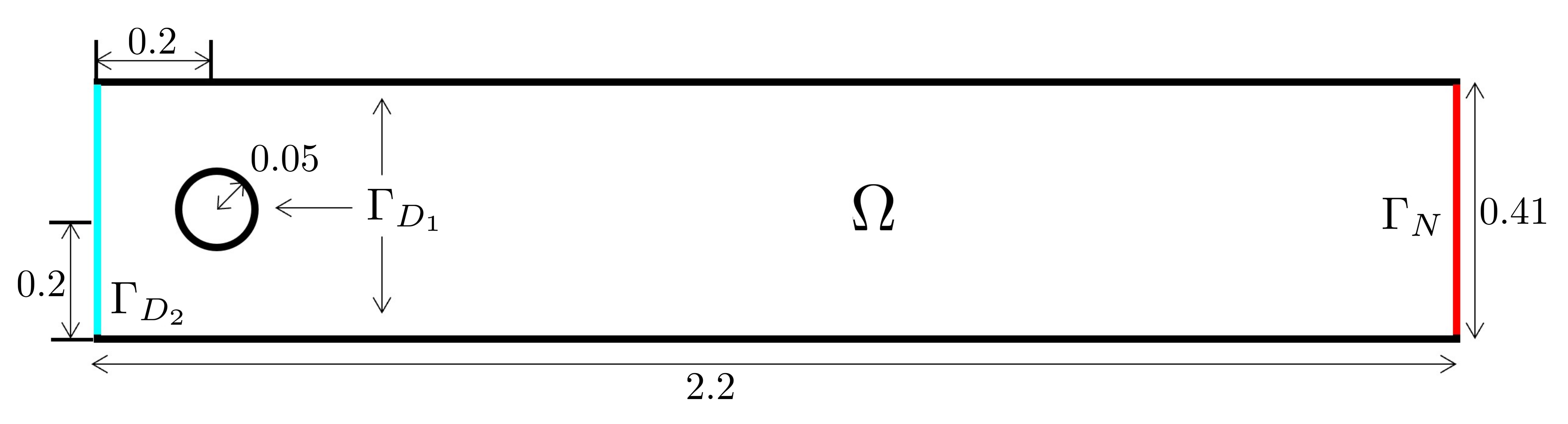}
    \caption{Geometry for the 2-D channel flow around a cylinder. All lengths are measured in meters.}
    \label{fig:NS_geometry}
\end{figure}

\subsection{Multi-fidelity setting}
\begin{figure}[t]
\centering
\subfigure[LF coarse mesh consisting of 7789 triangular elements and 3899 nodes.]{
	\label{subfig:LF_mesh}
	\includegraphics[width=0.74\textwidth]{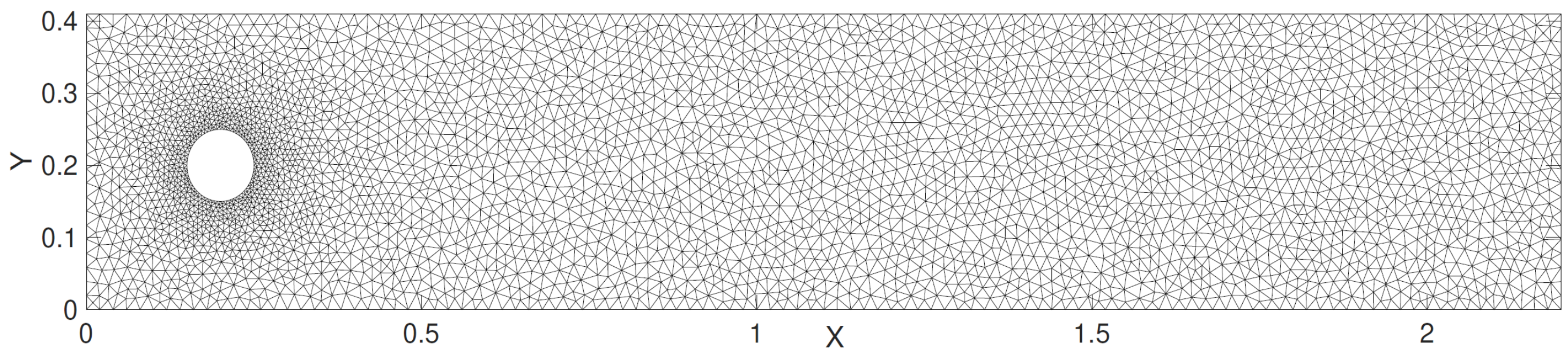} } 
\centering
\subfigure[HF fine mesh consisting of 16478 triangular elements and 8239 nodes.]{
	\label{subfig:HF_mesh}
	\includegraphics[width=0.74\textwidth]{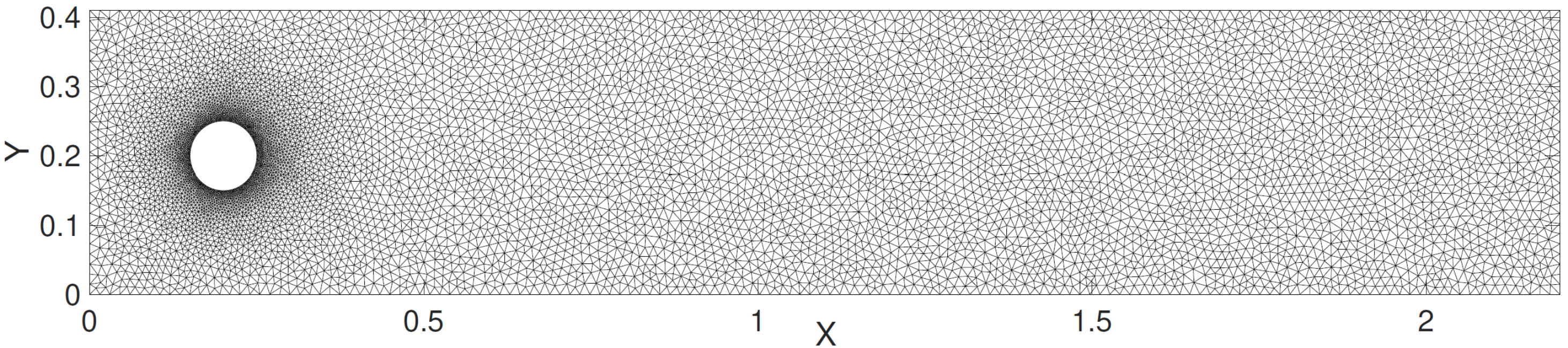} } 
\caption{Computational meshes used for the finite element discretization of the Navier-Stokes equations \eqref{eq: NS}.}
\label{fig: meshes}
\end{figure}

Both the HF and LF models are constructed through the numerical approximation of \eqref{eq: NS} using the MATLAB library redbKIT \cite{redbKIT}, which exploits finite elements and the backward differentiation formula for the spatial and temporal discretizations, respectively.  
The two fidelity levels are distinguished by the mesh and time step sizes. In the HF model we adopt a time step $\Delta t_\texttt{HF} = 0.01$ and a fine computational mesh, while we use $\Delta t_\texttt{LF} = 0.02$ and a coarse mesh in the LF model. The two meshes are displayed in Fig. \ref{fig: meshes}.

As for the data sets, we first consider uniform grids of $N^{\mu}_{\texttt{LF}}=19$ and $N^{\mu}_{\texttt{HF}}=10$ Reynolds number values over the parameter interval $\mathcal{P} = [70,160]$. For each parameter value, we include $N^\text{T}_{\texttt{HF}} = N^\text{T}_{\texttt{LF}} = 26$ uniform steps over the time interval $[t_0,T] =[14.5,15.0]$, and compute the corresponding values of drag and lift coefficients.
In the LSTM models, we consider full-length time sequence by setting the batch subsequence length $K=26$. NNs are trained to approximate the drag and lift coefficients as functions of time $t$ and Reynolds number $Re$. {\color{black}{Moreover, the test set comprises the HF evaluations of drag and lift coefficients for 19 Reynolds number values and 26 time instants that are uniformly spaced in $\mathcal{P}$ and $[t_0, T]$, respectively.}}

Next, we repeat the analysis with a smaller HF training set consisting of $N_\texttt{HF}^\mu=6$ equally spaced Reynolds number values, while keeping $N^\text{T}_\texttt{HF} = 26$, to assess the robustness of the NN models. 
{\color{black}{In this work, we assume that LF data are cheap/easy to obtain in a sufficiently large amount, so that the mapping from the time-parameter inputs to LF outputs can be approximated accurately. A general guidance for choosing the number of LF data is to verify that the LF training and testing errors are controlled within the same order of magnitude, i.e., a good generalization on the LF level is achieved. Hence, what primarily limits the model performance is either the \textit{quality} of LF approximation or the \textit{quantity} of HF data. In this example, we analyze the latter case while fixing the LF model, because our focus is on exploring the models’ prediction and extrapolation capabilities when the HF data are limited. However, we refer to \cite{guo2022multi} for a detailed discussion about the impact of 
varied LF model quality on the feed-forward versions of intermediate, 2-step and 3-step architectures in a parametric differential problem solved through a reduced basis method.}}

\begin{figure}[ht]
	\centering
	\subfigure[Lift coefficient - LF solution]
	{\includegraphics[width=0.335\linewidth]{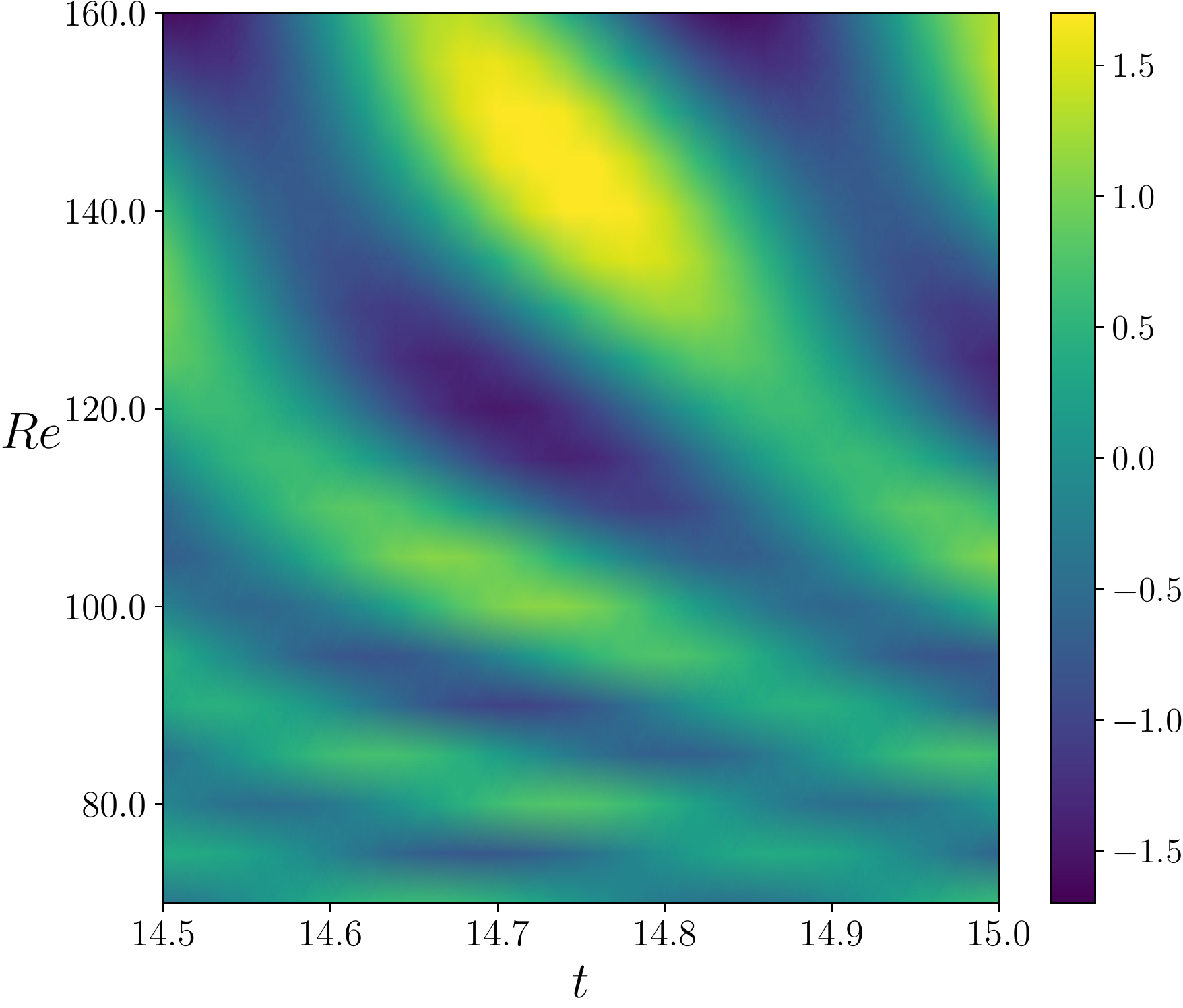}}
	\qquad
	\subfigure[Lift coefficient - HF solution]
	{\includegraphics[width=0.335\linewidth]{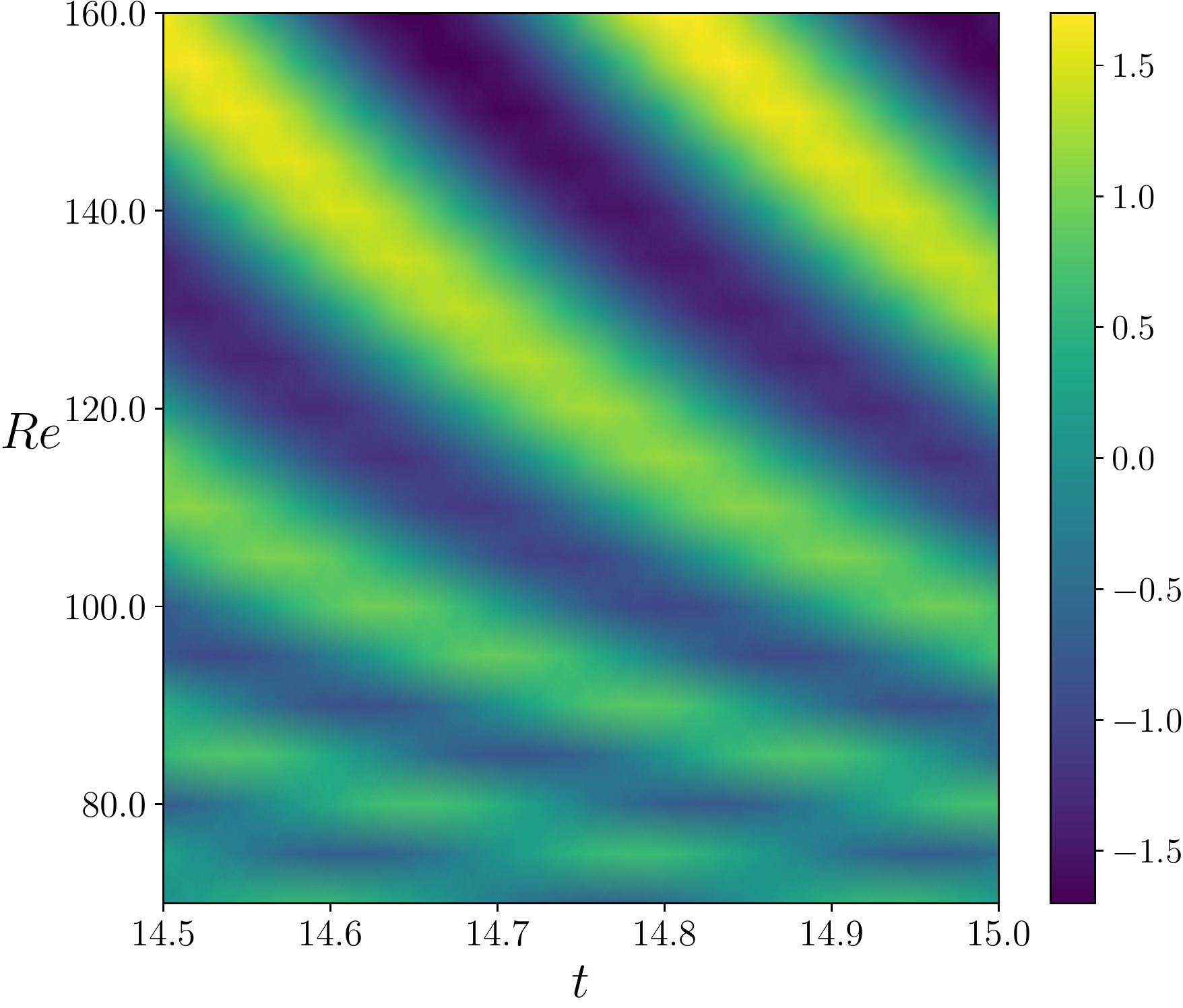}}
    
    \subfigure[Drag coefficient - LF solution]
	{\includegraphics[width=0.335\linewidth]{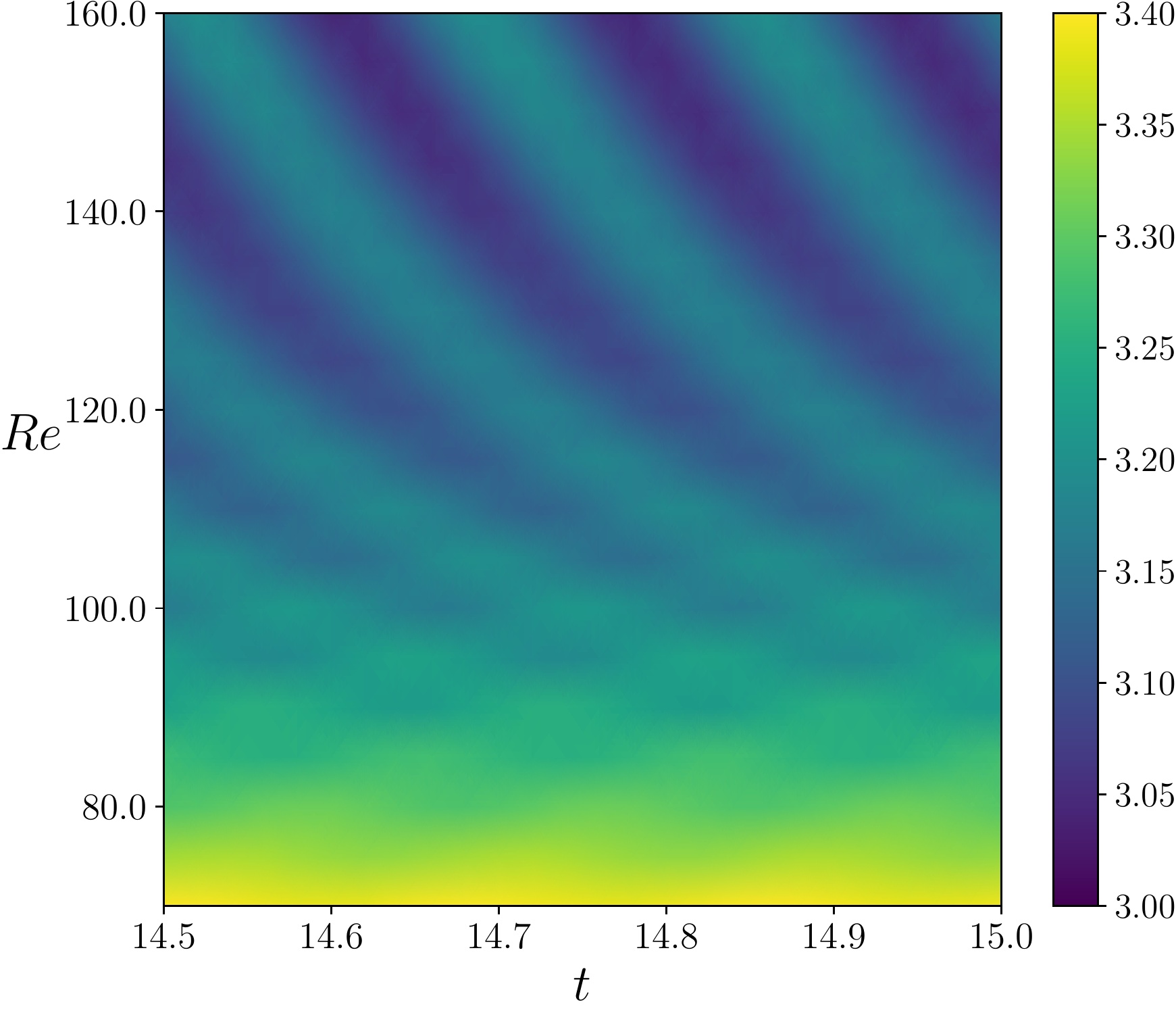}}
	\qquad
	\subfigure[Drag coefficient - HF solution]
	{\includegraphics[width=0.335\linewidth]{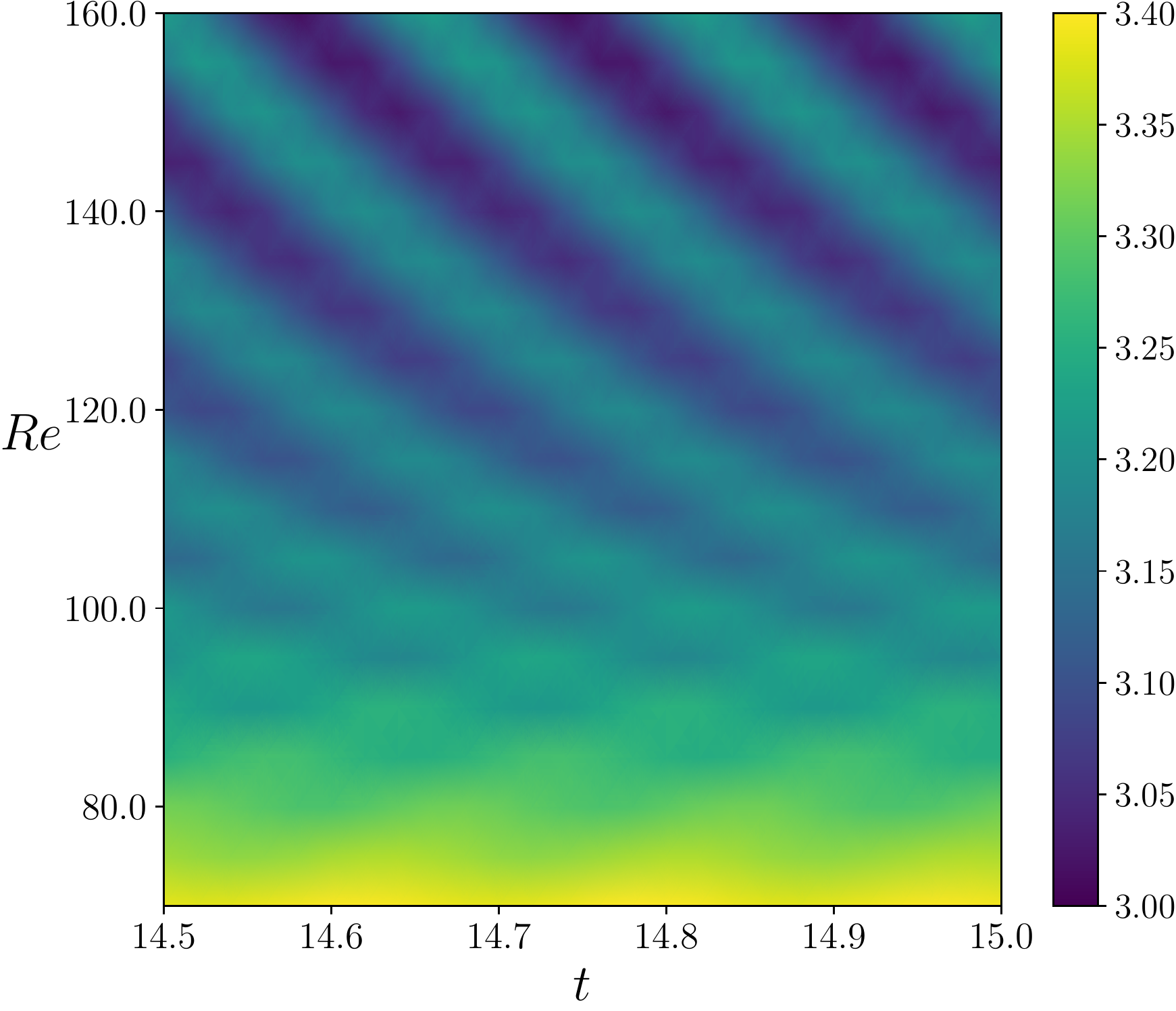}}
	\caption{LF and HF solutions of the lift (above) and drag (below) coefficients.}
    \label{fig: drag_lift_models}
\end{figure}

The LF and HF solutions are shown in Fig. \ref{fig: drag_lift_models}. We note that the drag and lift coefficients exhibit oscillatory patterns over time, varying with respect to the Reynolds number. This may be problematic for the feed-forward NNs that treat time as a generic input entry.

\subsection{Results and discussions 1: Interpolation over the time-parameter domain}
 
\begin{figure}[h!]
	\centering
	\hspace*{-3mm}
	\subfigure[LF feed-forward]
	{\includegraphics[width=0.259\linewidth]{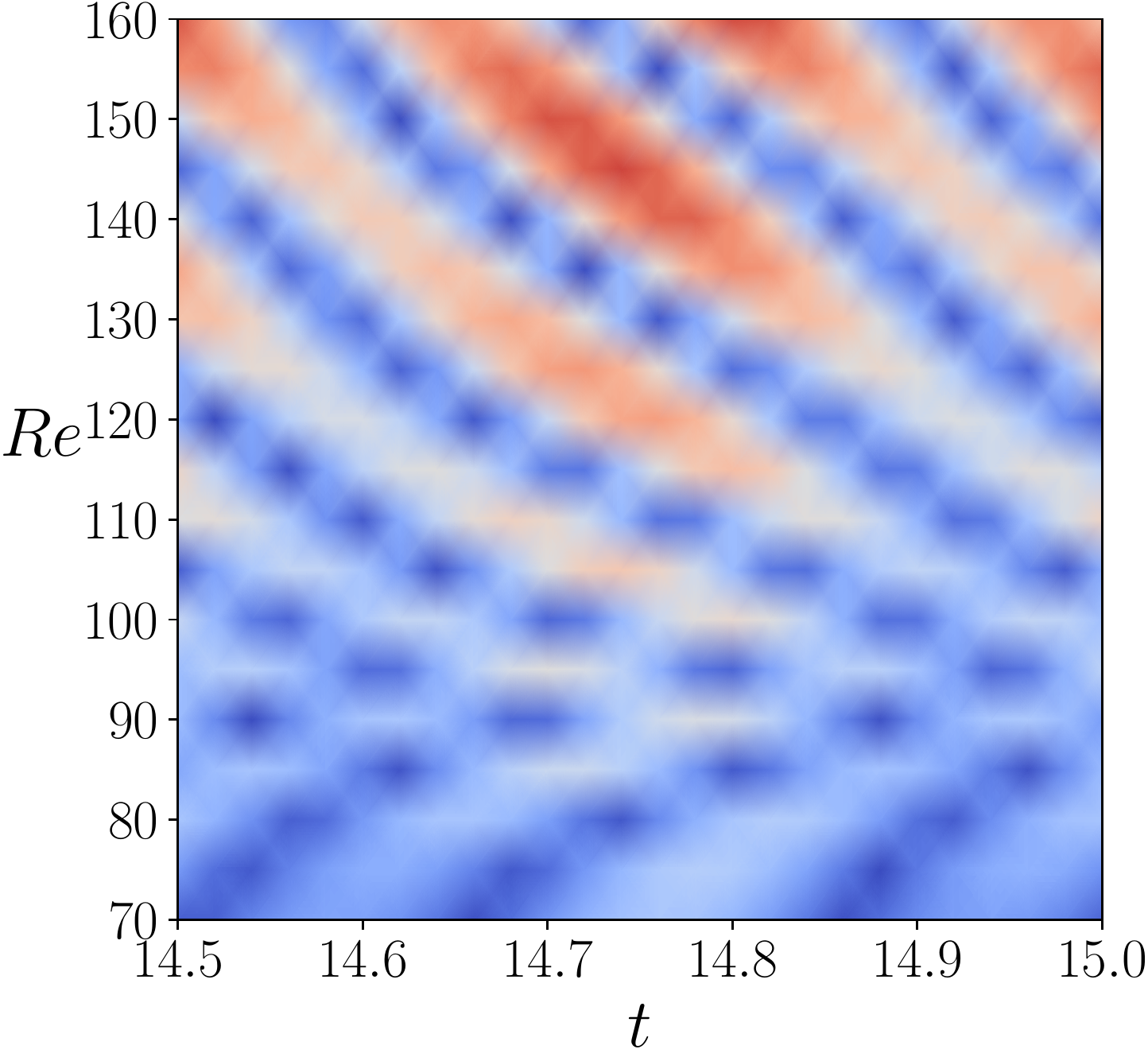}
	\label{fig: err_lift_LF_FF}}
	\subfigure[HF feed-forward]
	{\includegraphics[width=0.229\linewidth]{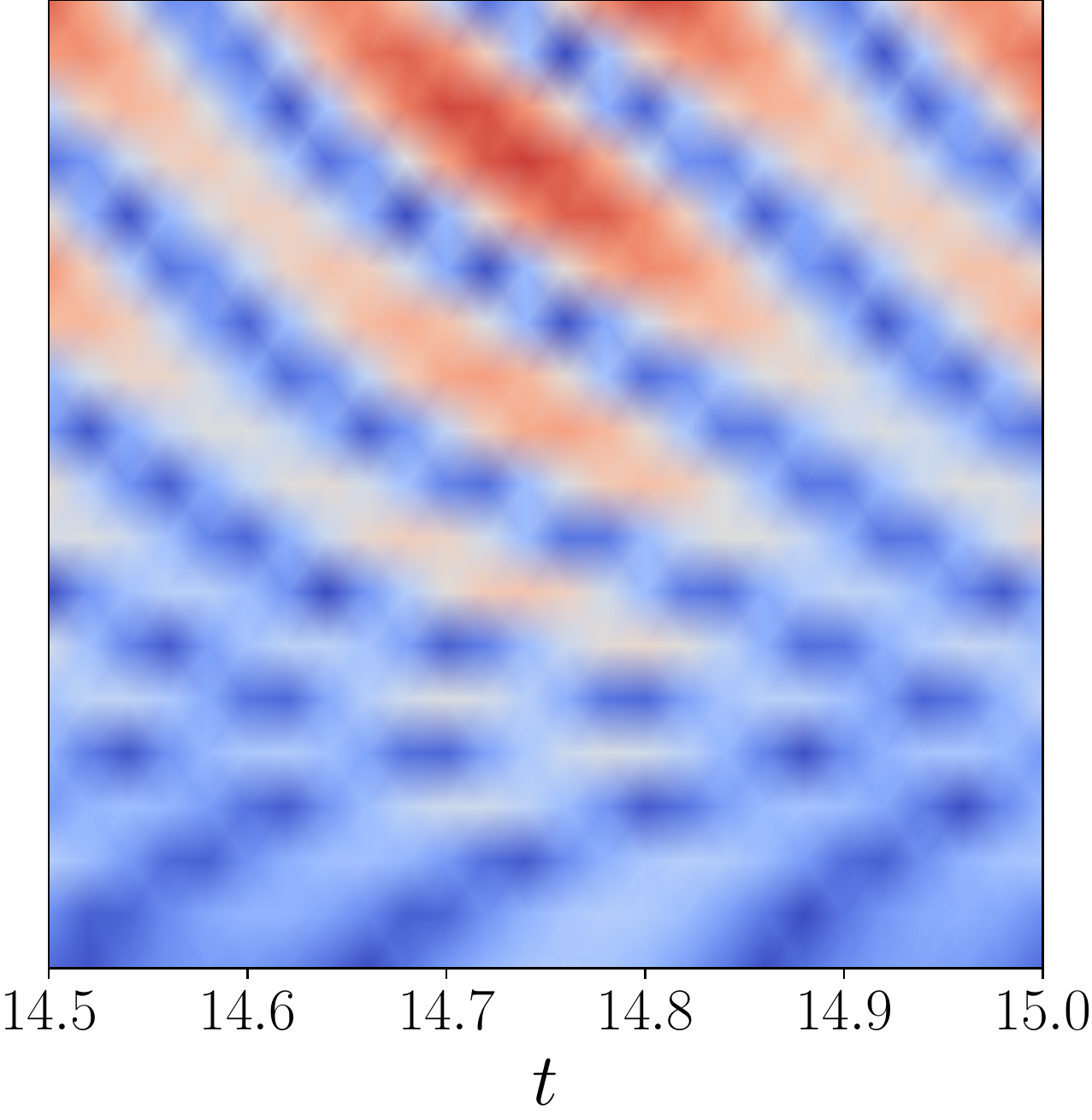}
	\label{fig: err_lift_HF_FF}}
	\subfigure[LF LSTM]
	{\includegraphics[width=0.229\linewidth]{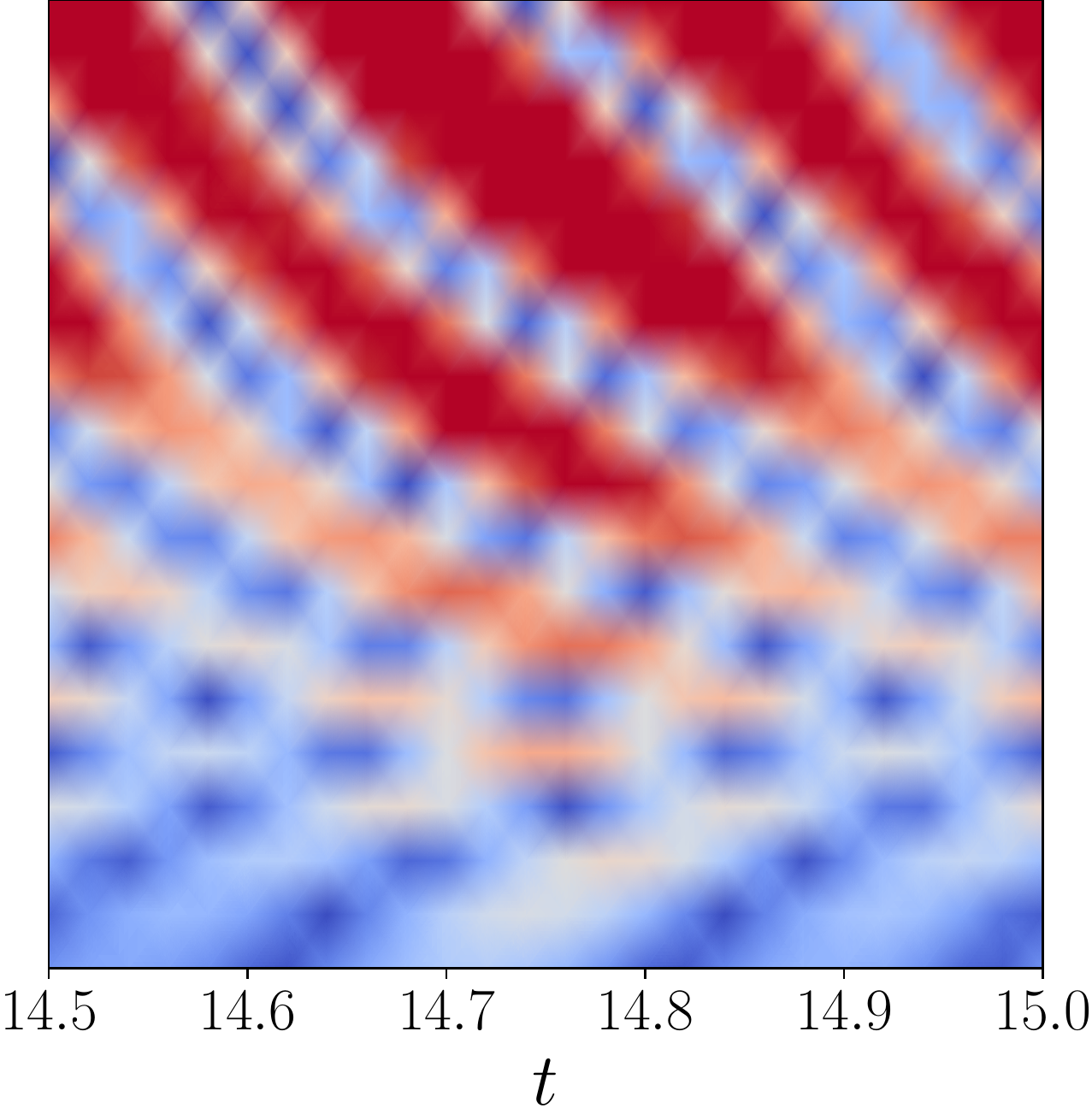}
	\label{fig: err_lift_LF_lstm}}
	\subfigure[HF LSTM]
	{\includegraphics[width=0.242\linewidth]{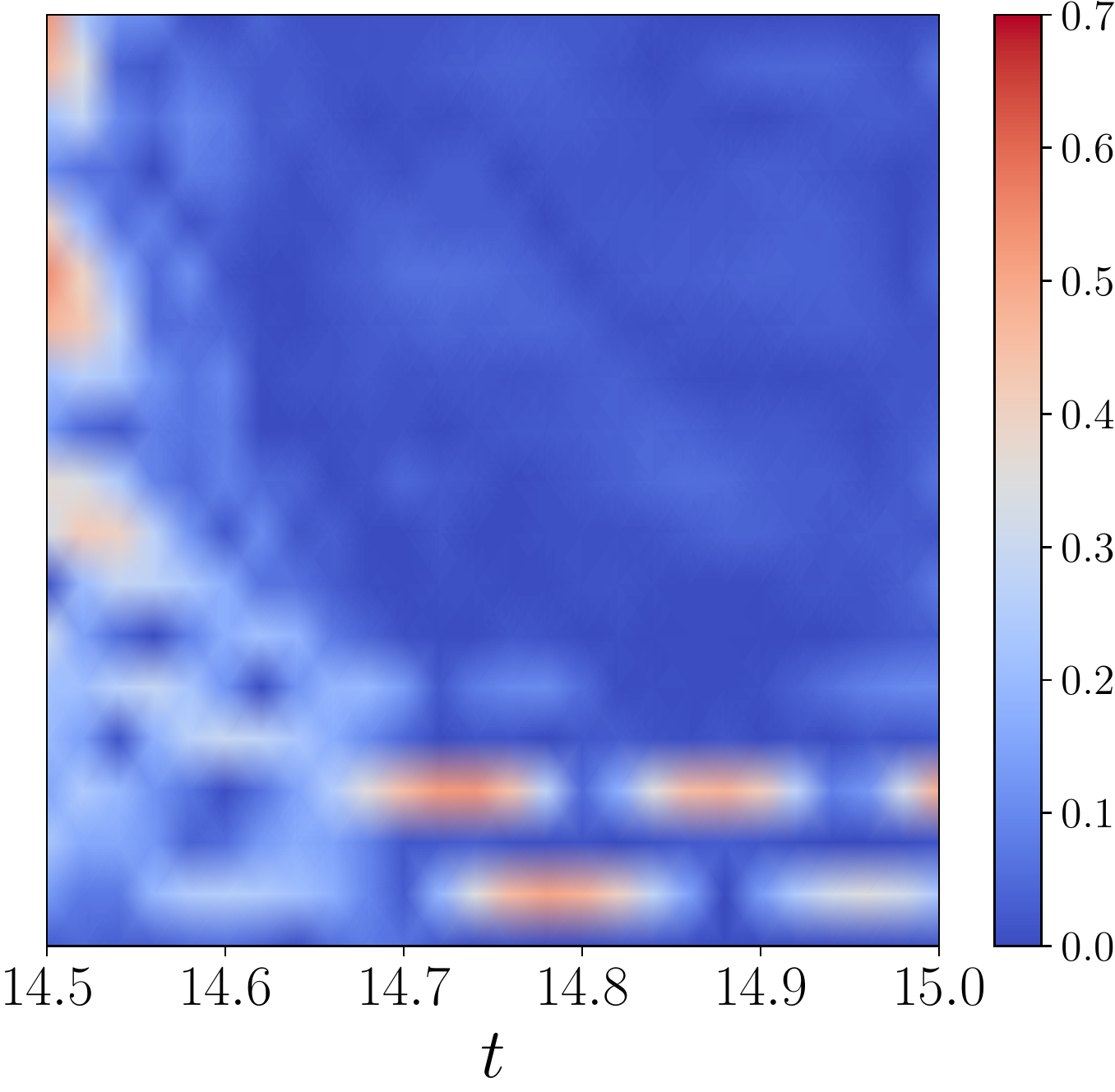}
	\label{fig: err_lift_HF_lstm}}
	\vspace{-3pt}
    \hspace*{-3mm}
	\subfigure[MF 3-step feed-forward]
	{\includegraphics[width=0.259\linewidth]{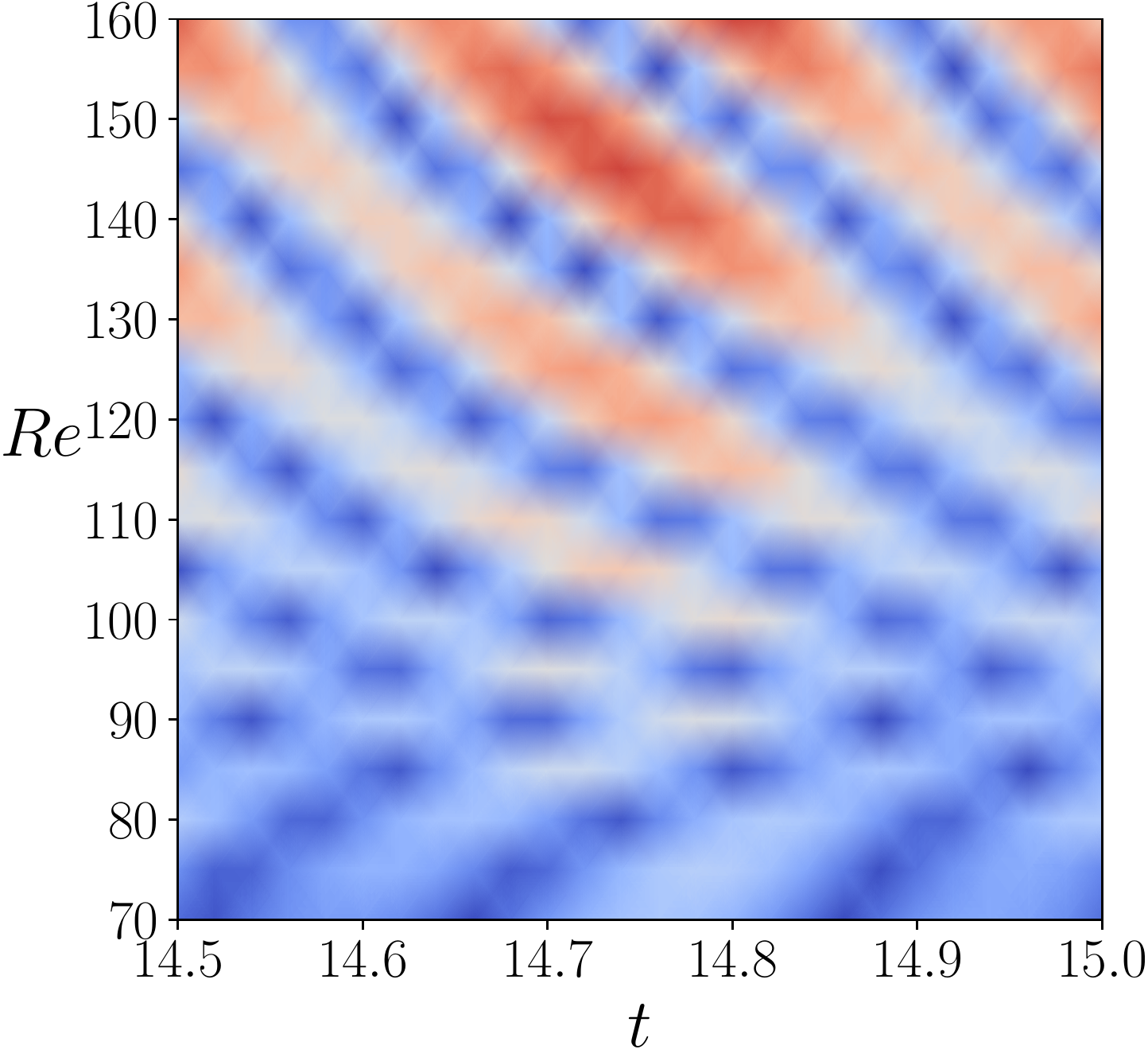}
	\label{fig: err_lift_MF_FF}}
	\subfigure[MF Intermediate]
	{\includegraphics[width=0.229\linewidth]{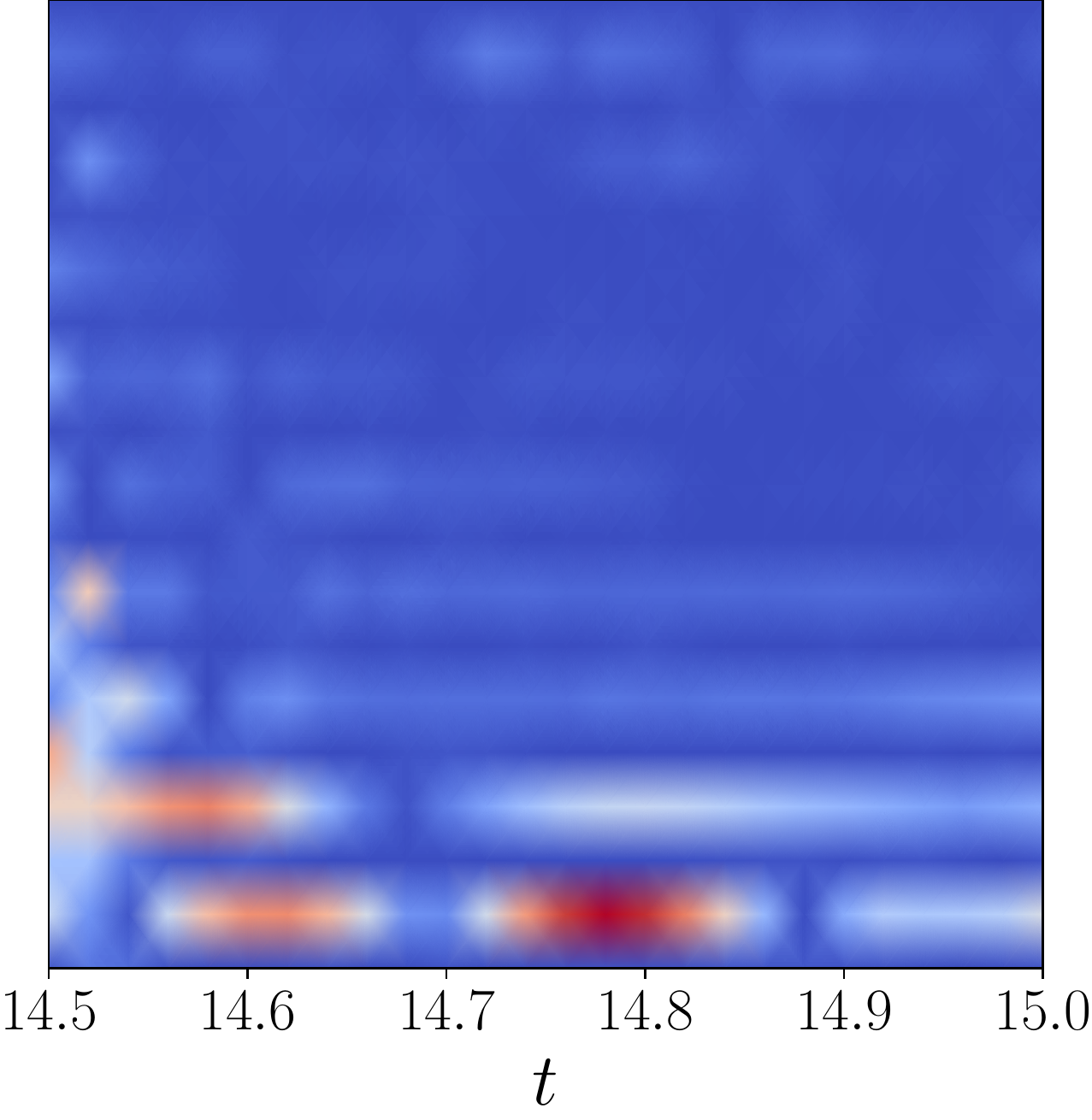}
	\label{fig: err_lift_inter}}
	\subfigure[MF 2-step LSTM]
	{\includegraphics[width=0.229\linewidth]{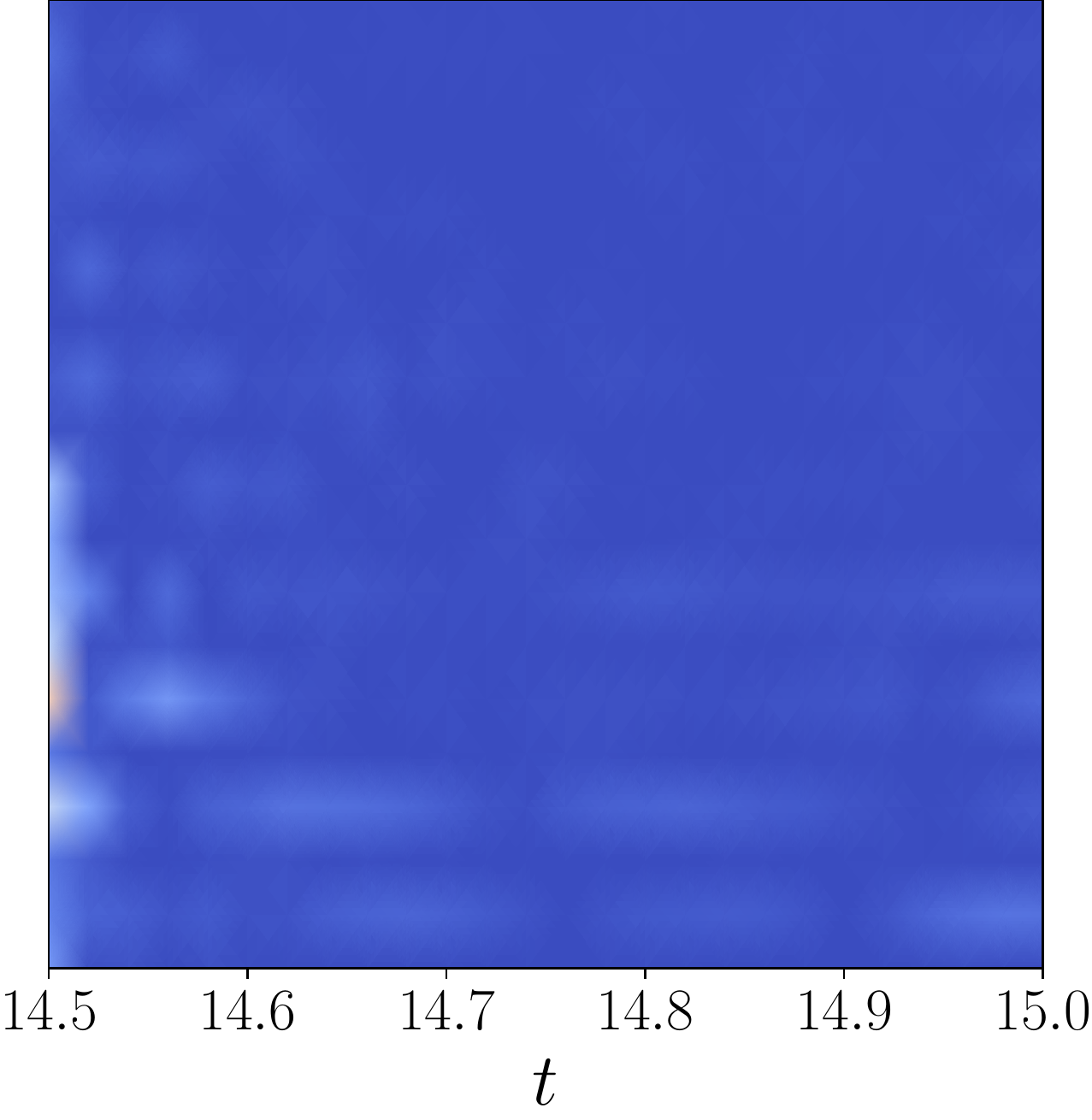}
	\label{fig: err_lift_2step}}
	\subfigure[MF 3-step LSTM]
	{\includegraphics[width=0.242\linewidth]{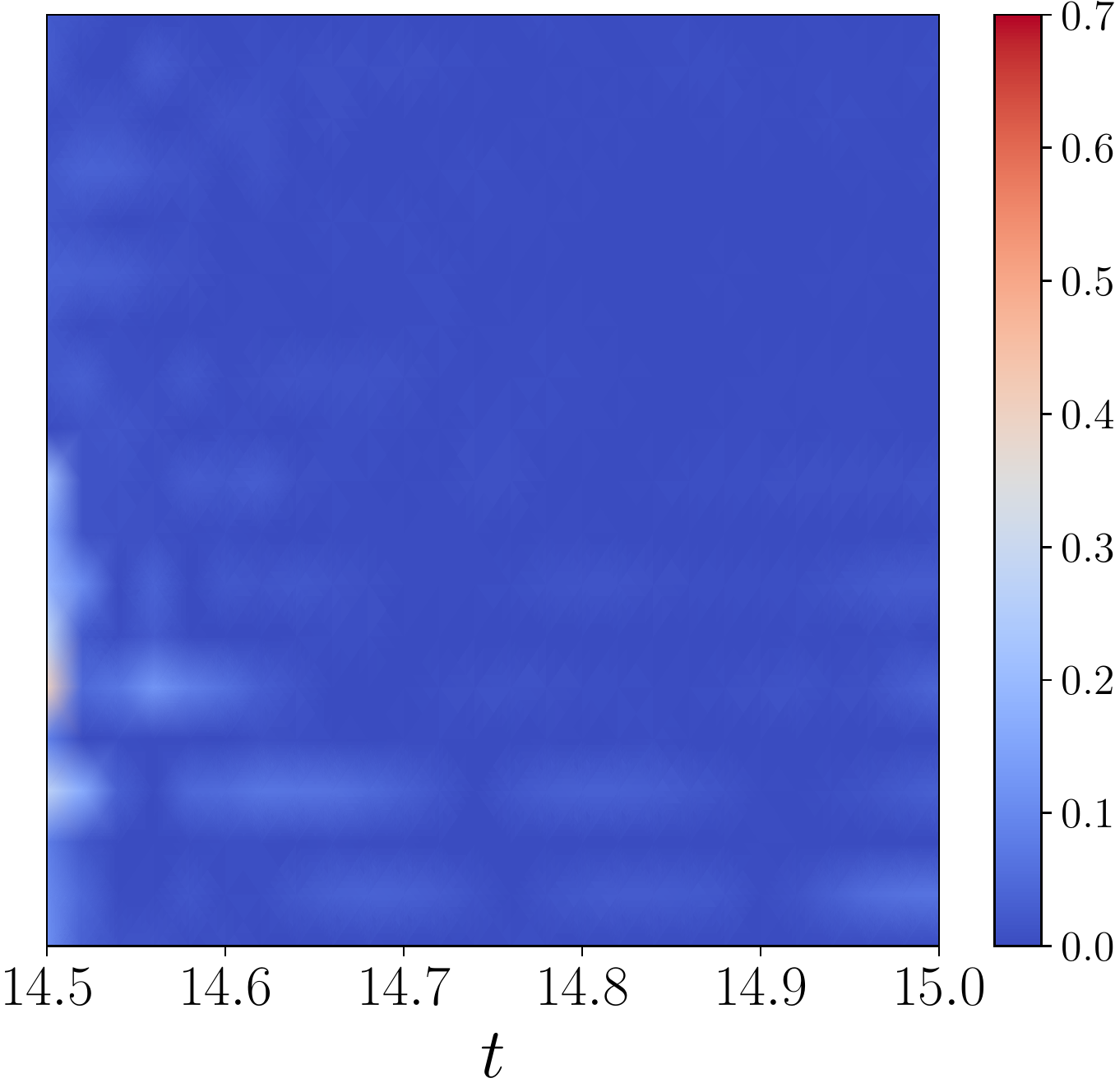}
	\label{fig: err_lift_3step}}
    \vspace{-5pt}
	\caption{Absolute values of the errors in different NN model predictions of the lift coefficient with $N^\mu_\texttt{HF} = 10$ HF training time series. }
    \label{fig: errors_lift}
\end{figure}

\begin{figure}[h!]
	\centering
	\hspace*{-3mm}
	\subfigure[LF feed-forward]
	{\includegraphics[width=0.259\linewidth]{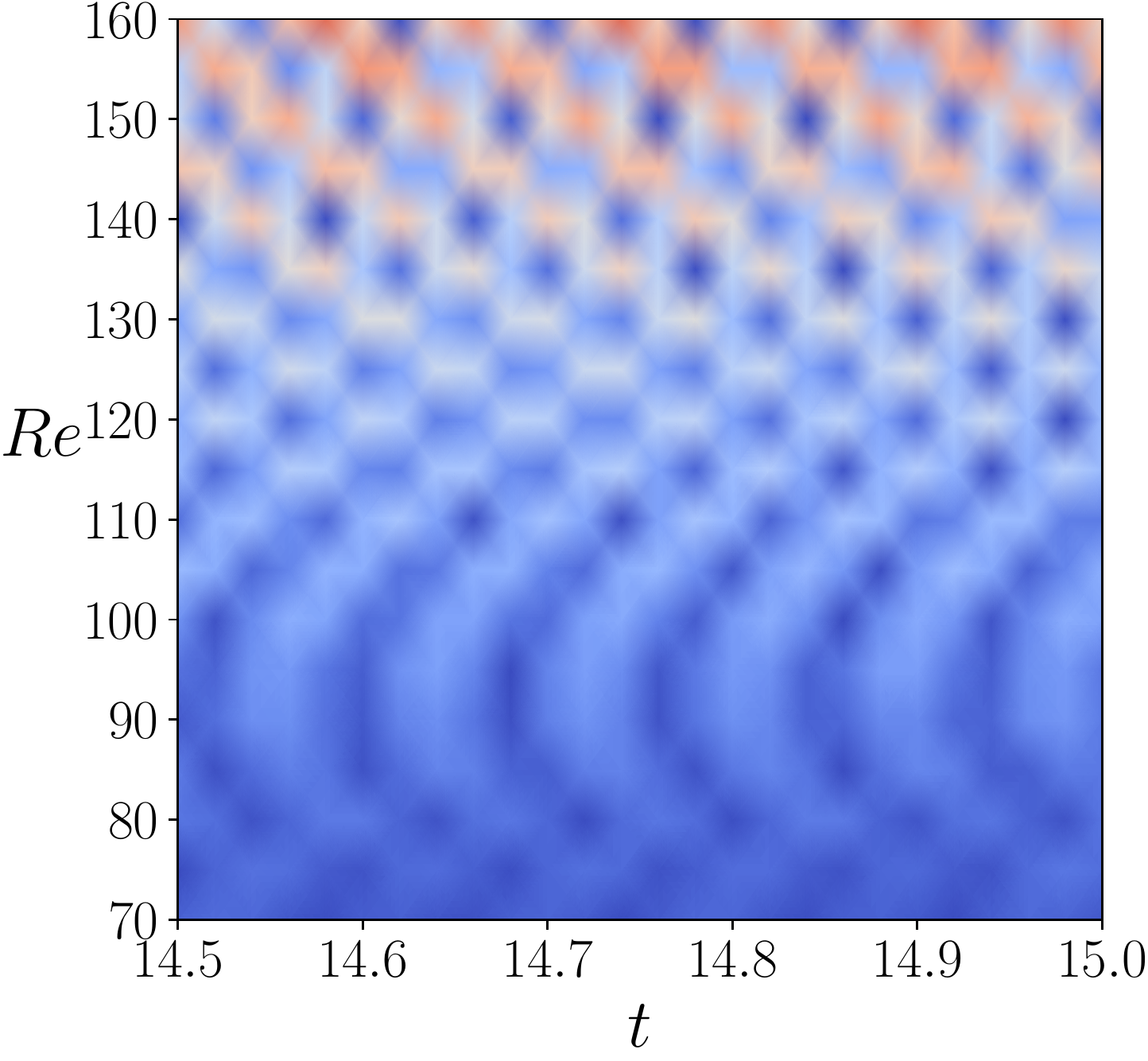}}
	\subfigure[HF feed-forward]
	{\includegraphics[width=0.229\linewidth]{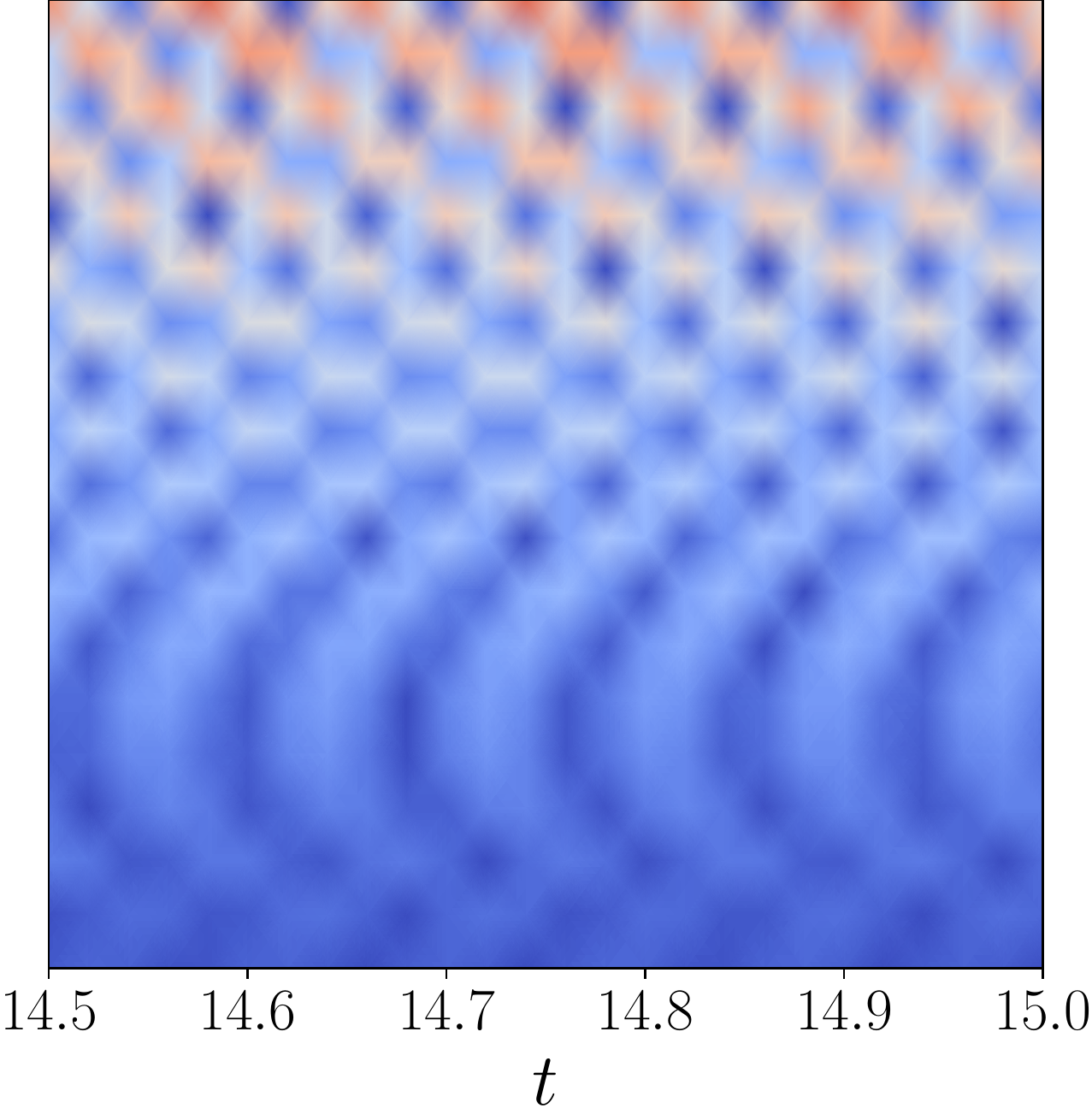}}
	\subfigure[LF LSTM]
	{\includegraphics[width=0.229\linewidth]{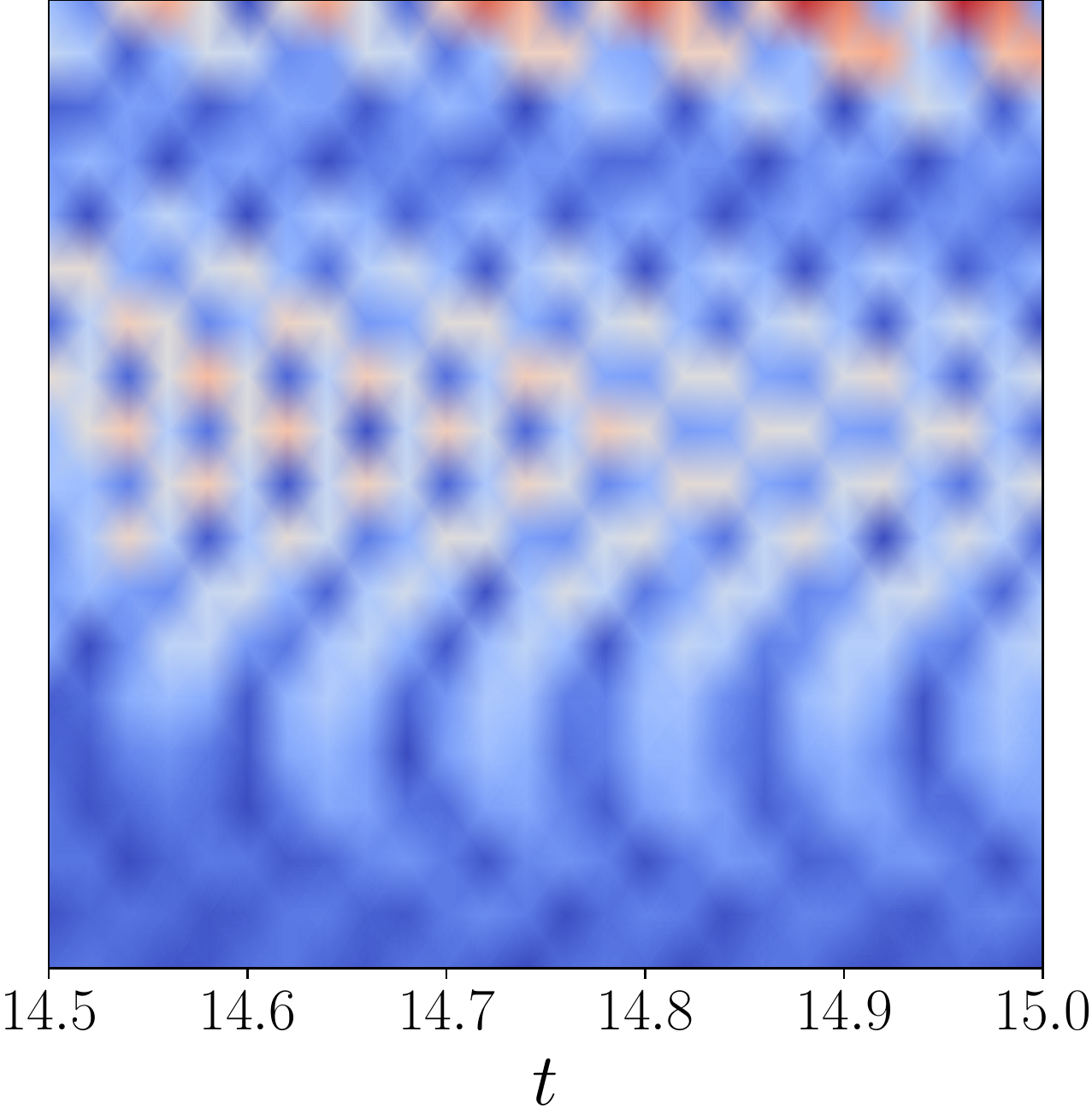}}
	\subfigure[HF LSTM]
	{\includegraphics[width=0.248\linewidth]{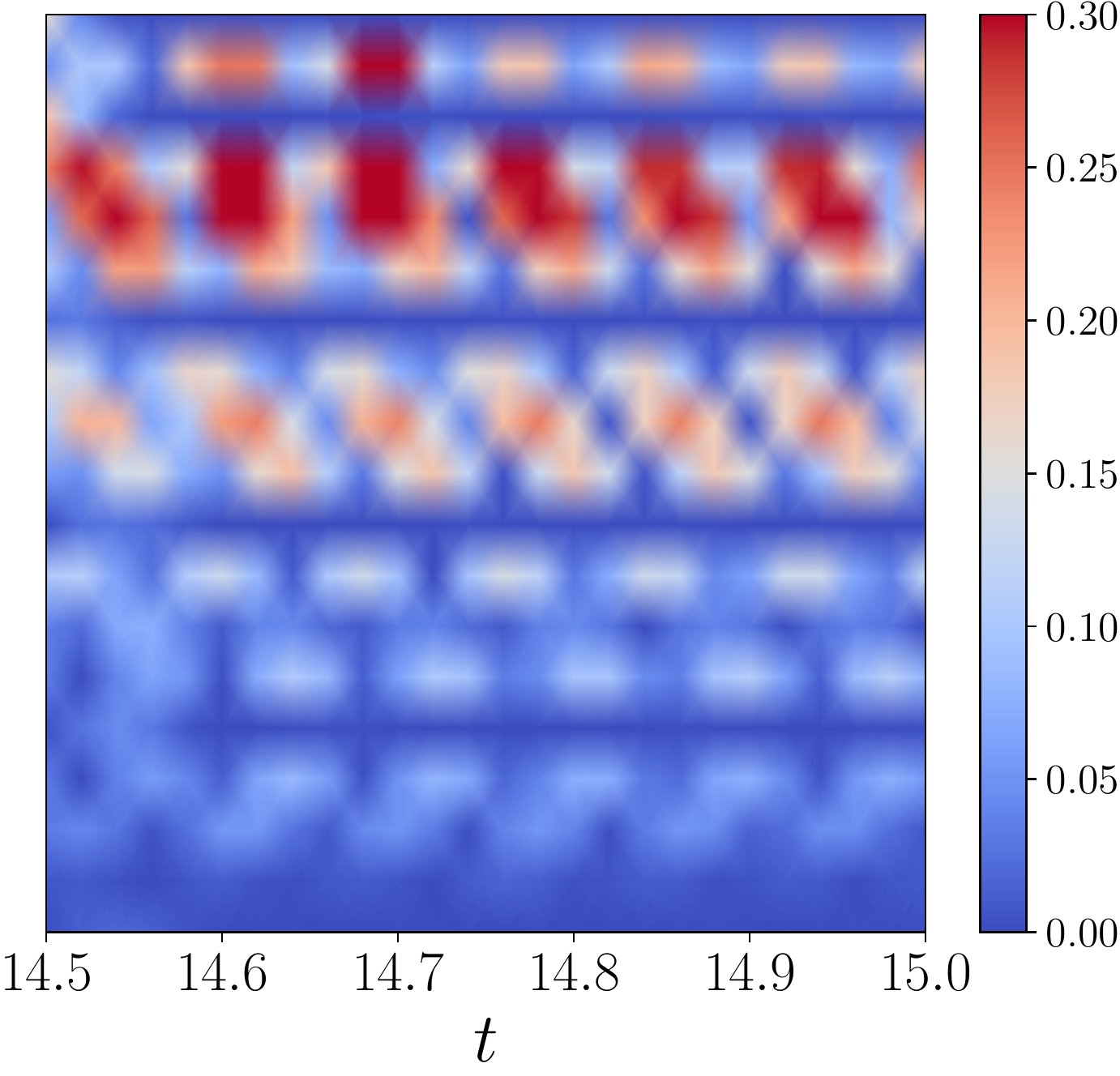}}
	\vspace{-3pt}
    \hspace*{-3mm}
	\subfigure[MF 3-step feed-forward]
	{\includegraphics[width=0.259\linewidth]{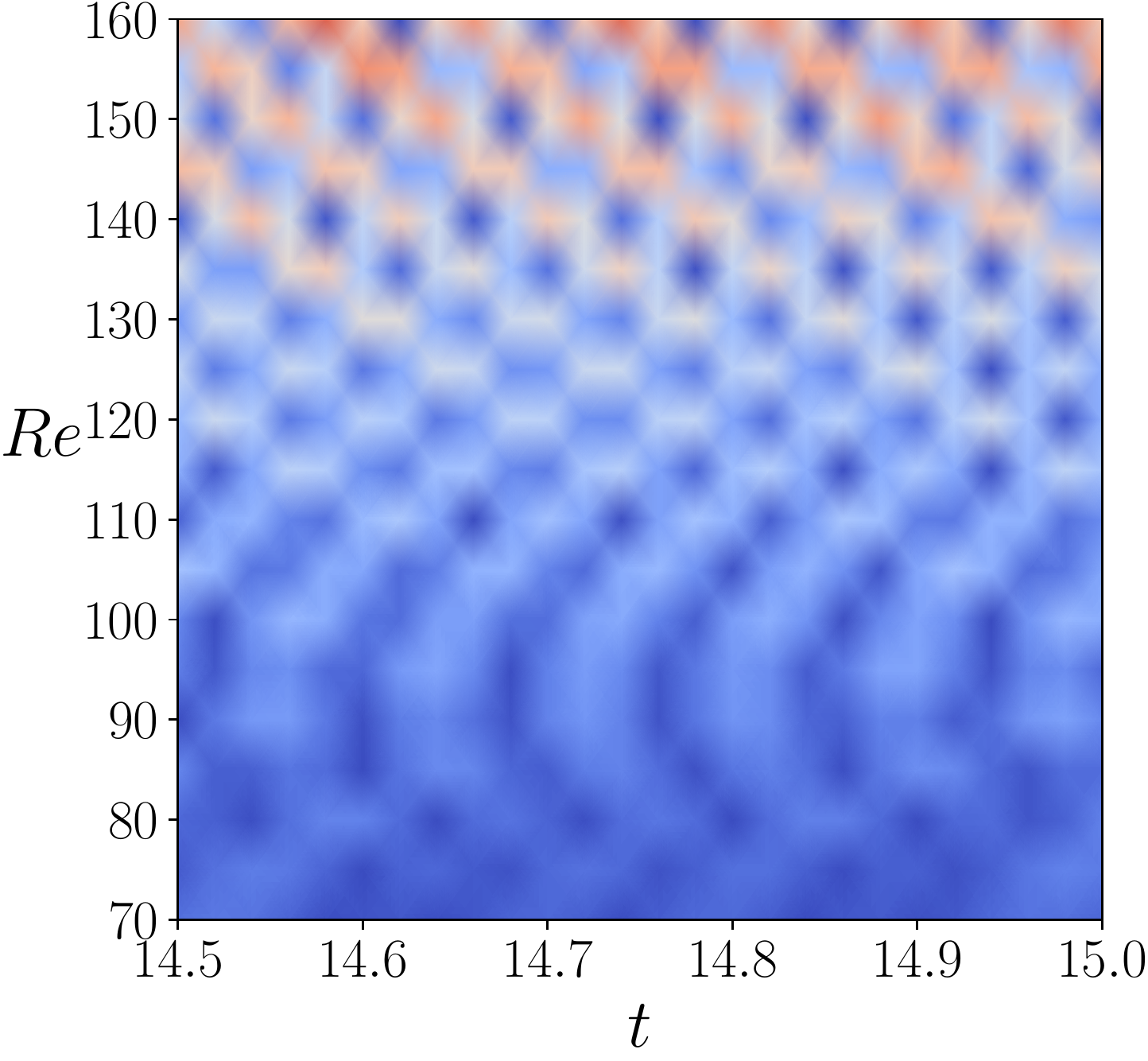}}
	\subfigure[MF Intermediate]
	{\includegraphics[width=0.229\linewidth]{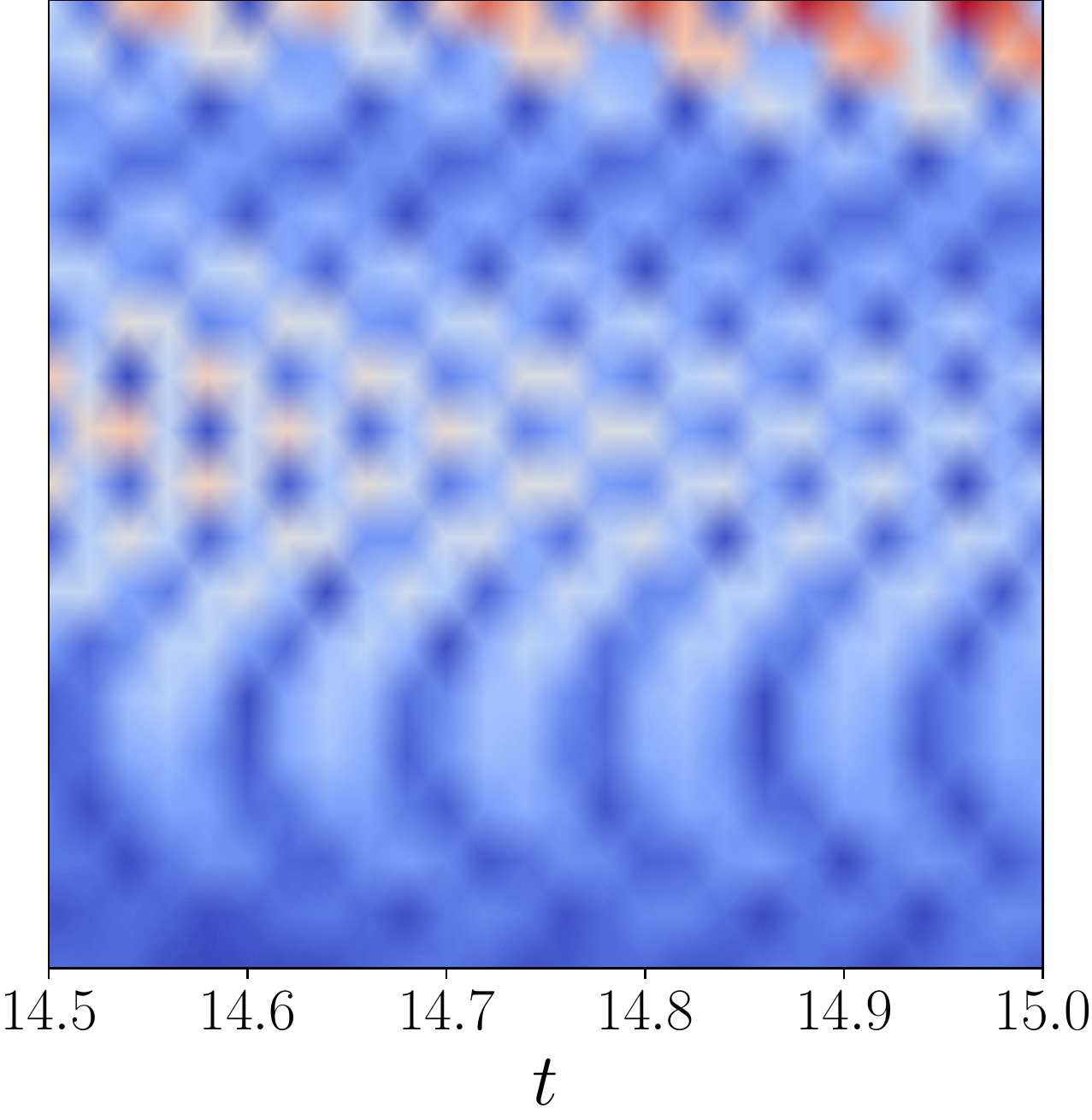}}
	\subfigure[MF 2-step LSTM]
	{\includegraphics[width=0.229\linewidth]{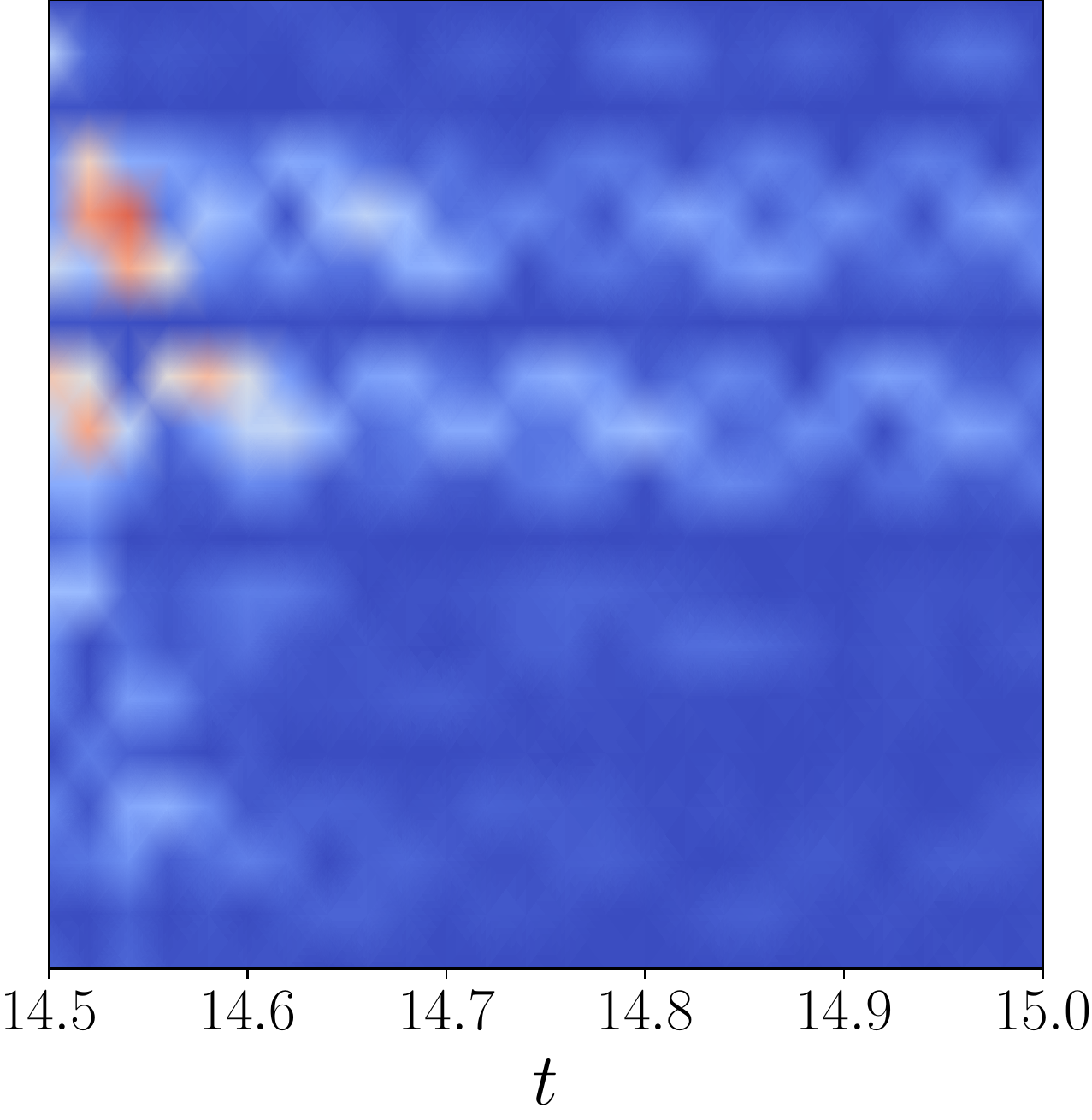}}
	\subfigure[MF 3-step LSTM]
	{\includegraphics[width=0.248\linewidth]{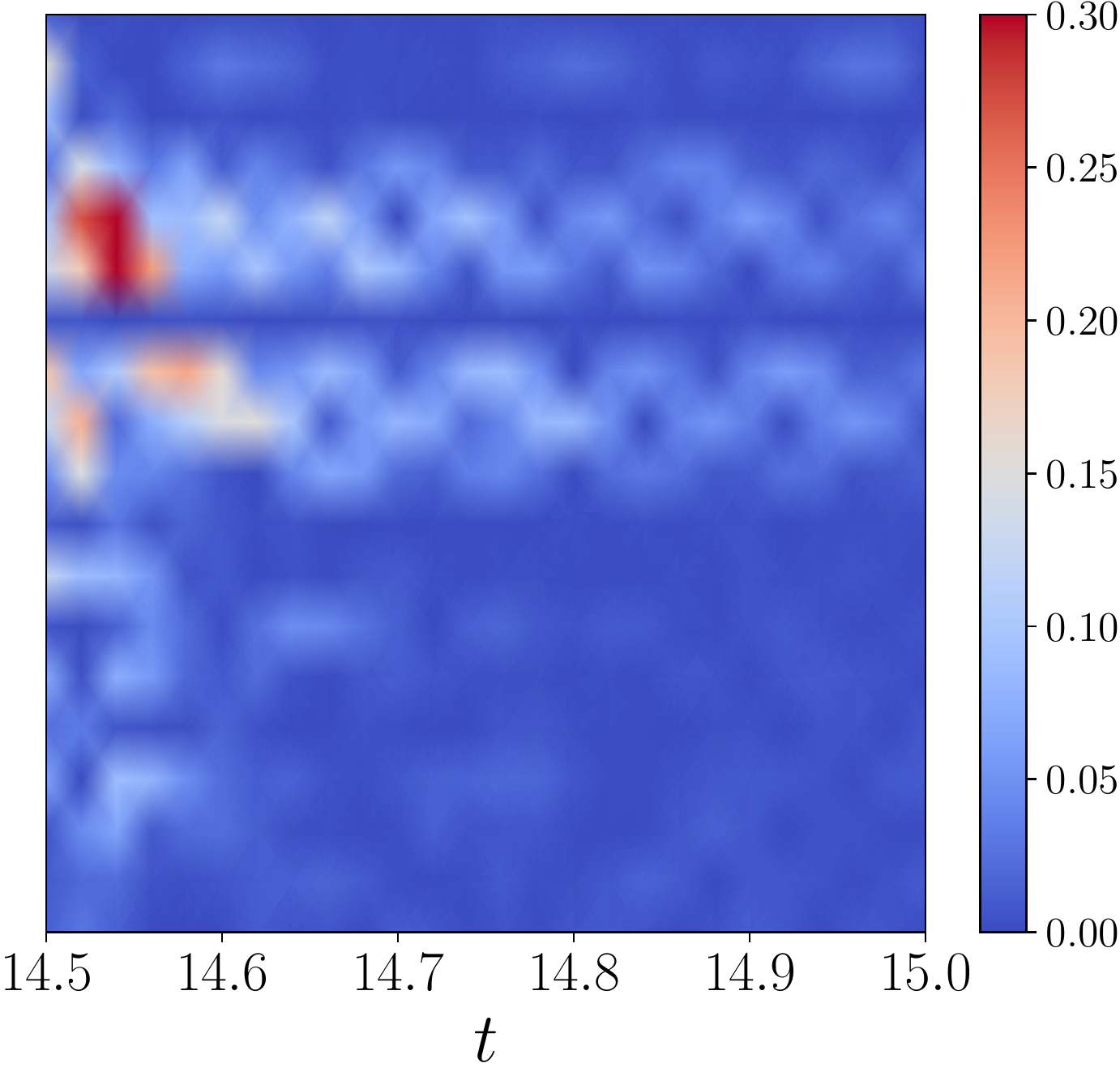}}
    \vspace{-5pt}
	\caption{Absolute values of the errors in different NN model predictions of the drag coefficient with $N^\mu_\texttt{HF} = 6$ HF training time series. }
    \label{fig: errors_drag}
\end{figure}

Similar to example (I), we train eight different NN models on the aforementioned data sets, and compare the proposed MF LSTM models with both the single-fidelity regressions and the MF feed-forward networks. 
In Tables \ref{tab: SF_NS} and \ref{tab: MF_NS}, we collect the prediction errors on the test set from all eight NN models. Our MF LSTM models achieve the best prediction accuracy. In particular, the 2-step and 3-step LSTM models outperform the others, as highlighted from the errors reported in the Table.  
The discrepancy between the HF solutions and network model predictions is displayed in Figs. \ref{fig: errors_lift} and \ref{fig: errors_drag}, the former for the lift with $N^\mu_\texttt{HF} = 10$ and the latter for the drag with $N^\mu_\texttt{HF} = 6$.

We notice that the feed-forward networks, which treat the time and parameter as equal input entries,  fail to approximate the oscillatory temporal patterns in the output quantities (see (a), (b), and (e) of Figs. \ref{fig: errors_lift} and \ref{fig: errors_drag}).
On the other hand, the single-fidelity LSTM regressions (subfigures (c) and (d)) manage to capture the oscillations; however they are still unable to guarantee high-quality predictions either due to low data quality or limited data availability.
The proposed MF LSTM models (subfigures (f), (g), and (h)), especially the 2-step and 3-step LSTM, achieve good predictive accuracy all over the time-parameter domain $\mathcal{P}\times [t_0,T]$.
\begin{table}[h!]
	\renewcommand{\arraystretch}{1.15}
	\centering%
		\caption{Test mean square errors (MSE) of the single-fidelity models in example (II).}
	\hspace*{-5pt}\begin{tabular}{c|c|cccc}
		\hline%
		\textbf{Output quantity }& \textbf{\#HF data} & \textbf{LF feed-forward} &  \textbf{HF feed-forward} &  \textbf{LF LSTM} &  \textbf{HF LSTM} \\
		\hline\hline
		Lift coefficient & 6 & $8.64 \times 10^{-2} $ & $1.12 \times 10^{-1}$ & $2.45\times10^{-1}$ & $1.12 \times 10^{-1}$\\
		Lift coefficient & 10 & $''$ & $8.63 \times 10^{-2}$ & $''$ & $1.92 \times 10^{-2}$ \\
		 Drag coefficient & 6 & $9.45\times 10^{-3}$ & $9.43 \times 10^{-3}$ & $8.77 \times 10^{-3} $& $1.33 \times 10^{-2}$\\
		 Drag coefficient & 10 & $''$ & $9.41 \times 10^{-3} $ & $''$ & $1.63 \times 10^{-3}$\\
		\hline
	\end{tabular}
	\label{tab: SF_NS}
	\vspace{-4pt}
\end{table}

\begin{table}[h!]
	\renewcommand{\arraystretch}{1.15}
	\setlength{\tabcolsep}{4.5pt}
	\centering%
	\caption{Test mean square errors (MSE) of the multi-fidelity models in example (II).}
	\hspace*{-6pt}\begin{tabular}{c|c|cccc}
		\hline%
		\textbf{Output quantity} & \textbf{\#HF data} & \textbf{3-step feed-forward} &  \textbf{Intermediate} &  \textbf{2-step LSTM} & \textbf{3-step LSTM}\\
		\hline\hline
		Lift coefficient & 6 & $8.67 \times 10^{-2}$ & $2.83 \times 10^{-2}$ & $\bm{1.14\times10^{-2}}$ & $2.24 \times 10^{-2}$\\
		Lift coefficient & 10 & $8.63 \times 10^{-2}$ & $6.64 \times 10^{-3}$ & ${1.25\times10^{-3}}$ & $\bm{1.22 \times 10^{-3}}$ \\
		 Drag coefficient & 6 & $9.44 \times 10^{-3}$ & $7.59 \times 10^{-3}$ & $\bm{1.58 \times 10^{-3}}$ & ${2.28 \times 10^{-3}}$\\
		 Drag coefficient & 10 & $9.41 \times 10^{-3}$ & $7.99\times 10^{-3}$& ${2.91 \times 10^{-4}}$ & $\bm{2.86 \times 10^{-4}}$\\
		\hline
	\end{tabular}
    \label{tab: MF_NS}
\end{table}

\begin{figure}[t]
	\centering
	\subfigure[LF LSTM -- lift]
	{\includegraphics[width=0.336\linewidth]{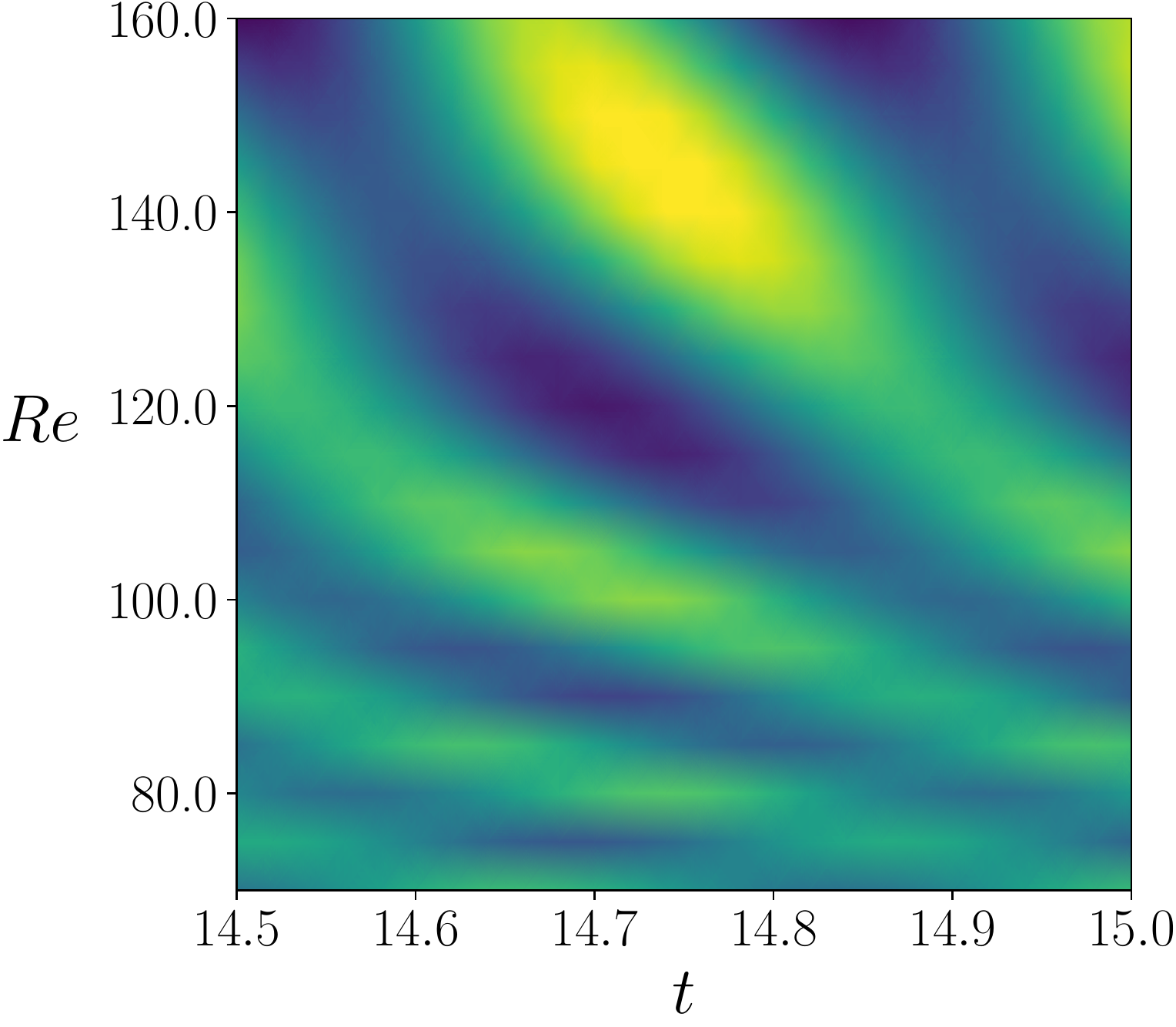}}
	\quad
	\subfigure[HF LSTM -- lift]
	{\includegraphics[width=0.281\linewidth]{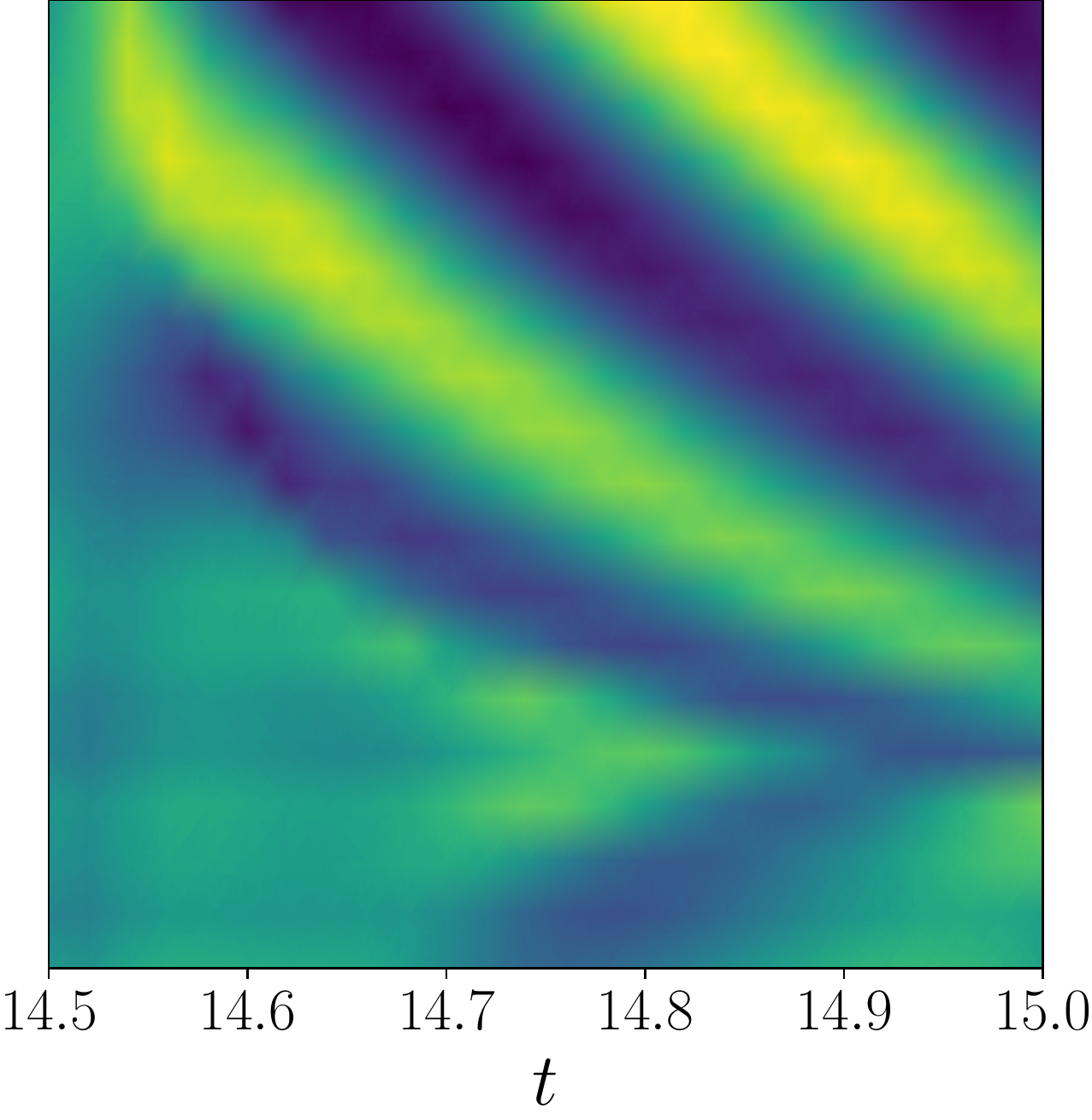}}
	\quad
	\subfigure[MF 3-step LSTM -- lift]
	{\includegraphics[width=0.308\linewidth]{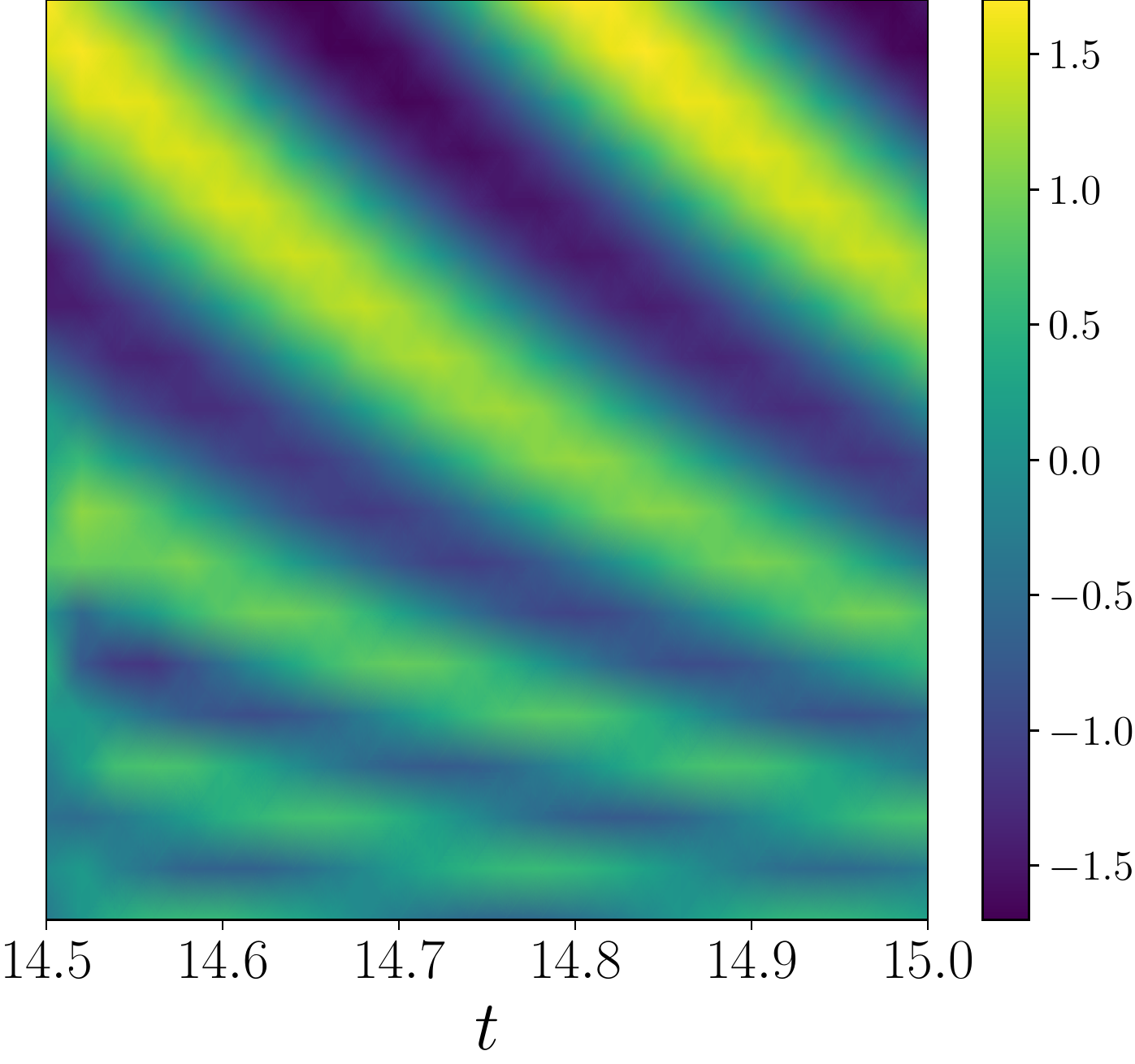}}

	\centering
	\subfigure[LF LSTM -- drag]
	{\includegraphics[width=0.336\linewidth]{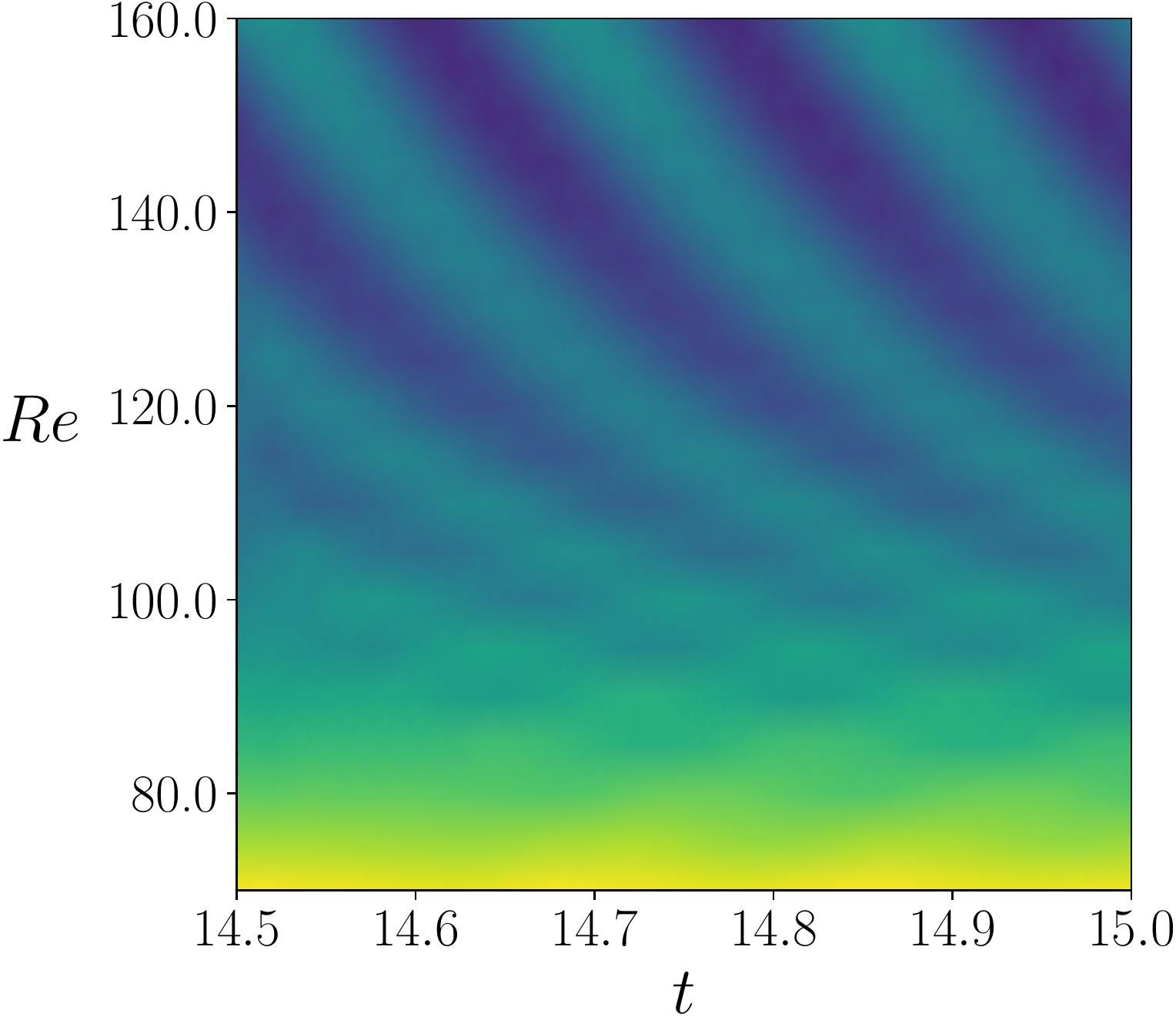}}
	\quad
	\subfigure[HF LSTM -- drag]
	{\includegraphics[width=0.281\linewidth]{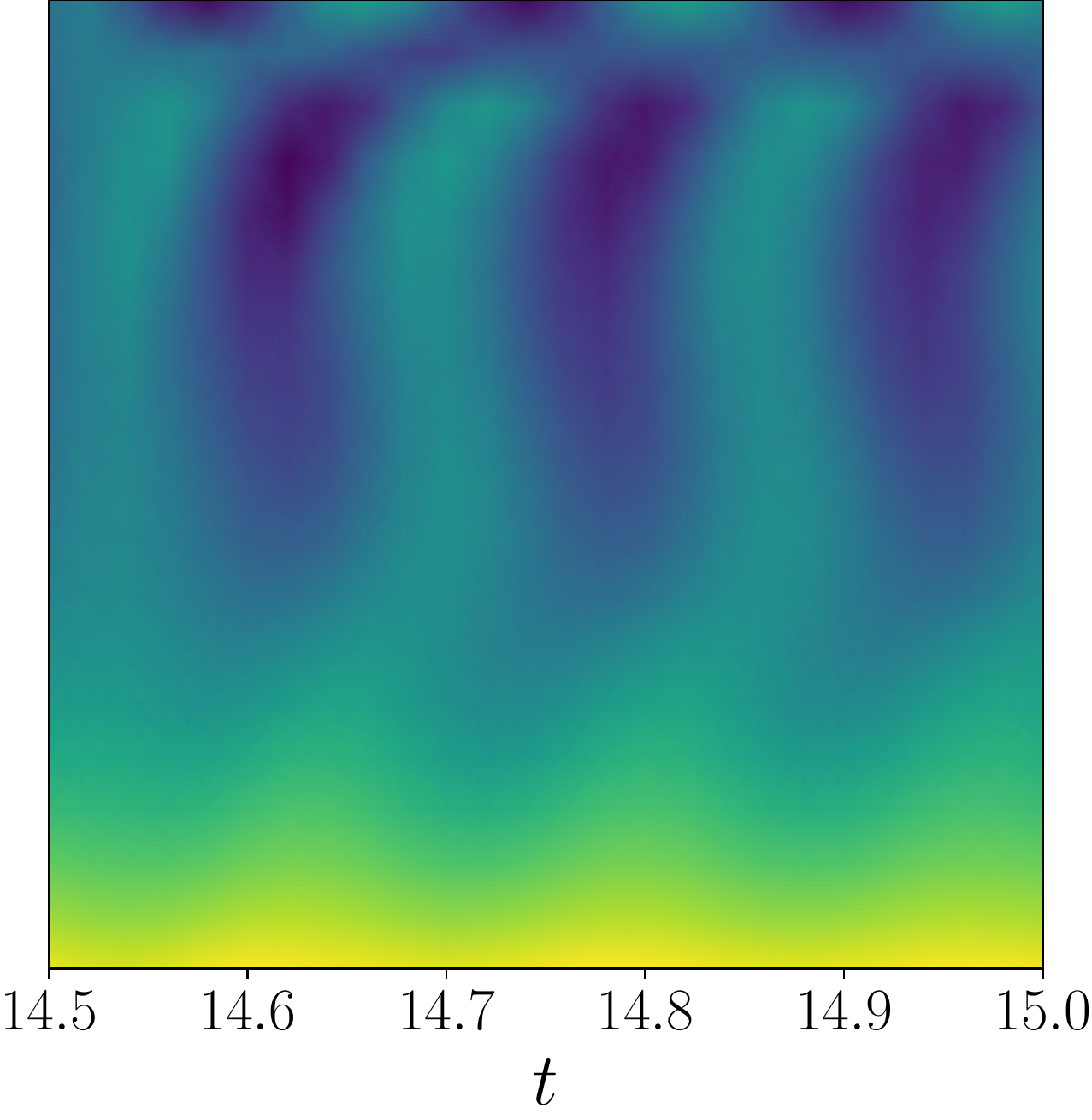}}
	\quad
	\subfigure[MF 2-step LSTM -- drag]
	{\includegraphics[width=0.308\linewidth]{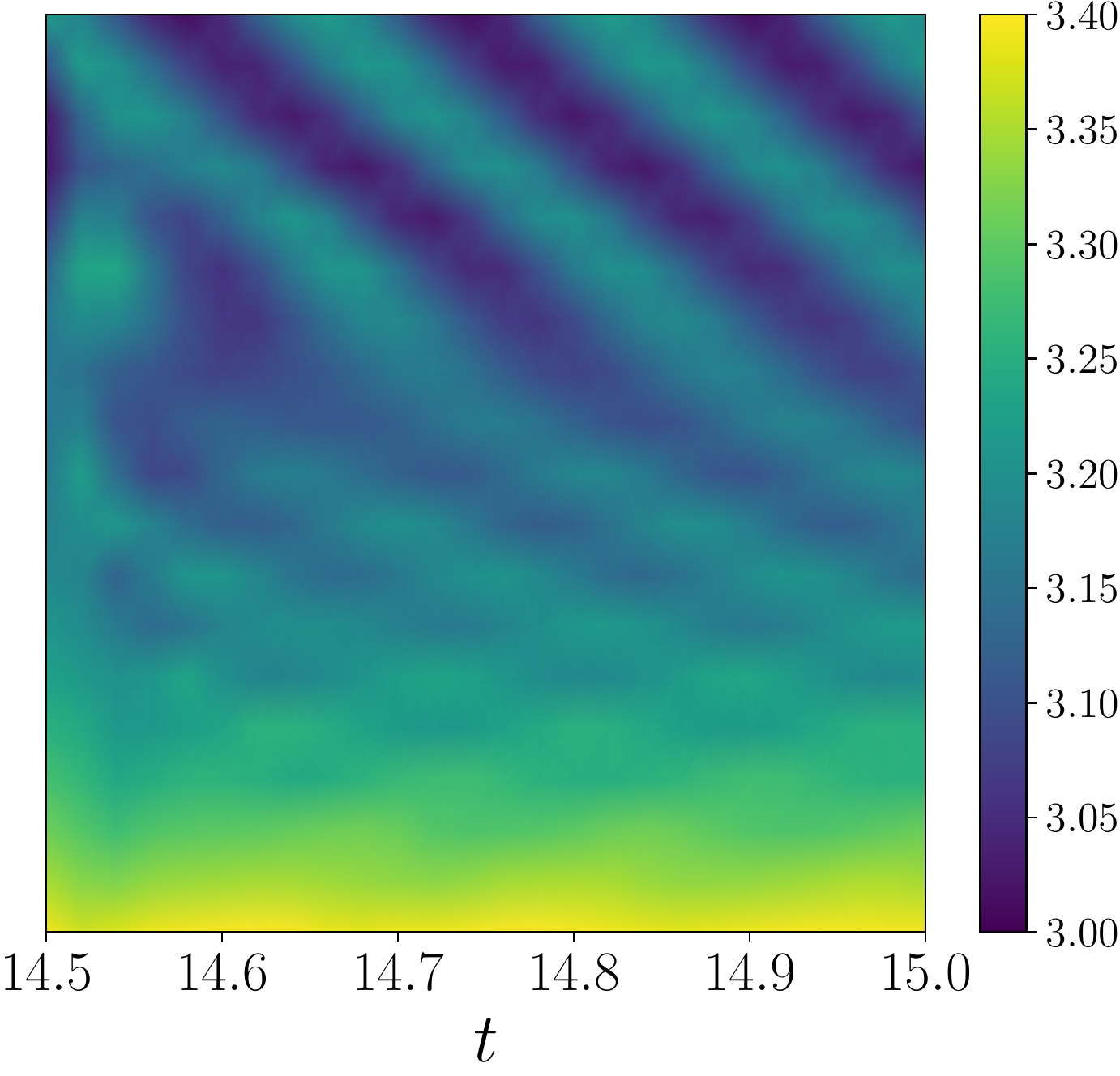}}

	\caption{Single-fidelity (LF and HF) and the best performing MF LSTM regressions for  the lift coefficient with $N^\mu_\texttt{HF}=10$ (top) and the drag coefficient with $N^\mu_\texttt{HF}=6$ (bottom).  }
    \label{fig: drag_lift_pred}
\end{figure}

In addition, we show in Fig.~\ref{fig: drag_lift_pred} the predictions from the single-fidelity (LF and HF) LSTM models and the best performing MF models, i.e., the 3-step (resp. 2-step) LSTM for the lift (resp. drag) coefficient with $N^\mu_\texttt{HF} = 10 $ (resp. $N^\mu_\texttt{HF} = 6 $).
Comparing these figures to the LF and HF solutions that generate the data (see Fig. \ref{fig: drag_lift_models}), it is reasonable that the LF LSTM reconstruction based on a large data set is in good agreement with the LF ground truth, yet still differs significantly from the HF solution because of the coarse spatio-temporal discretization, and that the HF LSTM reconstruction exhibits a better fit to the HF ground truth yet is still not ideal due to the limited data coverage. By taking advantage of both the large amount of LF data and the limited HF measurements, the MF results, however, achieve a close-to-HF accuracy in learning the temporal patterns with the LSTM layers.

\subsection{Results and discussions 2: Prediction forward in time}

So far, we have shown that the MF LSTM networks allow to accurately generalize/interpolate within the training region $\mathcal{P}\times [t_0,T]$. 
In this subsection we investigate the generalization capability of predicting forward in time, i.e., for the future states beyond the training range, which is known to be challenging for data-driven models.

\begin{figure}[b!]
    \vspace{10pt}
	\centering
	\subfigure[Lift coefficient ($N^\mu_\texttt{HF} = 10$)]
	{\includegraphics[width=0.47\linewidth]{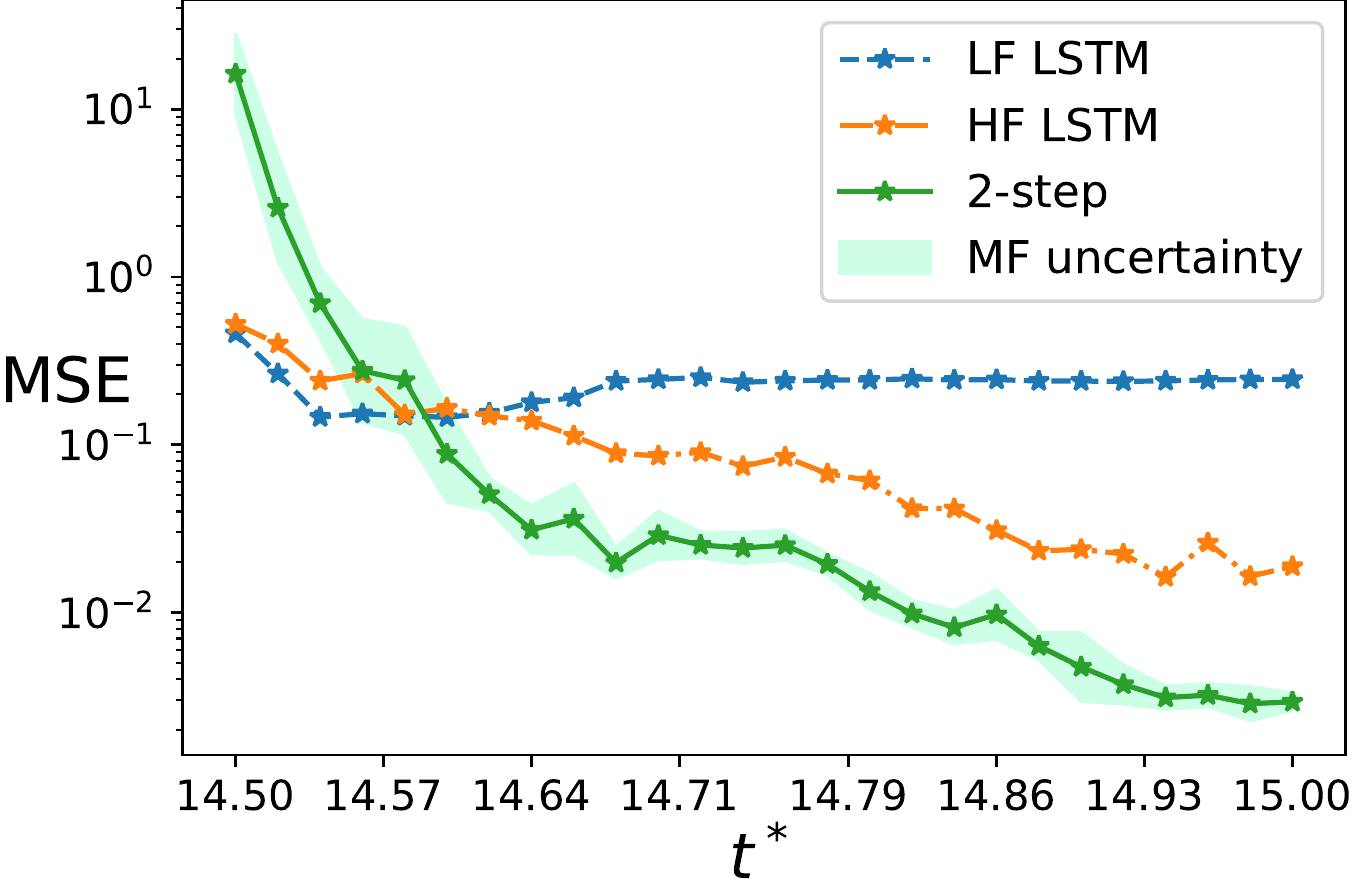}}
	\qquad
	\subfigure[Drag coefficient ($N^\mu_\texttt{HF} = 6$)]
	{\includegraphics[width=0.47\linewidth]{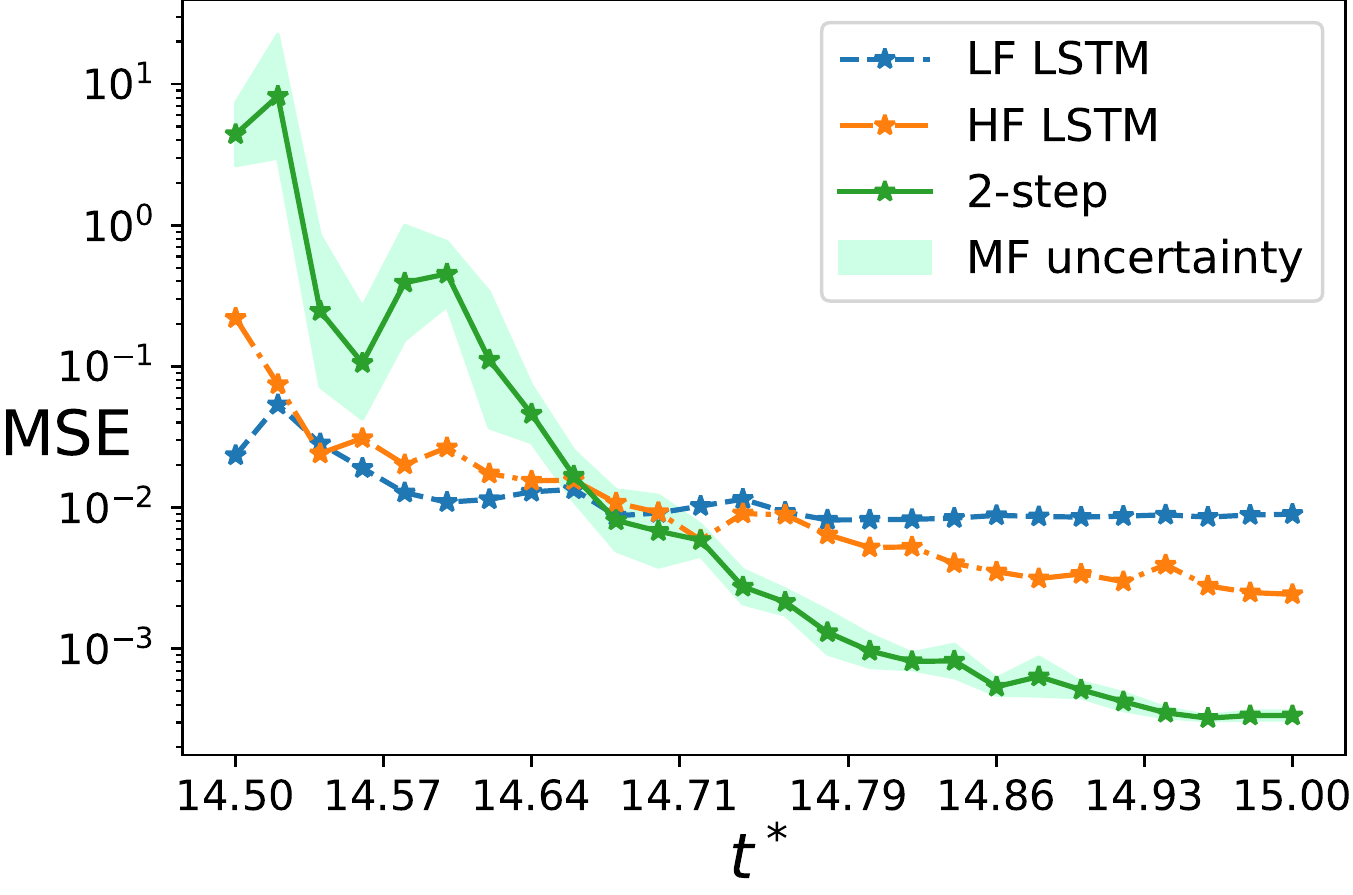}
	\label{fig: error_forward_drag}}
	\caption{{\color{black}{Test error on the whole time interval $[t_0, T] = [14.5,15.0]$ versus the final training time $t^*$. The testing time window $[t_0, T]$ remains fixed, while training time length ($t^*-t_0$) is increased as $t^*$ goes from $t_0$ to $T$. For each value of $t^*$, the second network of 2-step model is trained $20$ times with random initialization, and the corresponding model predictions are used to compute uncertainty bounds on the test set through an ensemble-based technique. The test error mean (green) and $\pm$ one standard deviation (shaded region) over the samples in the ensemble are shown, providing an uncertainty quantification for the MF predictions.}}}
    \label{fig: error_forward}
\end{figure}

\begin{figure}[b!]
	\centering
	\subfigure[Lift coefficient ($N^\mu_\texttt{HF} = 10$)]
	{\includegraphics[width=0.45\linewidth]{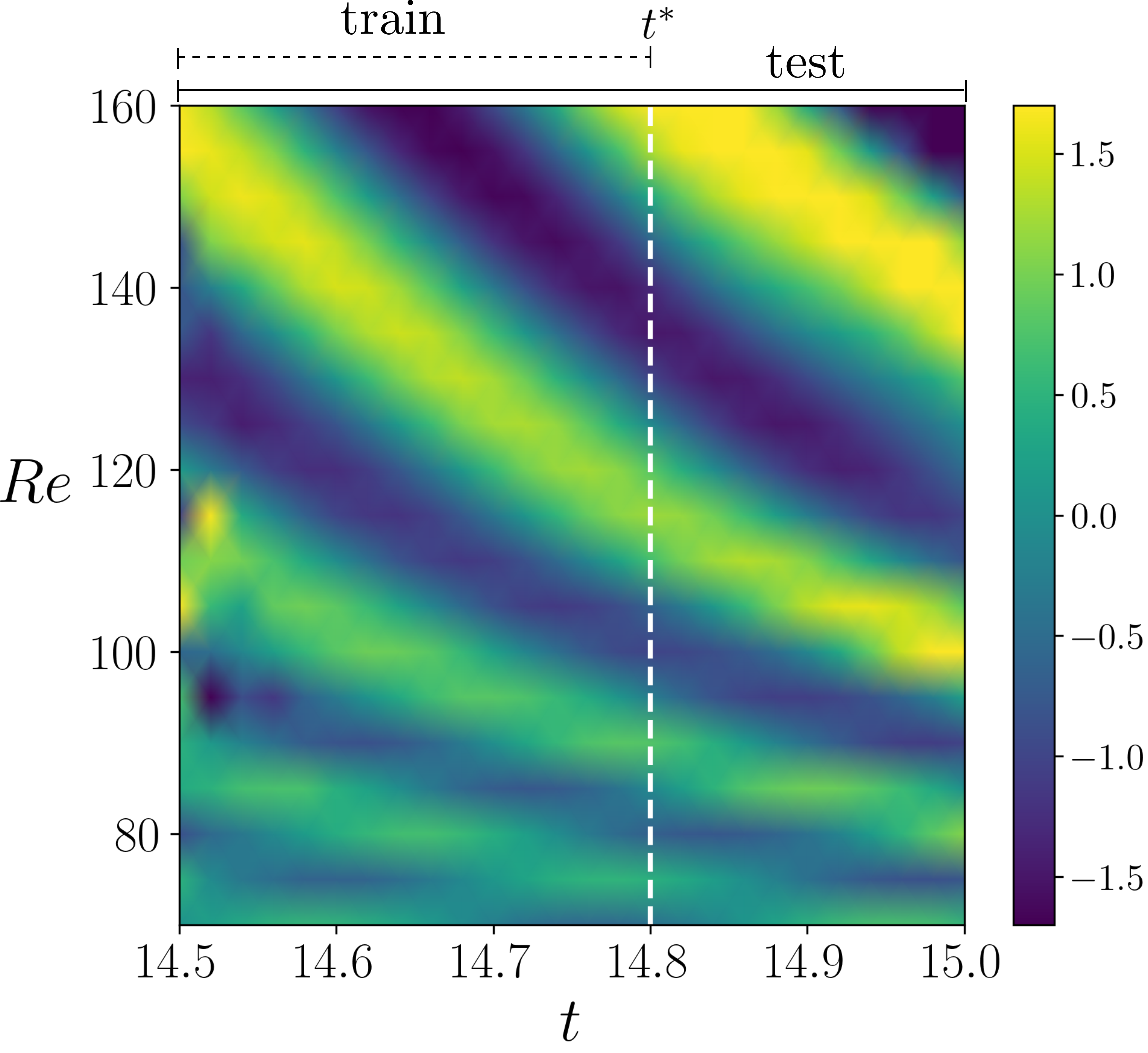}
	 \label{fig: lift_forward}}
	 \qquad
	\subfigure[Drag coefficient ($N^\mu_\texttt{HF} = 6$)]
	{\includegraphics[width=0.45\linewidth]{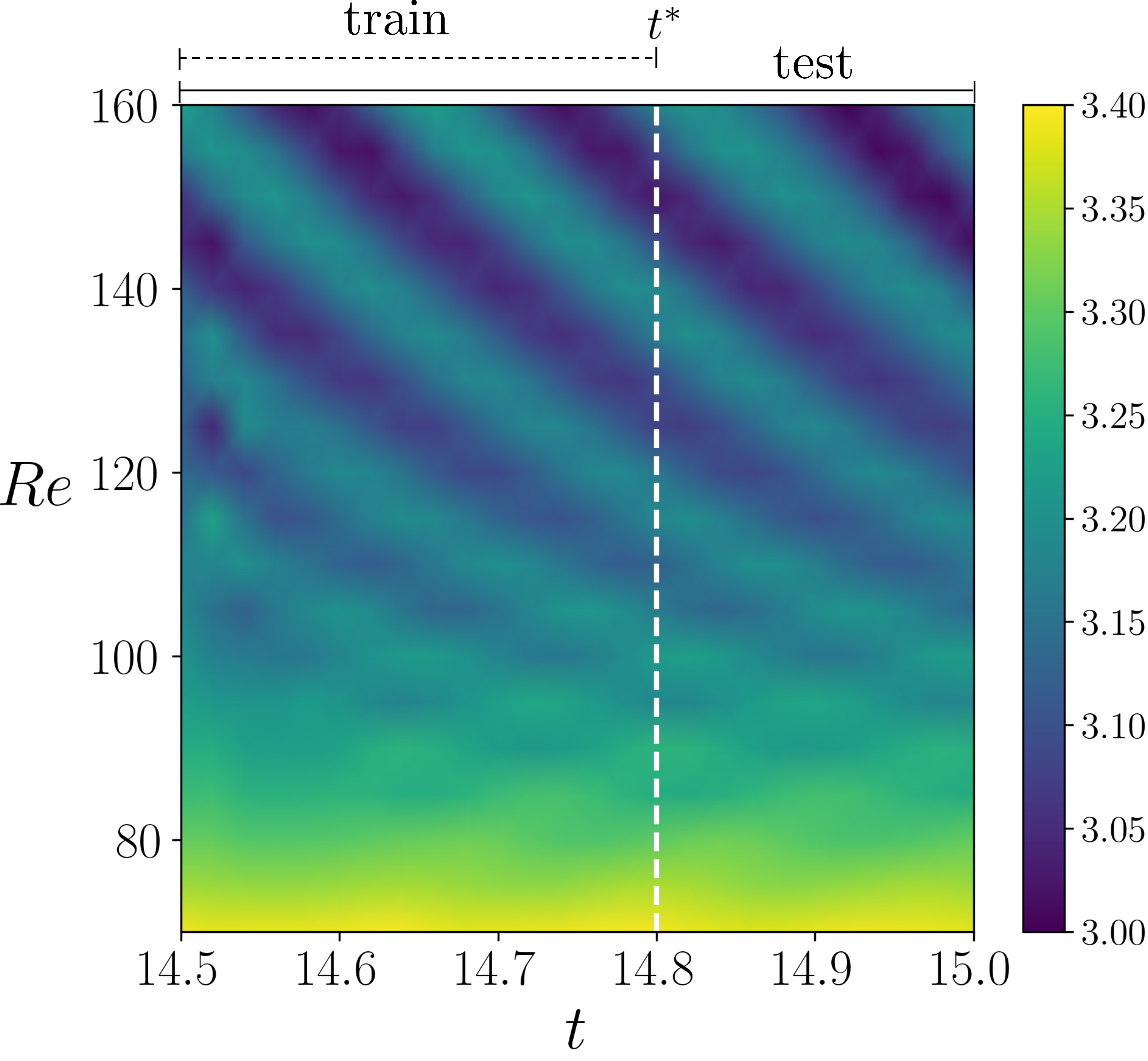}
	\label{fig: drag_forward}}
	\caption{\textit 2-step LSTM model prediction for the lift (left) and drag (right) coefficients. The model is trained with data up to $t^* = 14.80$ (white dashed line) and tested over the whole time interval $[14.5,15.0]$.}
	\label{fig: pred_forward}
\end{figure}

Starting with the same training data set as in the previous subsection, we progressively move backward the final training time $t^*$, from $t^* = T$ to $t^* = t_0$, while using the same instances of the Reynolds number. Though the training time interval is shortened, we still predict over the entire $[t_0,T]$, which includes both reconstruction and extrapolation in time. 
In Fig. \ref{fig: error_forward}, we report the prediction errors with respect to $t^*$ through different LSTM network models.
Given sufficient time steps to effectively merge LF and HF data, the MF LSTM model shows a consistent advantage over the single-fidelity models, both for the lift and drag coefficients (see Fig. \ref{fig: error_forward}).
For example, in the drag case (Fig. \ref{fig: error_forward_drag}), a smaller test error is achieved with the 2-step LSTM trained on half of the whole time interval ($t^* = 14.75$) than the single-fidelity LSTM, either HF or LF, with the full training set over time ($t^* = T = 15.0$). This suggests that we can improve the predictive accuracy by using MF LSTM with only more than half of the training data, which may enable a significant reduction in the cost of HF data generation. 
{\color{black}{Moreover, a description of uncertainty in the model predictions can be obtained with an ensemble-based technique \cite{dietterich2000ensemble}. Here, we train the last-step ($NN_\texttt{HF}$) of the 2-step LSTM model multiple times with the same data and randomly generated network initialization, and then compute the statistical moments of the corresponding samples of predictions. This allows to construct uncertainty bounds for MF approximation, see Fig. \ref{fig: error_forward}. The low predictive uncertainty further highlights the consistent robustness of the proposed 2-step MF model as the extrapolation time window varies. }}
In Fig. \ref{fig: pred_forward} we depict the approximation of both the lift and drag coefficients by the 2-step LSTM model with $t^* = 14.80$. The solution shows a very good agreement with the HF ground truth, even over the time interval $[t^*,T]$ where we have no training data at all -- neither HF nor LF. 

Results have highlighted the effectiveness of the proposed MF LSTM networks in time-parameter-dependent surrogate modeling. For output quantities that present an oscillatory behaviour, the proposed models allow to learn the temporal pattern and predict the evolution of future states. Good generalization has been seen not only in the parameter ranges where few HF data are available, but also at future time instances outside the coverage of training data.

\section{Numerical example (III): Lotka-Volterra system}
\label{section: LV}
In this section we consider a Lotka–Volterra system that describes a nonlinear, three-species prey-predator interaction:
\begin{equation}
    \begin{cases} 
    \frac{dy_1}{dt}(t) = y_1(t)(\mu - 0.1 y_1 (t) -0.5 y_2(t) - 0.5  y_3 (t)), & t \in (0,T),\\ 
    \frac{dy_2}{dt}(t) = y_2(t)(-\mu + 0.5 y_1 (t) - 0.3 y_3 (t)), & t \in (0,T),\\ 
    \frac{dy_3}{dt}(t) = y_3(t)(-\mu + 0.2 y_1 (t) + 0.5 y_2 (t)), & t \in (0,T),\\ 
    y_i(0) = 0.5, & \forall i = 1,2,3.
    \end{cases}
\label{eq: LV}
\end{equation}

Here we employ the proposed MF LSTM models to approximate vector-valued output quantities $\vb{y}(t;\mu)=[y_1(t;\mu), y_2(t;\mu), y_3(t;\mu)]^\text{T} \in \mathbb{R}^3$, i.e., the number of individuals in each population/species. In this MF regression task, we aim at estimating the system solution $\vb{y}$ up to time $T = 15.0$, as the parameter $\mu$ changes in $\mathcal{P}=[1.0,3.0]$. 
Such a regression task is more challenging than those in the previous examples because the vector-valued quantity of interest $\vb{y}$ exhibits aperiodic oscillations.  Moreover, the parametric variation in $\mu$ induces different aperiodic patterns in the amplitude and frequency of the oscillations. In this case, we only let the HF data cover a limited time window until $T_\texttt{HF} < T$, and intend to infer the aperiodic evolution of $\vb{y}$ from the LF data that cover the entire domain $\mathcal{D} = \mathcal{P} \times [0,T]$.

\subsection{Multi-fidelity setting}
\begin{figure}[h!]
	\centering
	\subfigure[$\mu = 1.0$]
	{\includegraphics[width=0.455\linewidth]{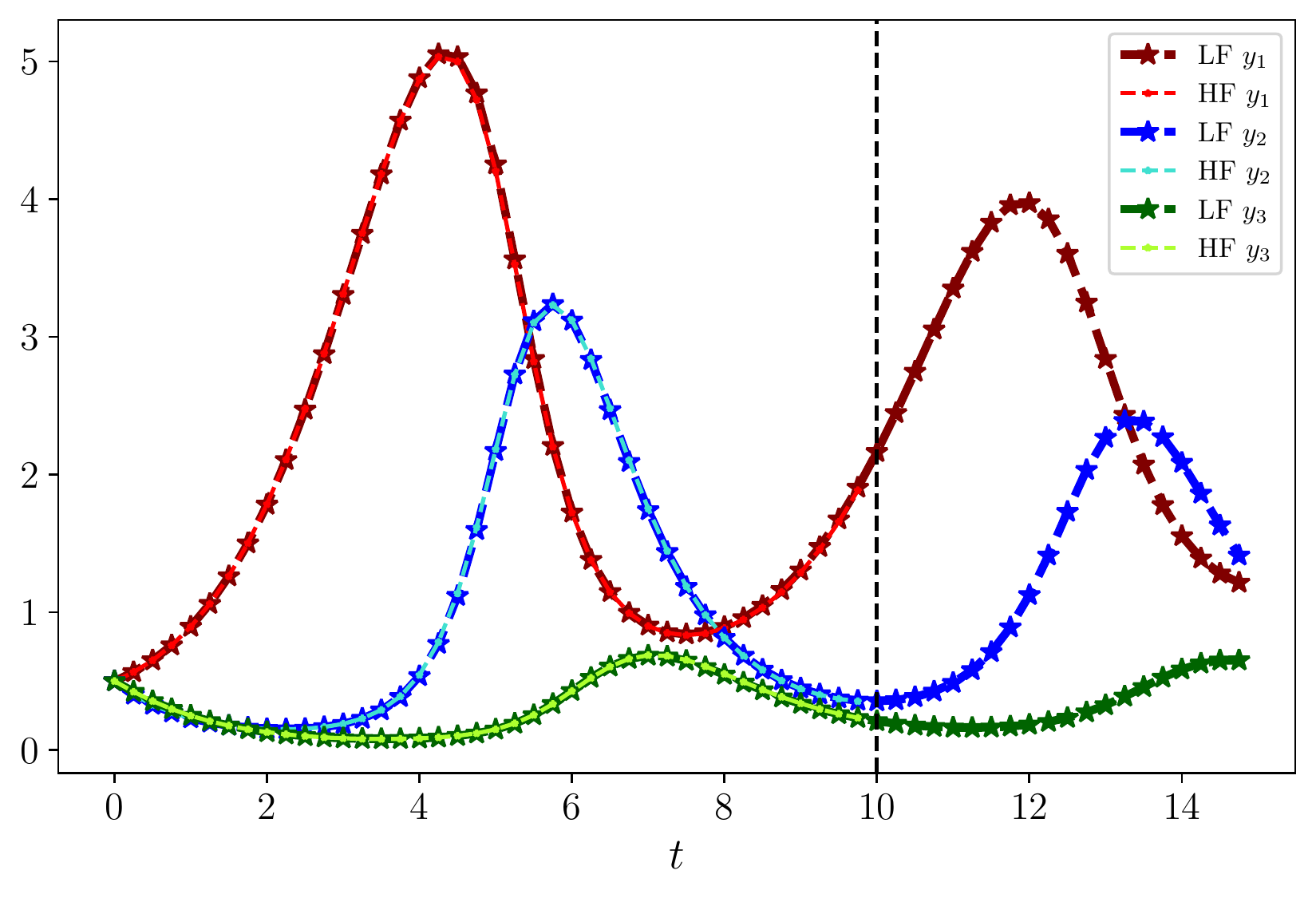}}
	\quad
	\subfigure[$\mu = 3.0$]
	{\includegraphics[width=0.465\linewidth]{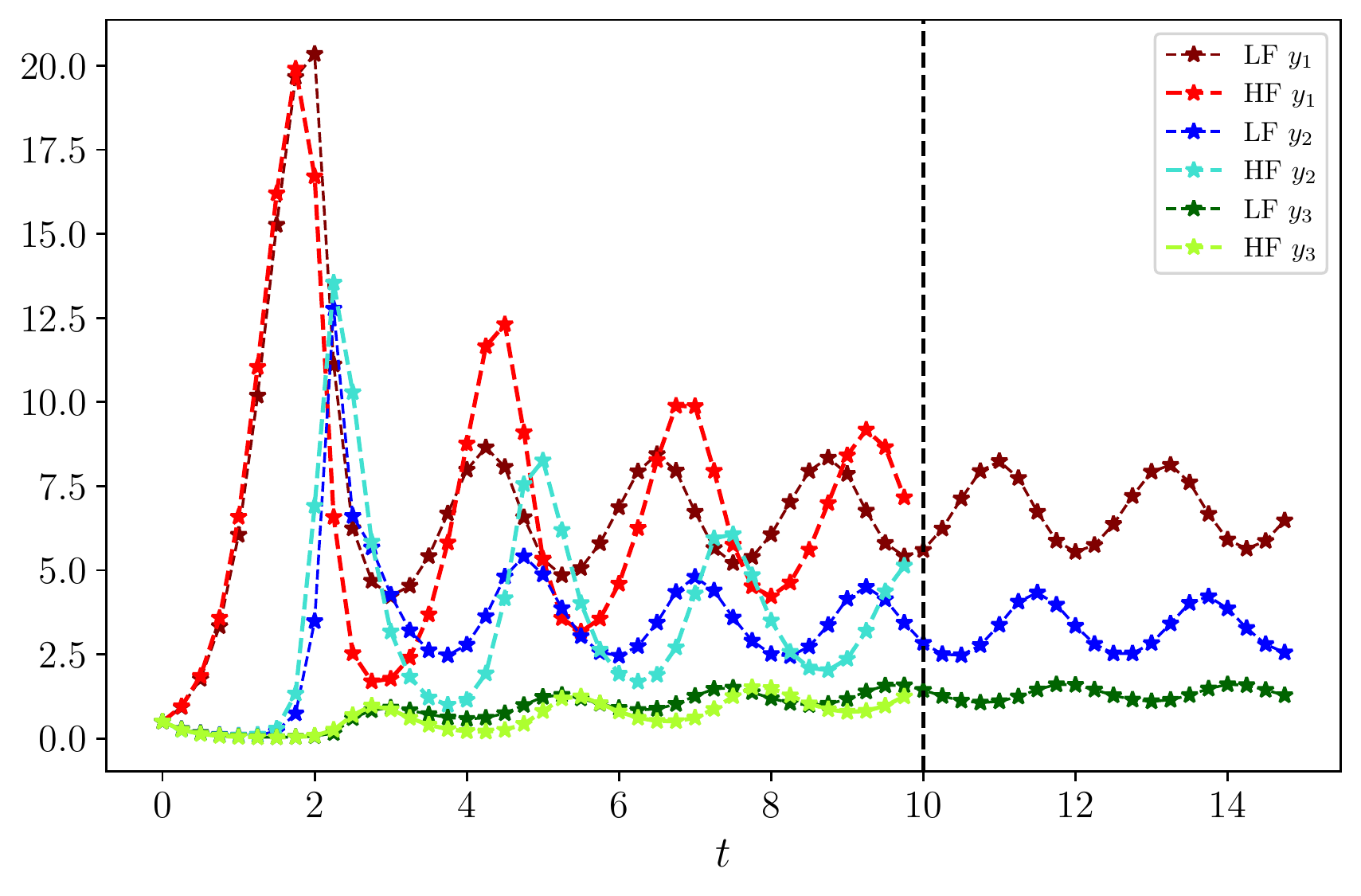}}
	\caption{HF and LF training data for two parameter instances. The parameter $\mu$ regulates both the amplitude and frequency of the oscillations of $\vb{y}$.  As $\mu$ grows, the amplitude and frequency increase, so does the discrepancy between the two data levels.}
    \label{fig: LV train data}
\end{figure}

The training data are generated by the second-order Runge-Kutta (RK2) scheme with time steps $\Delta t_\texttt{LF} = 0.25$ and $\Delta t_\texttt{HF} = 0.0025$ for the LF and HF models, respectively. 
The LF data are taken at $N_\texttt{LF}^\mu = 20$ uniform $\mu$-values over $\mathcal{P} = [1.0,3.0]$. For each $\mu$-instance, we collect the LF training data along with the time integration with $\Delta t_\texttt{LF}=0.25$ up to the final time $T_\texttt{LF} = T= 15$. Similarly, the HF data are taken at $N_\texttt{HF}^\mu = 10$ uniform $\mu$-values. Though the HF time integration is evaluated with $\Delta t_\texttt{HF}= 0.0025$, we only collect the HF training data every $\Delta t_\texttt{LF}=0.25$ to be consistent with the LF data, and stop the data collection at $t = T_\texttt{HF}=10$. Thus, no HF data are available in the time interval $[10.0,15.0]$. For the LSTM training, we group the data into batch subsequences of length $K = 25$. 
Fig. \ref{fig: LV train data} shows the training data for $\mu = 1.0$ and $3.0$, the minimum and maximum $\mu$-values. We note that the frequency of output oscillations increases with $\mu$. This implies that the LF coarse time discretization impacts the data quality more significantly for a larger value of $\mu$. In fact, at $\mu = 3.0$ we spot a much larger discrepancy between the HF and LF data in comparison to the case of $\mu = 1.0$. 
We generate the test set similarly to the HF training data, but include $N_\text{test} = 30$ $\mu$-values in $\mathcal{P}$ and reach the final time  $T = 15$.

\subsection{Results and discussions}
\begin{figure}[h!]
	\centering
	\subfigure[Single-fidelity LF LSTM]
	{\includegraphics[width=0.325\linewidth]{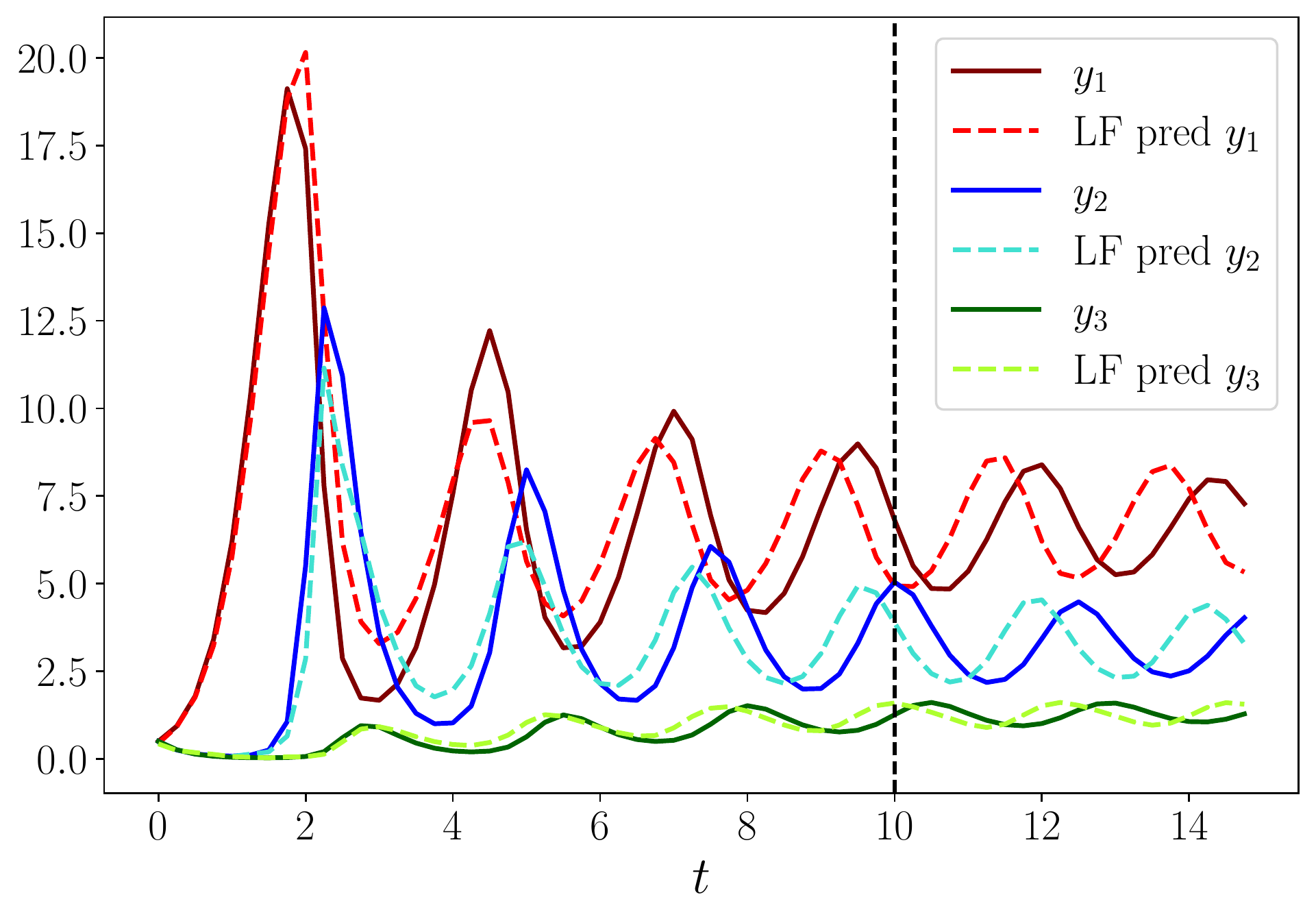}}
	\subfigure[Single-fidelity HF LSTM]
	{\includegraphics[width=0.325\linewidth]{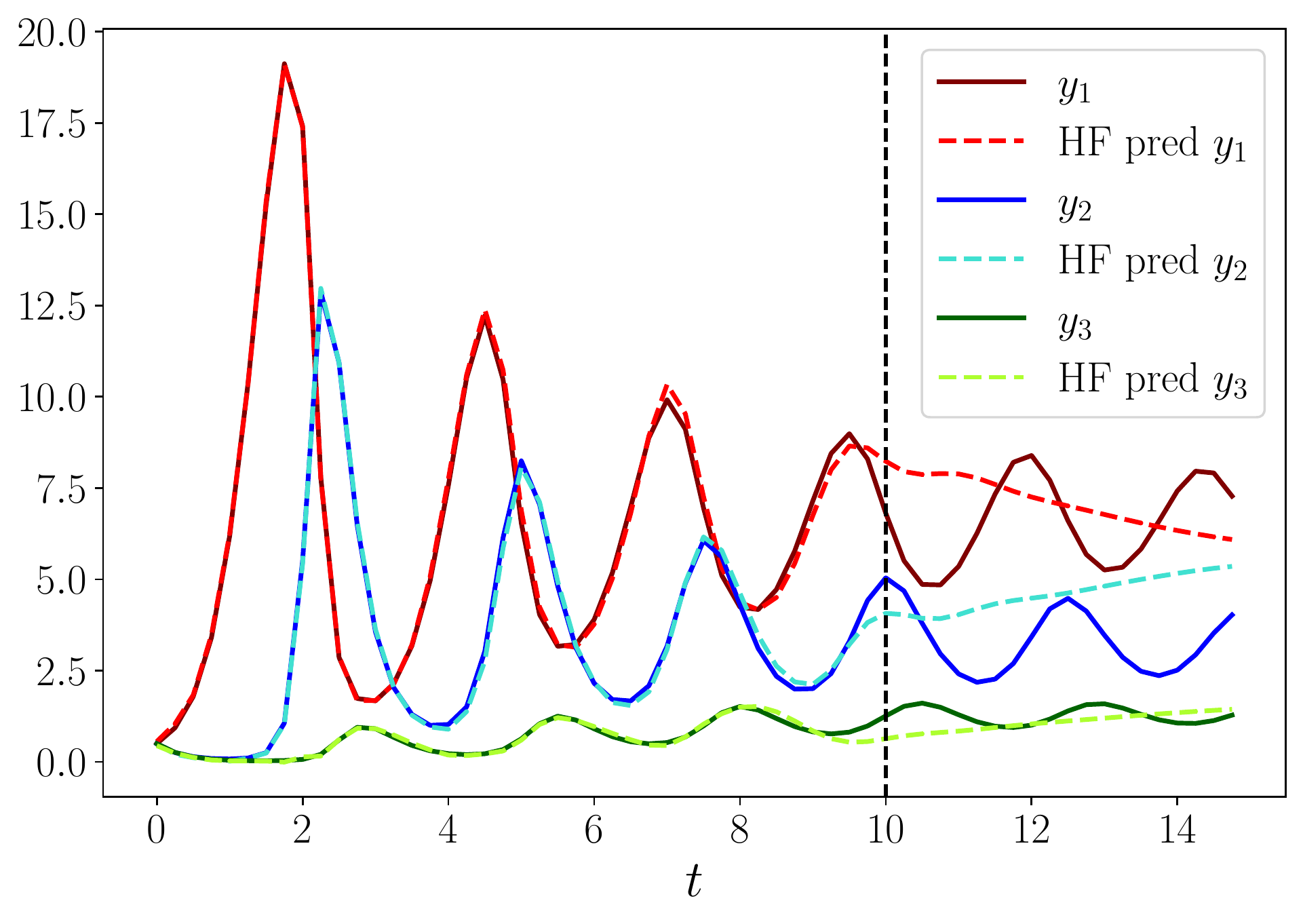}}
    \subfigure[Multi-fidelity 2-step LSTM]
	{\includegraphics[width=0.325\linewidth]{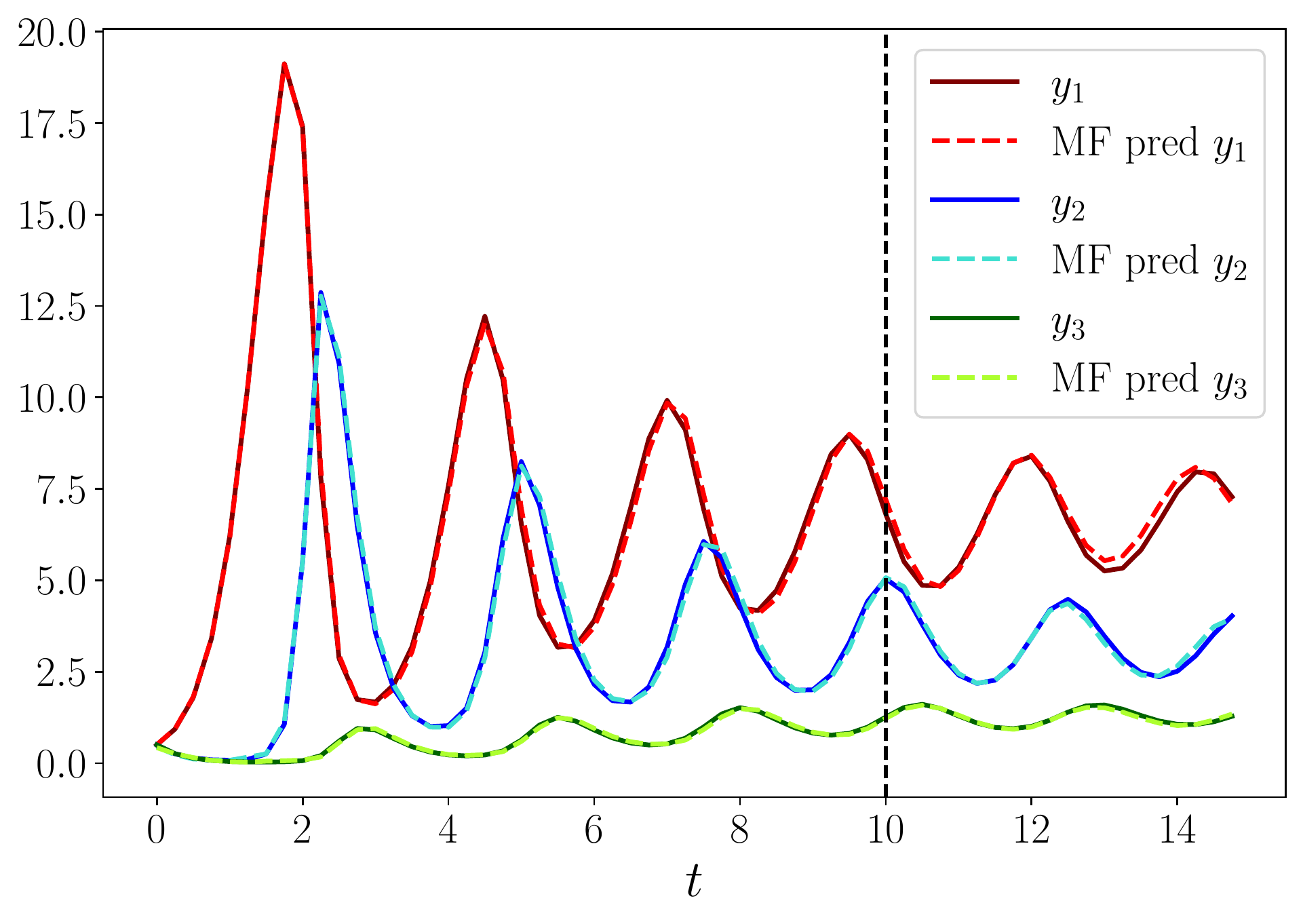}}
	\caption{Model predictions by the single- and multi-fidelity LSTM networks compared to the HF solution ($\mu = 2.93$). The black dashed line at $T_\texttt{HF} = 10$ indicates the maximum time instant of the HF training data, i.e., no HF information is available over $t\in [T_\texttt{HF},T] =[10,15]$.}
    \label{fig: LSTM LV}
\end{figure}

Here we compare the single-fidelity and MF LSTM models on the given data sets. For the sake of brevity, we only present the results by the 2-step architecture among the MF LSTM networks.
Table \ref{tab: LV tab} collects the prediction MSE on the test set that covers the whole domain $\mathcal{D}$, and we can clearly observe that the MF model outperforms the single-fidelity ones. Furthermore, Fig. \ref{fig: LSTM LV} shows the comparison between the model predictions and the HF reference solution at a testing parameter instance $\mu = 2.93$.

The LF LSTM model fails to guarantee a low test error due to the poor data quality stemming from the coarse approximation. 
The HF LSTM model achieves a good accuracy in the range $[0,T_\texttt{HF}]=[0,10]$, but cannot generalize to $t > T_\texttt{HF}$ where no HF data are available for any $\mu$-value. This implies that even though HF LSTM can well interpolate over the parameter domain, it fails to extrapolate over time.  This is understandable as the limited HF information is not sufficient for capturing the aperiodic oscillations or supporting the predictions of future states. 
On the other hand, our 2-step LSTM model performs well: it manages to both interpolate for the parameter $\mu$ and predict forward in time. Here the full-domain LF data coverage over $\mathcal{D}$ plays a critical role in enabling the MF time extrapolation, as the MF LSTM captures the aperiodic oscillations by well approximating the correlation between the two fidelity levels.

\begin{table}[h!]
	\renewcommand{\arraystretch}{1.3}
	\centering%
	\caption{Test errors for the single- and multi-fidelity models in example (III). }	
		\begin{tabular}{c|ccc}
	\hline%
		{\textbf{Model}} &  {LF LSTM} &  {HF LSTM} &  {MF 2-step LSTM}\\
		\hline\hline
		{\textbf{Test MSE}} &  $0.354$ &  $0.853$ &  $\mathbf{0.043}$\\
		\hline
	\end{tabular}
	\label{tab: LV tab}
\end{table}

{\color{black}{
\section{Extension to more fidelity levels}
\label{section: Extension}
Although this work has been presented with bi-fidelity data, the proposed models in Sect.~\ref{sect: mf_models} can be extended to more than two fidelity levels. Let us consider $M$ datasets $\{\mathcal{T}_m = \{\vb{x}_m, \vb{y}_m\}\}_{m=1}^M$ given by a hierarchy of fidelity levels, sorted by increased accuracy of the data.

One can upgrade the proposed \textit{multi-level} models in the following ways:
\begin{itemize}
    \item \textit{Series multi-level}. The number of ``steps'' is increased by concatenating multiple NNs -- one for each fidelity level. This leads to a sequence of networks $\{NN_m\}_{m=1}^M$. The first network is fed with $\mathcal{T}_1$ to learn $\vb{f}_1$, and then, sequentially, the \textit{m}-th network $NN_m$ is trained to approximate $\vb{f}_m$ between the inputs $[\vb{x}_m, \vb{f}_{m-1}(\vb{x}_m), \ldots \vb{f}_1(\vb{x}_m)]^\text{T}$ (all of them or a subset) and the outputs $\vb{y}_m$.
    \item \textit{Parallel multi-level}. Each network of $\{NN_m\}_{m=1}^{M-1}$ is trained independently on $\mathcal{T}_m$ to approximate $\vb{f}_m$. The last network $NN_M$ is fed with the inputs $[\vb{x}_M, \vb{f}_{M-1}(\vb{x}_M), \ldots \vb{f}_1(\vb{x}_M)]^\text{T}$ (including all the other networks' outputs) and the outputs $\vb{y}_M$ to learn the HF function $\vb{f}_M$.
    \item A combination that mixes the \emph{series} and \emph{parallel} settings of networks.
\end{itemize}
The \textit{series} approach should be advantageous over the \textit{parallel} approach in the case where data come from sources of similar nature (e.g., refinement of numerical discretization) and each fidelity level adds further information to the previous. Thus, it is meaningful to sequentially incorporate the features from previous fidelity levels into the prediction of the current level. On the other hand, the \textit{parallel} strategy is more suitable for the case with multiple LF data sources that are not strongly correlated but individually provide useful information for the final HF prediction, so they are all together included in the inputs of $NN_M$. It is worth noting that the \textit{parallel} strategy is expected to present better computational flexibility, because the first $M-1$ NNs can be trained independently and/or in parallel.

The ``intermediate'' model can be upgraded beyond bi-fidelity by locating the multiple levels of LF outputs in sequential hidden layers. 
Despite a similarity with the \textit{series multi-level} extension, the upgraded ``intermediate'' model is determined by a single NN training whose loss function is a convex combination of all the losses of individual fidelity levels with coefficients $\{\alpha_m \in \mathbb{R}^+: 1\leq m \leq M, ~\sum_{m=1}^{M}\alpha_m =1\}$, similar to \eqref{eq: intermediate_loss} in Sect.~\ref{sect: mf_models}. To regulate the contribution of each loss term, these coefficients can either be set manually or determined via the hyperparameter optimization.

It is sometimes the case that we may update the data for one of the fidelity levels or insert an additional level to an already trained MF model. The \textit{parallel multi-level} approach hence only requires the network training for the updated/added data levels and the retraining of the last network, while freezing the other networks. With the ``intermediate'' setting or the \textit{series multi-level} extension, however, the whole model must be retrained.
}}

\section{Concluding remarks}
\label{section: conclusions}
In this work we present several novel methods for the estimation of time-parameter-dependent quantities of interest using multi-fidelity techniques with LSTM neural networks. The proposed techniques are shown to be advantageous over both the single-fidelity neural networks and the multi-fidelity ones without LSTM units, especially when approximating time-dependency. The multi-fidelity LSTM strategy enables accurate estimation of the target quantities' time evolution at a reasonable computational cost, as it leverages many low-fidelity, easily attainable data and only requires a limited number of high-fidelity, expensive data. The non-intrusive nature of the proposed models guarantees a wide applicability, which has been exemplified by a diverse collection of engineering applications. On the other hand, the multi-fidelity LSTM models show excellent predictive capabilities and generalization performance both in time and over parametric variation, which is typically deemed challenging for data-driven, non-intrusive surrogate models.

The proposed multi-fidelity LSTM models can be straightforwardly extended to more than two fidelity levels, as well as to the involvement of multiple physical systems, for which we can collect multi-fidelity data from several interacting systems governing the quantities of interest. Another promising, but more challenging, future application is the full-field approximation of PDE solutions. In this regard, an viable strategy is to consider the time-parameter-dependent coefficients of a low-dimensional reduced basis (e.g., through proper orthogonal decomposition on a set of HF snapshot solutions) as output quantities for the multi-fidelity LSTM models.

\section*{Acknowledgment}
The second author is financially supported by Sectorplan Bèta (the Netherlands) under the focus area \emph{Mathematics of Computational Science}. The third author acknowledges the support from Fondazione Cariplo under the Grant n. 2019-4608.  The authors would like to express their appreciation to Dr. Stefania Fresca for the fruitful discussions and for her help with numerical implementations.

\setcounter{table}{0}
\renewcommand{\thetable}{A.\arabic{table}}

\newpage
\section*{Appendix: hyperparameter summary}
Features of an NN model that cannot be optimized during the training process are called \emph{hyperparameters}. They are either related to the network structure (e.g., the numbers of layers and nodes) or associated with the training (e.g., the learning rate). The performance of an NN highly relies on the hyperparameter choices, and thus it is important to find an optimal set of hyperparameter values. 
In this work, we employ a Bayesian hyperparameter optimization technique that minimizes an objective function $\mathcal{O}$ in a multi-dimensional domain $\Lambda = \Lambda_1 \times \Lambda_2 \times \cdots \times \Lambda_N$ of all $N$ hyperparameters. Here $\Lambda_i$ represents the range of the $i$-th hyperparameter value. $\mathcal{O}$ is defined as the cross-validation error on the training set, because generating a HF validation set would be expensive and data reuse thus becomes necessary. 
Specifically, we use the Python package Hyperopt \cite{bergstra2013hyperopt} (see \cite{bergstra2011algorithms, hyperopt} for more details). 

{\color{black}{
To ensure a fair comparison among different models, the same choice of domain $\Lambda$ and objective function $\mathcal{O}$ is considered in the hyperparameter optimization for each model, unless otherwise specified. The hyperparameters are estimated by a common Bayesian optimization method \cite{bergstra2013hyperopt, bergstra2011algorithms}, which balances between the minimization of error $\mathcal{O}$ and the exploration of domain $\Lambda$. In the presented numerical examples, the performance of different models is assessed with optimized hyperparameter values.
}}

Tables \ref{tab: HP_Pulse}, \ref{tab: HP_Lift}, \ref{tab: HP_Drag} and \ref{tab: HP_LV} collect the optimized HPs in all the numerical examples.
For notation, we let $\eta$ denote the learning rate, and let `depth'  and `width' represent the numbers of layers and nodes, respectively. For the multi-level models we enumerate the level steps with the variable `step'. Although the intermediate model consists of a single NN, we use the variable 'step' to denote the `LF' (resp. `HF') portion of the network, i.e., the first (resp. second) part of the network dedicated to the estimation of the LF (resp. HF) output (see Fig. \ref{subfig:Intermediate} in Section \ref{section: NNs}). We recall that $\alpha$ is the coefficient regulating the contributions of the two fidelity levels to the loss function in the intermediate LSTM model. Moreover, all the NNs use the hyperbolic tangent activation function.

	\begin{table}[h!]
	\setlength{\tabcolsep}{4.5pt}
 		\caption{Optimized HP values in example (I) (Section \ref{section: Pulse}).}
		\renewcommand{\arraystretch}{1.3}
		\centering
		\begin{tabular}{c|c|cccccc}
			\hline 
    		Model &  Step & \begin{tabular}{@{}c@{}}Depth $\times$ width \vspace{-4pt} \\ LSTM\end{tabular}   & \begin{tabular}{@{}c@{}}Depth $\times$ width \vspace{-4pt} \\  Dense \end{tabular}   & Optimizer & $\eta$ & Batch size & $\alpha$ \\
			\hline 
			LF feed-forward & - & - & 4 $\times$ 61 & Adam & 2.54$\times 10^{-2}$ & 13 & -\\ 
			HF feed-forward & - & - & 3 $\times$ 52 & Adam & 1.98$\times 10^{-2}$ & 3 & - \\ 
			LF LSTM & - & 1 $\times$ 98 & 0 $\times$ 0  & Adam & $7.14 \times 10^{-3}$ & 1 & -\\ 
			HF LSTM & - & 2 $\times$ 94 & 0 $\times$ 0  & Adamax & $2.94 \times 10^{-2}$ & 4 & -\\ 
			3-step feed-forward & 1 & - & 4 $\times$ 61 & Adam & 2.54$\times 10^{-2}$ & 13 & -  \\
     & 2 &   -  & 3 $\times$ 41  &  Adamax & $3.89 \times 10^{-3}$  & 67  & -  \\
    & 3 & - &  1 $\times$ 32  &  Adam & $1.78 \times 10^{-4}$  &  151  & - \\ 
			2-step LSTM & 1 & 1 $\times$ 98 & 0$\times$ 0  & Adam & $7.14 \times 10^{-3}$ & 1 & -  \\ 
			     & 2 & 2 $\times$ 98 & 2 $\times 16$ & Adamax & $1.78 \times 10^{-3}$ & 3  & - \\ 
			3-step LSTM & 1 & 1 $\times$ 98 & 0 $\times$ 0  & Adam & $7.14 \times 10^{-3}$ & 1 & - \\
			     & 2 & 4 $\times$ 20 & 0 $\times$ 0 & Adam & $1.80 \times 10^{-2}$ & 3  & - \\
                & 3 & 3 $\times$ 82 & 3 $\times 20$ & Adamax & $6.98 \times 10^{-3}$ & 4  & - \\
            Intermediate & LF & 3 $\times$ 62 & 0 $\times$ 0 & Adam & $3.81 \times 10 ^{-4}$ & 27  & 0.51\\
            & HF & 2 $\times 72$ & 1    $\times$ 118 & " & " & " & "\\
			\hline
		\end{tabular}
		\label{tab: HP_Pulse}
	\end{table}

	\begin{table}[h!]
	\setlength{\tabcolsep}{4.5pt}
 		\caption{Optimized HP values for the lift coefficient in example (II) (Section \ref{section: NS}).}
		\renewcommand{\arraystretch}{1.3}
		\centering%
		\begin{tabular}{c|c|cccccc}
			\hline 
    		Model &  Step & \begin{tabular}{@{}c@{}}Depth $\times$ width \vspace{-4pt} \\ LSTM\end{tabular}   & \begin{tabular}{@{}c@{}}Depth $\times$ width \vspace{-4pt} \\  Dense \end{tabular}   & Optimizer & $\eta$ & Batch size & $\alpha$  \\
			\hline 
			LF feed-forward & - & - & 2 $\times$ 20  & Adamax  & $3.52 \times 10^{-4}$  & 177& -  \\ 
		    HF feed-forward & - & - & 2 $\times$ 34  & Adamax  & $1.49 \times 10^{-4}$ & 149 & -  \\ 
			LF LSTM & - & 4 $\times$ 128 & 1 $\times$ 46  & Adam & $1.96 \times 10^{-3}$ & 1& - \\ 
			HF LSTM & - & 4 $\times$ 24 & 0 $\times$ 0  & Adamax & $5.46 \times 10^{-2}$ & 6& - \\ 
			3-step feed-forward & 1 & - & 2 $\times$ 20  & Adamax  & $3.52 \times 10^{-4}$  & 177& -   \\
     & 2 &   -  & 4 $\times$ 62  &  Adamax & $5.52 \times 10^{-4}$ &    266 & -  \\
    & 3 &  - &   1 $\times$ 44  & Adam  & $8.73 \times 10^{-3}$  &   128 & - \\ 
			2-step LSTM & 1 & 4 $\times$ 128 & 1 $\times$ 46  & Adam & $1.96 \times 10^{-3}$ & 1& - \\
			     & 2 & 3 $\times$ 64 & 2 $\times$ 26 & Adam & $2.66 \times 10^{-3}$ & 5 & -  \\ 
			3-step LSTM & 1 & 4 $\times$ 128 & 1 $\times$ 46  & Adam & $1.96 \times 10^{-3}$ & 1& - \\
			     & 2 & 3 $\times$ 18 & 0 $\times$ 0 & Adamax & $1.17 \times 10^{-2}$ & 6  & - \\
                & 3 & 4 $\times$ 24 & 0 $\times$ 0 & Adamax & $6.20 \times 10^{-4}$ & 8  & - \\ 
            Intermediate & LF & 1 $\times$ 106 & 0 $\times$ 0 & Adam & $2.65 \times 10 ^{-2}$ & 16 & 0.64\\
            & HF & 3 $\times 20$ & 0    $\times$ 0 & " & " & " & "\\
			\hline
		\end{tabular}
		\label{tab: HP_Lift}
	\end{table}

	\begin{table}[h!]
	\setlength{\tabcolsep}{4.5pt}
	\vspace*{-4mm}
 		\caption{Optimized HP values for the drag coefficient in example (II) (Section \ref{section: NS}).}
		\renewcommand{\arraystretch}{1.3}
		\centering%
		\begin{tabular}{c|c|cccccc}
			\hline 
    		Model &  Step & \begin{tabular}{@{}c@{}}Depth $\times$ width \vspace{-4pt} \\ LSTM\end{tabular}   & \begin{tabular}{@{}c@{}}Depth $\times$ width \vspace{-4pt} \\  Dense \end{tabular}  & Optimizer &$ \eta$ & Batch size & $\alpha$ \\
			\hline 
			LF feed-forward & - & - & 4 $\times$ 52   & Adam  & $ 8.91 \times 10^{-2}$ & 162 & - \\ 
			HF feed-forward & - & - & 1 $\times$ 38   &  Adam  &  $1.71 \times 10^{-4}$ &  74  & -  \\ 
			LF LSTM & - & 3 $\times$ 36 & 1 $\times$ 112  & Adam & $1.25 \times 10^{-3}$ & 1 & - \\ 
			HF LSTM & - & 3 $\times$ 36 & 0 $\times$ 0  & Adamax &$ 2.60 \times 10^{-2}$ & 6 & - \\ 
			3-step feed-forward & 1 & - & 4 $\times$ 52   & Adam  & $ 8.91 \times 10^{-2}$ & 162 & - \\
     & 2 &  -  & 2 $\times$ 22   & Adamax  & $5.74 \times 10^{-2}$  &  158  & -  \\
    & 3 & - &  1 $\times$ 39 & Adamax  & $3.32 \times 10^{-3}$  & 107   & -  \\ 
			2-step LSTM & 1 & 3 $\times$ 36 & 1 $\times$ 112  & Adam & $1.25 \times 10^{-3}$ & 1 & - \\
			     & 2 & 3 $\times$ 78 & 1 $\times$ 50 & Adam & $2.87 \times 10^{-3}$ & 10  & -  \\ 
			3-step LSTM & 1 & 3 $\times$ 36 & 1 $\times$ 112  & Adam & $1.25 \times 10^{-3}$ & 1 & - \\
			     & 2 & 3 $\times$ 34 & 0 $\times$ 0 & Adam & $2.06 \times 10^{-3}$ & 10  & -  \\
                & 3 & 1 $\times$ 24 & 0 $\times $ 0 & Adamax & $5.05 \times 10^{-2}$ & 8  & -  \\ 
            Intermediate & LF & 2 $\times$ 46 & 0 $\times$ 0 & Adam & $2.87 \times 10 ^{-2}$ & 11 & 0.1 \\
            & HF & 3 $\times$ 110 & 2    $\times$ 22 & " & " & " & "\\
        \hline
		\end{tabular}
		\label{tab: HP_Drag}
	\end{table}

	\begin{table}[h!]
	\setlength{\tabcolsep}{4.5pt}
	\vspace*{-4mm}
 		\caption{HP values in example (III) (Section \ref{section: LV}). The same structure is used in all the networks, showing that the MF approach is robustly advantageous over the single-fidelity ones, regardless of the hyperparameter tuning.}
		\renewcommand{\arraystretch}{1.3}
		\centering%
		\begin{tabular}{c|c|ccccc}
			\hline     		
			Model &  Step & Depth $\times$ width LSTM & Depth $\times$ width Dense  & Optimizer &$ \eta$ & Batch size \\
			\hline 
			LF LSTM & - & 3 $\times$ 64 & 1 $\times$ 32  & Adamax & $1.00 \times 10^{-3} $& 50\\ 
			HF LSTM & - & 3 $\times$ 64 & 1 $\times$ 32  & Adamax & $1.00 \times 10^{-3}$ & 50\\ 
			2-step LSTM & 1 & 3 $\times$ 64 & 1 $\times$ 32  & Adamax & $1.00 \times 10^{-3} $& 50 \\
			     & 2 & 3 $\times$ 64 & 1 $\times$ 32 & Adamax & $5.00 \times 10^{-3}$ & 100  \\
			\hline
		\end{tabular}
		\label{tab: HP_LV}
	\end{table}

\clearpage
\bibliographystyle{abbrv}
\bibliography{references.bib}

\end{document}